\documentclass[a4paper,11pt]{amsart}
\usepackage{amssymb}
\usepackage{xcolor}
\usepackage{hyperref}
\usepackage{enumitem}

\usepackage[margin=2.9cm]{geometry}

\newtheorem{Main}{Theorem}

\newtheorem{Mainprop}[Main]{Proposition}

\usepackage{color}
\usepackage{caption}
\usepackage{subcaption}

\usepackage{xcolor}
\usepackage[T1]{fontenc}
\usepackage[draft]{todonotes}

\definecolor{orange}{rgb}{1,0.5,0}

\usepackage[backend=biber, maxnames=4, style=alphabetic,
            maxalphanames=3, minalphanames=1,
]{biblatex}
\DeclareLabelalphaTemplate{%
  \labelelement{%
    \field[strwidth=3,strside=left,ifnames=1]{labelname}%
    \field[strwidth=1,strside=left]{labelname}%
  }%
}
\DeclareFieldFormat{extraalpha}{%
  \ifnumless{\thefield{year}}{2000}%
    {\the\numexpr\thefield{year}-1900}%
    {\the\numexpr\thefield{year}-2000}}

\usepackage{sidecap}

\usepackage{comment}

\usepackage{mathrsfs}
\DeclareMathAlphabet{\mathpzc}{OT1}{pzc}{L}{it}

\theoremstyle{definition}
\newtheorem{definition}{Definition}[section]
\newtheorem{theorem}[definition]{Theorem}

\newtheorem{proposition}[definition]{Proposition}

\newtheorem{lemma}[definition]{Lemma}
\newtheorem{example}[definition]{Example}

\newtheorem{remark}[definition]{Remark}

\newtheorem{claim}{Claim}
\newtheorem*{claim*}{Claim}

\numberwithin{equation}{section}
\usepackage{thmtools, thm-restate}

\def\dist{\mathrm{dist}}

\def\cP{\mathcal{P}}
\def\geq{\geqslant}
\def\leq{\leqslant}
\def\R{\mathbb{R}}
\def\T{\mathbb{T}}
\def\P{\mathbb{P}}

\newcommand{\Sc}{{\mathcal S}}
\def\E{\mathbb{E}}
\def\Z{\mathbb{Z}}
\def\N{\mathbb{N}}
\def\F{\mathbb{F}}

\def\Q{\mathbb{Q}}
\def\cF{\mathcal F}
\newcommand{\be}{\begin{equation}}
\newcommand{\ee}{\end{equation}}
\newcommand{\RQ}{\R \setminus \Q}
\newcommand{\eps}{\varepsilon}

\newcommand{\rv}{\mathcal}
\newcommand{\norm}[1]{\|{#1}\|}
\newcommand{\snorm}[1]{(({{#1}))}}
\newcommand{\ssnorm}[1]{\langle\langle{{#1}\rangle\rangle}}
\newcommand{\floor}[1]{\left\lfloor #1 \right\rfloor}

\newcommand{{\M}}{\mathcal{M}}
\newcommand{\cU}{\mathcal{U}}
\newcommand{\ph}{\varphi}

\title[A limit theorem for deviation of ergodic integrals on flows]{A limit theorem for deviation of ergodic integrals on a class of flows on the two-torus}
\author[M Lorenzo--Laguno]{Miriam Lorenzo--Laguno}
\date{}
\begin{document}

\begin{abstract} We study the fluctuations of ergodic integrals for  
{locally Hamiltonian flows on the two-torus with a single non-degenerate saddle (Arnol'd flows)}.
We sample a flow from a suitably random family, together with a base point $p$ chosen uniformly at random on that flow's domain. For a smooth observable $f$, we show that the centered Birkhoff integral $ \int_{0}^t f(\phi_s p)\,ds-t\bar f$, where $\bar f$ is the spatial mean of $f$, normalized by $(\log t)^2$, converges in distribution to the product of two independent random variables. 
Using special representation of locally Hamiltonian flows, we reduce the problem to studying ergodic sums of functions with logarithmic singularities over circle rotations.
Then, we prove a Cauchy limit theorem for such ergodic sums, which is a result of independent interest.
\end{abstract}

\maketitle

\section{Introduction and main results}
\subsection{Arnol'd flows with a single saddle}
Smooth area-preserving flows on surfaces are an extensively studied class of dynamical systems. They are low-dimensional examples of systems with non-trivial ergodic, mixing, and spectral properties.

Let $\mathcal{M}'$ be a compact, connected, oriented smooth surface equipped with a smooth non-degenerate area form $\omega$. Any smooth closed 1-form $\lambda$ on $\mathcal{M}'$ determines a unique smooth vector field $W: \mathcal{M}' \to T\mathcal{M}'$, given by $\iota_W(\omega)=\lambda$, where $\iota$ denotes the contraction operator. The vector field $W$ (or equivalently $\lambda$) generates an $\omega$-preserving flow $\phi_\R=(\phi_t)_{t \in \mathbb{R}}$. These are known as \emph{locally Hamiltonian flows}. Indeed, locally there exists a smooth real-valued function $H$ such that $\lambda=dH$. In local coordinates $(x,y)$ for which $\omega = dx \wedge dy$, the vector field reads $W=\dot x\partial_ x+\dot y\partial_y$ where $\dot{x} = \frac{\partial H}{\partial y},\dot{y} = -\frac{\partial H}{\partial x}$.

\begin{figure}[]
    \centering
    \includegraphics[width=0.55\linewidth]{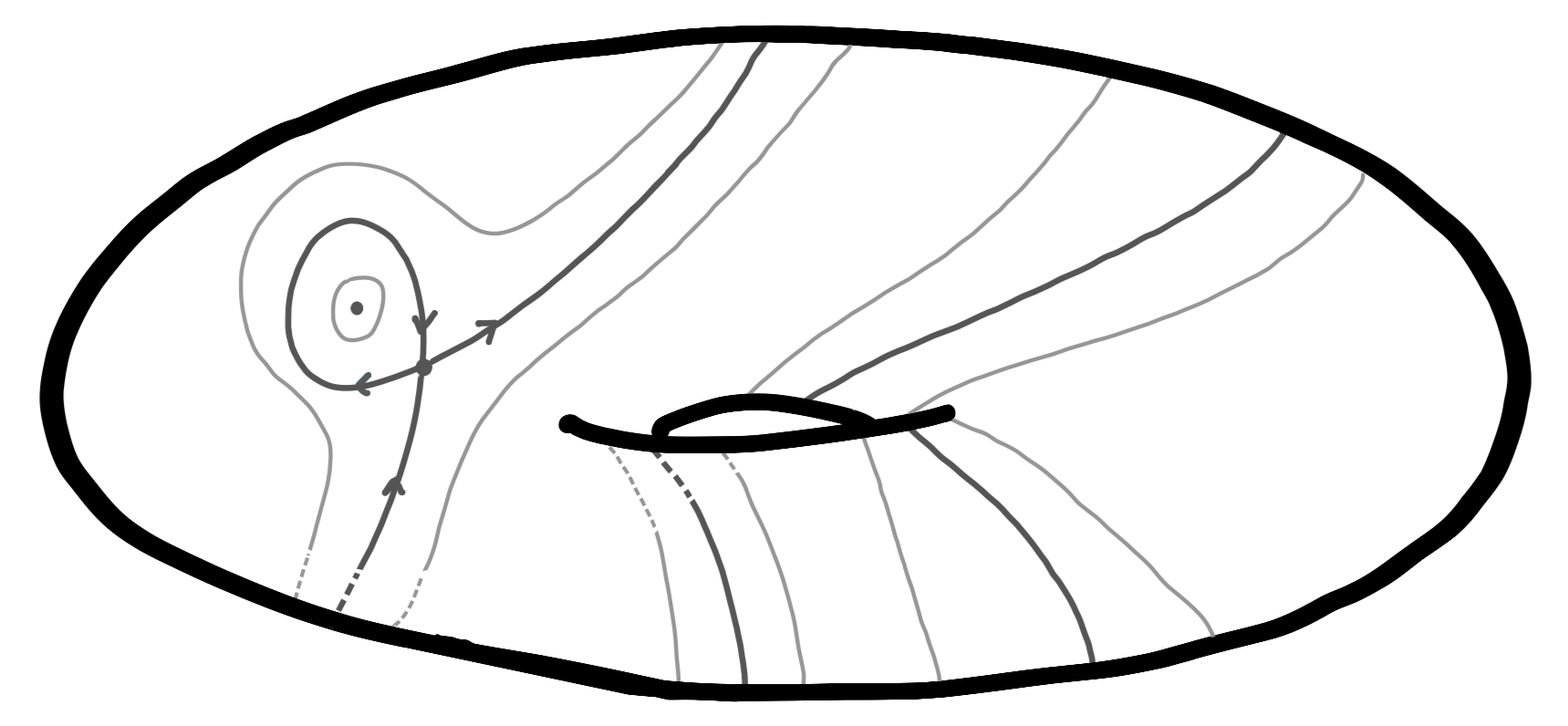}
    \caption{An Arnold' flow in $\T^2$.}
    \label{fig:Af}
\end{figure}

\textit{Arnold flows}, first studied by Arnold in \cite{Arnold1991}, are among the most fundamental examples of locally Hamiltonian flows. In this setting, the underlying manifold is the two-torus $\mathbb{T}^2$, and we assume that $W$ has exactly two zeros: a non-degenerate ($\det DW\neq 0$) saddle and a center lying inside a homoclinic loop homologous to zero (see Figure \ref{fig:Af}). This loop divides $\mathbb{T}^2$ into two $\phi_\R$-invariant components: the one bounded by the saddle loop is foliated by periodic orbits, and we denote the closure of its complement by $\M$. We refer to the flow $(\phi_t\vert_{\M})_{t\in\R}$ as an Arnol'd flow, and we also denote it by $\phi_\R$.

Arnol'd flows are \textit{mixing }for every \textit{irrational frequency} (see Section \ref{sec:flowalpha}). This was first shown by Sinai--Khanin \cite{SinaiKhanin1992} for Lebesgue-almost every irrational frequency and was later improved by Kochergin \cite{Kochergin2003,Kochergin2004} to all irrational frequencies.

In particular, for such flows we have \textit{equidistribution of typical trajectories}: for any integrable function $f$, the Birkhoff integrals satisfy $\frac{1}{t}\int_{0}^t f(\phi_sp)ds\to\frac{1}{\omega(\M)}\int_{\M}fd\omega$ pointwise almost-surely as $t\to\infty$, where $\omega(\M)$ is the area of $\M$. This raises the question of the nature of its fluctuations: how large are the \textit{deviations of the ergodic averages} and how are they asymptotically distributed? 

More precisely, we study the centered Birkhoff integral
\be\label{eq:ergodicavg}
\int_{0}^t f(\phi_s p)\,ds-\frac{t}{\omega(\M)}\int_{\M} f\,d\omega,
\ee
as $t\to\infty$. We show that, provided $f$ is smooth enough and the initial point $p$ is chosen uniformly in $\M$ (with respect to $\omega$), the appropriate scaling of the fluctuations is $(\log t)^2$, for suitable random families of flows (made precise in Definition \ref{def:good}). After this normalization, \eqref{eq:ergodicavg} converges in distribution to the product of a Cauchy random variable and a random variable quantifying the difference between the spatial average of $f$ and its value at the saddle. 

As will be shown, this question is tightly related to limit theorems for ergodic sums over rotations, 
in particular, for functions with (asymmetric) logarithmic singularities.

\subsection{Supension flows}\label{sec:susflow}
Let $T$ be a measure-preserving endomorphism on an interval $[a,b)$ and let $\varphi:[a,b)\to\R^{+}\cup\{\infty\}$, satisfying $\inf_{x \in[a,b)} \varphi(x) > 0$. The \textit{suspension flow} $(T_t)_{t\in\mathbb{R}}$ over $T$ built under the roof function $\ph$ is a flow defined on the space
\[
X = \left\{(x,s)\in[a,b)\times\mathbb{R} : 
0\leq s<\varphi(x)\right\}
\]
with the identification $(x,\varphi(x))\sim (Tx,0)$ at the points where $\varphi(x)<\infty$, and acts by vertical translation
\[
T_t(x,s)=(x,s+t)
\]
followed by the canonical identifications whenever the second coordinate exits the interval $[0,\varphi(x))$.

Let $\T$ denote the torus, which we identify with $[0,1)$. For every $\alpha\in\R$, the rotation map $R_\alpha:\T\to\T$ is defined by $R_\alpha x=x+\alpha\bmod1$ for every $x\in\T$.

\subsection{Arnol'd flows as suspension flows}\label{sec:flowalpha} Consider $\T^2$ equipped with a smooth non-degenerate area form $\omega$. Let $\phi_\R$ be an Arnol'd flow on ${\M }\subset\T^2$ generated by a smooth vector field $W$.
It is a classical result (see \cite{Ravotti2017,KOZ2026} and the references therein) that $\phi_\R$ is measure-theoretically isomorphic to the suspension flow $(T_t)$ over a rotation $R_\alpha$ on $\T = [0,1)$, with roof function $\ph$ with logarithmic singularities and $\inf_\T \ph >0$; via a map $\Phi \colon X \to \M$ that is measure-preserving up to a linear rescaling of the $x$-coordinate\footnote{That is, $\Phi^{-1}$ is the composition of a measure-preserving map onto a suspension space over $[0,c)$, $c>0$, with a linear rescaling of the base coordinate $x\mapsto x/c$. In particular, $\int_\T\ph(x)dx=\omega(\M)/c$.}. 

We call $\Sigma :=\Phi (\T\times\{0\})\subset{\M }$ the base transversal section. Then, $\Phi(x,\cdot)$ denotes the orbit under $\phi_\R$ of $\Phi(x,0)\in\Sigma$, and $\ph(x)$ is the return time of $\Phi(x,0)$ to $\Sigma$.

Let $p_0$ be the (non-degenerate) saddle of $\phi_\R$. The rotation parameter $\alpha$ can be computed directly from the $1$-form $\lambda$ corresponding to the flow. Indeed, let $z_0$ and $z_0'$ denote the first intersections of the stable and unstable separatrixes of the saddle $p_0$ with $\Sigma$, respectively. Then the rotation parameter is given by (see Section \ref{sec:passtime})
\begin{equation}
\label{eq:VF-RotNumb}    
\alpha = \left(\int_{\Sigma}\lambda\right)^{-1}
\int_{z_0}^{z_0'}\lambda.
\end{equation}

\subsection{Roof function and passage constant}\label{sec:passtime0} Here, we present some formulas that will be useful for the statement of the main results. We refer the reader to Section \ref{sec:passtime} for further exposition and proofs.

Let $T(c)$ be the passage time around a fixed neighborhood of the saddle $p_0$ at Hamiltonian energy level $c$--close. Define the passage constant  (see Section \ref{sec:passtime})
\be\label{eq:passtime}
    c_0  = \lim_{c\to 0}\frac{T(c)}{|\log c|}=\frac{1}{\sqrt{-\text{det}DW(p_0)}}.
\ee
Note that $c_0$ depends only on $\det DW(p_0)$.

The parameter in \eqref{eq:passtime} determines the return time to $\Sigma$ of $\Phi(x,0)$ up to a bounded variation error. Indeed, the return time function is given by 
\begin{equation}
        \label{eq:def_phi}
        \ph (x) = \psi (x)+ 2c_0|\log x|+c_0|\log(1-x)|,
    \end{equation}
where $\psi $ is of bounded variation, and the choice between $c_0$ and $2c_0$ depends on whether the orbit passes once or twice near $p_0$ before returning to the Poincaré section $\Sigma$, the latter corresponding to the side of the loop homologous to zero (see Figure \ref{fig:Af}).

\begin{remark} Given that a vector field $W$ generates an Arnold flow on $\T^2$, all $\phi_t$, $p_0$, $\M$, $\alpha$, $\ph$, $\psi$ $c_0$, $T_t$, $X$, $\Phi$, and $\Sigma$ may depend on $W$. We denote this dependence by $\cdot^W$.\end{remark}

\subsection{Main result}
Let $\cF$ be a family of vector fields on $\T^2$ and let $\hat\eta$ be a probability measure on $\cF$. 
\begin{definition} \label{def:good} We say that $(\cF,\hat\eta)$ is \textit{good} if 
\begin{itemize}
    \item[a)] Every $W\in\cF$ generates an Arnol'd flow on $\T^2$.
    \item[b)] Under the identification via $\Phi^W$ of $(\M^W,\phi_t^W)$ and $(X^W,T_t^W)$ for each $W$ described in Section \ref{sec:flowalpha}, let $\alpha(W)$ be the rotation number of the base map $R_\alpha$ (given by \eqref{eq:VF-RotNumb}). Then as $\rv W\sim\hat\eta$, $\rv A=\alpha(\rv W)$ has an absolutely continuous distribution with respect to Lebesgue$(\T)$.
\end{itemize}

\end{definition}

Let $\epsilon>0$ and $f\in C^{2+\epsilon}(\T^2)$. For each $W\in \cF$, define the mean of $f$ in $\M^W$ as
    \be\label{eq:meanf}\bar f^W:=\frac{1}{\omega(\M^W)}\int_{\M^W}f(p)\,dp.\ee
We want to study the scale and the limiting distribution as $t\to\infty$ of the deviations from the mean of the Birkhoff integral
\[
\int_0^{t} f \left(\phi^W_sp\right)ds,
\]
under $\rv W\sim\hat\eta$ and uniform $p\in{\M^W}$ (with respect to $\omega$). Our main result, Theorem \ref{thm:main}, states that the right scale is $(\log t)^2$ and the limiting distribution is the product of the random variable $c_0^W(f(p_0^{\rv W})-\bar f^{\rv W})$ and an independent Cauchy distribution, provided $(\cF,\hat\eta)$ is \textit{good}.

\begin{Main}\label{thm:main}
Let $\cF$ be a family of vector fields on $\T^2$, let $\hat\eta$ be a probability measure on $\cF$, let $\epsilon>0$, and let $f\in C^{2+\epsilon}(\T^2)$. Suppose $(\cF,\hat\eta)$ is \textit{good} in the sense of Definition \ref{def:good}.
Then, taking $\rv W\sim\hat \eta$ in $\cF$, and then $\rv P$ uniform with respect to $\omega$ in each $\M^W$, we have
\[
\frac{\int_0^{t} f \left(\phi^{\rv W}_s\rv P\right)ds - {t}\bar f^{\rv W}}{(\log t)^2}\ \xrightarrow[t\to\infty]{\text{dist}}\ \left(f(p_0^{\rv W})-\bar f^{\rv W}\right)\frac{c_0^{\rv W}}{4\pi}\,\mathfrak{C},
\]
where $\mathfrak{C}$ is the standard Cauchy distribution, independent of $\rv W$, and $\bar f^W$ is as in \eqref{eq:meanf}.

Equivalently, for every $r\in\R$,
\[
\lim_{t\to\infty}\P\left({4\pi}\frac{\int_0^{t} f \left(\phi^{\rv W}_s\rv P\right)ds - {t}\bar f^{\rv W}}{(\log t)^2}\leq r\right)
=\int_{-\infty}^{\infty}\left(\frac12+\frac{1}{\pi}\arctan\!\left(\frac{r}{|a|}\right)\right)\,d\nu(a),
\]
where $\nu$ denotes the law of the random variable $c_0^{\rv W}(f(p_0^{\rv W})-\bar f^{\rv W})$ induced by $\rv W\sim\hat\eta$ on $\R$, and, by convention, $\arctan(\pm\infty)=\pm\pi/2$, so that the integrand equals $\chi_{[0,\infty)}(r)$ at $a=0$.
\end{Main}

Therefore, the fluctuations of the Birkhoff integral of a function depend essentially on the value at the saddle relative to the mean, and on the time spent near the saddle relative to the time spent far from it, i.e., on the passage constant $c_0$ in \eqref{eq:passtime}.

Theorem \ref{thm:main} requires $(\cF,\hat\eta)$ to be \textit{good} in the sense of Definition \ref{def:good}. The main reason is that we leverage some Diophantine conditions of $\alpha$ (c.f. Section \ref{sec:proof1}) and some limit theorems of ergodic theorems over rotations (see Proposition \ref{prop:main_harm} and Lemma \ref{claim:boundsEsums} below, and Sections \ref{sec:knownres} and \ref{sec:proof2}) that fail to hold when $\alpha(\rv W)$ does not have an absolutely continuous distribution with respect to Lebesgue on $\T$ when $\rv W\sim\hat\eta$. 
Nevertheless, there are several examples satisfying the assumptions in Definition \ref{def:good}. We present some in Section \ref{sec:ex}.

As shown in Section \ref{sec:mainproof}, Theorem \ref{thm:main} is a consequence of the following two results.
\begin{Main}\label{thm:mainbis}
Let $\cF$ be a family of vector fields on $\T^2$, let $\hat\eta$ be a probability measure on $\cF$, let $\epsilon>0$, and let $f\in C^{2+\epsilon}(\T^2)$. Suppose $(\cF,\hat\eta)$ is \textit{good} in the sense of Definition \ref{def:good}. Moreover, suppose $f(p_0^{W})-\bar f^W\neq 0$, for $\bar f^W$ as in \eqref{eq:meanf}, holds $\hat\eta$-almost surely.
Then, taking $\rv W\sim\hat \eta$ in $\cF$, and then $\rv P$ uniform in each $\M^W$, we have
\[
\frac{\int_0^{t} f \left(\phi^{\rv W}_s\rv P\right)ds - {t}\bar f^{\rv W}}{{c_0^{\rv W}}\left(f(p_0^{\rv W})-\bar f^{\rv W}\right)(\log t)^2}\ \xrightarrow[t\to\infty]{\text{dist}}\ \frac{1}{4\pi}\,\mathfrak{C},
\]
where $\mathfrak{C}$ is the standard Cauchy distribution.
\end{Main}
\begin{Main}\label{thm:mainbisbis}
Let $\cF$ be a family of vector fields on $\T^2$, let $\hat\eta$ be a probability measure on $\cF$, let $\epsilon>0$, and let $f\in C^{2+\epsilon}(\T^2)$. Suppose $(\cF,\hat\eta)$ is \textit{good} in the sense of Definition \ref{def:good}. Moreover, suppose ${c_0^{\rv W}}(f(p_0^{W})-\bar f^W)\equiv a\in\R$ is constant $\hat\eta$-almost surely.
Then, taking $\rv W\sim\hat \eta$ in $\cF$, and then $\rv P$ uniform in each $\M^W$, we have
\[
\frac{\int_0^{t} f \left(\phi^{\rv W}_s\rv P\right)ds - {t}\bar f^{\rv W}}{(\log t)^2}\ \xrightarrow[t\to\infty]{\text{dist}}\ \frac{a}{4\pi}\,\mathfrak{C},
\]
where $\mathfrak{C}$ is the standard Cauchy distribution.
\end{Main}
For $a\in\R^*$, Theorem \ref{thm:mainbisbis} is a direct consequence of Theorem \ref{thm:mainbis}, but it also covers the case $a=0$. The proof of Theorem \ref{thm:mainbisbis} follows the same steps of that of Theorem \ref{thm:mainbis}, with the extra simplification that ${c_0^{\rv W}}(f(p_0^{\rv W})-\bar f^{\rv W})$ is no longer a random variable, so we do not need to normalize by it the left hand side of the limit. We prove Theorem \ref{thm:mainbis} and leave a sketch of the elementary modifications required to prove Theorem \ref{thm:mainbisbis} to Appendix \ref{app:1}.

\subsection{Birkhoff sums and logarithmic profiles} Theorem \ref{thm:mainbis}, hence Theorem \ref{thm:main}, relies on a limit theorem for Birkhoff sums of asymmetric logarithmic profiles: Proposition \ref{prop:main_harm}. 

Define, for $x\in(0,1)$,  
\be
\label{eq:sing_p}
h(x)= 2|\log x|+|\log(1-x)|-3.
\ee
Given $N\geq 0$, $\alpha\in\R$, $x\in \T$, and $g\in L^2(\T)$, define the Birkhoff sum of $g$
\[
    S(\alpha,g,x,N):=\sum_{m=0}^{N-1}g(R_\alpha^mx).
\]

We say that a family of scalar functions $(\rho_t)_{t\in I}$, $I$ any set, is uniformly $C^2(\T)$ if $\rho_t\in C^2(\T)$ for every $t\in I$ and $\sup_I\norm{\rho_t}_{C^2}<\infty$, where $\norm{\cdot}_{C^2}=\norm{\cdot}_{\infty}+\norm{D\cdot}_{\infty}+\norm{D^2\cdot}_{\infty}$ denotes the standard $C^2$ norm.

In the spirit of the main theorems in \cite{Kesten1962,DolgopyatFayad2020, DolgopyatSarig2019} (see Section \ref{sec:knownres}), we show the following.

\begin{Mainprop}\label{prop:main_harm}
Let $\rv A$ be distributed according to $\eta\in C^\infty(\T)$ and $\rv X|\rv A=\alpha\sim \rho_{\alpha}$, $\eta-$almost surely, with $(\rho_{\alpha})$ uniformly $C^2(\T)$. 
Then, for $h$ as in \eqref{eq:sing_p},
\begin{equation}\label{eq:convalpha}
\frac{S(\rv A,h,\rv X,N)}{(\log N)^2}\ \xrightarrow[N\to\infty]{\text{dist}}\ \frac{1}{4\pi} \mathfrak{C}
\end{equation}
where $\mathfrak{C}$ is the standard Cauchy distribution.
\end{Mainprop}

As it turns out, the randomization of both $\alpha$ and $x$ is crucial for the result. We make this explicit in the result below, which follows from classical results for the continued fraction expansion of irrational numbers and Khinchin's theorem. We include the proof in Appendix {\ref{sec:boundsEsums}}. Ultimately, the asymmetry of $h$ is what gives raise to the $(\log N)^2$ scaling.

\begin{lemma}
\label{claim:boundsEsums}
    For Lebesgue-almost every $\alpha\in(0,1)$, there exist an increasing sequence $N_k$ and a constant $C>0$ such that, for a Lebesgue-full measure set of $x$,
    \be \label{eq:small}
    \limsup_{k\to\infty}\frac{|S(\alpha,h,x,N_k)|}{\log N_k}
    \leq C.
    \ee
    Moreover, let $\psi(n)=\log n\log\log n\log\log \log n$. For Lebesgue-almost every $(\alpha,x)\in\T^2$, there exist an increasing sequence $N_k$ and $C'>0$ such that 
    \be \label{eq:large}
    \liminf_{k\to\infty}\frac{|S(\alpha,h,x,N_k)|}{\log N_k\psi(N_k)}
    \geq C'.
    \ee
\end{lemma}
In particular, we get
\[
\text{Leb}\left(\left\{(\alpha,x):\limsup_{N\to\infty}\frac{|S(\alpha,h,x,N)|}{(\log N)^2\log\log N}=\infty,\ \liminf_{N\to\infty}\frac{|S(\alpha,h,x,N)|\log\log N}{(\log N)^2}=0\right\}\right)=1.
\]

\subsection{Known results}\label{sec:knownres}
\subsubsection{On the growth of ergodic integrals of flows on surfaces}
The first results on the deviations of ergodic averages for smooth flows on surfaces are due to Forni \cite{Forni2002}. {However, his results apply only to mean zero observables that, crucially, vanish at the fixed points of the flow.}

Fr{\k{a}}czek and Ulcigrai subsequently lifted this restriction, extending this theory to observables that may be nonzero at the fixed points. Their result is stated in greater generality for locally Hamiltonian flows on surfaces of arbitrary genus without degenerate saddles in \cite{FraczekUlcigrai2024}. Below, we state the corresponding result for $\T^2$. Finally, Fr{\k{a}}czek and Kim \cite{FraczekKim2024} further extended the deviation theory to locally Hamiltonian flows with degenerate saddles.

\begin{theorem}\cite{FraczekUlcigrai2024}
For a full measure set (see Section \ref{sec:Kfc}) of Arnol'd flows on $\M\subset \T^2$ with saddle $p_0$, for any $\epsilon>0$, and for every $f\in C^{2+\epsilon}(\M)$, there exist an smooth cocycle $u:\mathbb{R}\times \M\to\mathbb{R}$ such that we have the asymptotic expansion
\[
\int_0^t f(\phi_s(p))\,ds-t\bar f=\left(f(p_0)-\bar f\right)u(t,p) + \operatorname{err}(f,t,p),\quad \operatorname{err}(f,t,p)=\int_0^t f_e(\phi_s(p))\,ds,
\]
where $\bar f$ is the mean of $f$ and $f_e\in C^{2+\epsilon}(\M)$ is a mean zero function with $f_e(p_0)=0$; and $u$, $\operatorname{err}$ are such that
\[
\limsup_{t\to+\infty}\frac{\log|u(t,p)|}{\log t}=0, \quad
\limsup_{t\to+\infty} \frac{\log\|u(t,\cdot)\|_{L^q(\M)}}{\log t} =0, \quad
\text{for a.e. }p\in \M\text{ and all }q\geq 1,
\]
\[
\limsup_{t\to+\infty}\frac{\log|\operatorname{err}(f,t,p)|}{\log t}\leq 0, \quad
\text{for a.e. }p\in \M.
\]
\end{theorem}

{Therefore, our results give a distributional description of the cocycle $u^W(t,p)$. Indeed, letting $\cF$ be a family of vector fields on $\T^2$ and $\hat\eta$ be a probability measure on $\cF$ such that $(\cF,\hat\eta)$ is \textit{good} in the sense of Definition \ref{def:good}, Theorems \ref{thm:mainbis} and \ref{thm:mainbisbis} yield
\[
\frac{u^{\rv W}(t,\rv P)}{c_0^{\rv W}(\log t)^2}\xrightarrow[t\to\infty]{\dist}\frac{1}{4\pi}\mathfrak{C}, \qquad \frac{\operatorname{err}^{\rv W}(f,t,\rv P)}{(\log t)^2}\xrightarrow[t\to\infty]{\dist}0,
\]
where the flow $\rv W$ is sampled on $\cF$ according to $\hat\eta$ and $\rv P$ is sampled randomly on $\M^W$ with respect to the area form $\omega$. Recall $\mathfrak{C}$ is the standard Cauchy distribution.
In particular, if $c_0^W=c_0\in\R^+$ is constant on $\cF$ (see Section \ref{sec:ex} for examples), we have 
\[
\frac{u^{\rv W}(t,\rv P)}{(\log t)^2}\xrightarrow[t\to\infty]{\dist}\frac{c_0}{4\pi}\mathfrak{C}.
\]}

\subsubsection{Limit theorems for ergodic sums over rotations}
As discussed in Section \ref{sec:flowalpha}, Arnold flows can be studied as suspension flows over a rotation $R_\alpha$ and under a roof function exhibiting logarithmic singularities. Consequently, as established in Sections \ref{sec:proof0} and \ref{sec:proof1}, the ergodic integrals of a function $f$ are closely related to ergodic sums of an associated function $\hat f$ over $R_\alpha$, also with logarithmic singularities. Our results therefore draw on several works on limit theorems for Birkhoff sums over circle rotations.

For any $g\in L^1(\T)$, define \[ S(\alpha,g,x,N):=\sum_{n=0}^{N-1}g(R_\alpha^nx). \] 
In a pioneering work, Kesten studied the distribution of $S(\alpha,g,x,N)$ for some functions with jump discontinuities, when $(x,\alpha)$ are uniformly distributed on $\T^2$. In \cite{Kesten1960,Kesten1962}, he showed that $S(\alpha,g,x,N)/\log N$ converges in distribution to a Cauchy distribution, for \[g(x)=\{x\}-1/2,\] where $\{x\}$ denotes the fractional part of $x\in\R$, and for \[g(x)=\chi_{[0,\beta)}(x)-\beta,\qquad \beta\in(0,1).\]

Later, Dolgopyat and Fayad \cite{DolgopyatFayad2014,DolgopyatFayad2020} obtained an analogous to Kesten's Cauchy limit for the discrepancy of indicator functions of boxes on higher-dimensional toral translations, identifying the correct scaling function as $(\log N)^d$, $d$ being the dimension of the torus, and the limit as a Cauchy distribution. The main mechanism of their results is a Poisson limit theorem on the space of unimodular lattices. As it will be addressed in Section \ref{sec:proof2}, it is also key for our results.

The next natural step is to study ergodic sums over $R_\alpha$ of functions with unbounded singularities. Sinai and Ulcigrai \cite{SinaiUlcigrai2008} proved convergence in distribution of $S(\alpha,g,x,N)/N$ for \[g(x)=(1-e^{2\pi i x})^{-1},\] to an unspecified distribution function, again when $(x,\alpha)$ is uniformly distributed in $\T^2$. Note that $g$ has a singularity of order $1/x$ at $x=0$.

These results illustrate how the nature of the singularities of $g$ may determine its normalization scale. In our setting, the relevant observable has an asymmetric logarithmic singularity (see \eqref{eq:def_phi} and \eqref{eq:sing_p})
\[h(x)=2|\log x|+|\log (1-x)|-3.\]
The limit theorem underlying Theorems \ref{thm:main} and \ref{thm:mainbis} is Proposition \ref{prop:main_harm}, which states that\[S(\alpha,h,x,N)/(\log N)^2,\] converges in distribution to a Cauchy distribution when $x$ and $\alpha$ are random with sufficiently regular distributions.

Logarithmically singular Birkhoff sums have also been studied for fixed rotation numbers; see, for example, Knill and Tangerman \cite{KnillTangerman2011} and Aistleitner, Technau and Zafeiropoulos \cite{AistleitnerTechnauZafeiropoulos2023}. In particular, the latter study Sudler products:
\[
P_N(\alpha)=\prod_{k=1}^N2|\sin(\pi k\alpha)|,
\]
whose logarithms are Birkhoff sums of $g(x)=\log |2\sin(\pi x)|$ on the $R_\alpha$-orbit of $\alpha\in\T$. For some choices of $\alpha$, they show 
\[
\limsup _{N\to\infty}P_N(\alpha)/N<\infty.
\]

Other authors have explored settings in which one or both of $\alpha$ and $x$ are fixed. For example, one can average with respect to the time variable, that is, taking $n$ uniformly on $\{1,...,N\}$, while fixing one or both of $\alpha$ and $x$. Results in this direction include those of Beck \cite{Beck1994}, Borda, Fr{\"u}hwirth and Hauke \cite{Borda2023}, and Dolgopyat and Sarig \cite{DolgopyatSarig2019}. The results in Section \ref{sec:proof2} (the proof of Proposition \ref{prop:main_harm}) are particularly inspired on the work of Dolgopyat and Sarig. 

\subsection*{Plan of the paper.}
The remainder of the paper is structured as follows. In Section \ref{sec:Af}, before turning to the proofs of the main theorems, we provide further motivation for the results by discussing Arnold flows in greater detail, proving \eqref{eq:VF-RotNumb} and \eqref{eq:passtime}, and presenting families of examples satisfying Definition \ref{def:good}.

In Section \ref{sec:mainproof} we prove Theorem \ref{thm:main} by showing it follows from Theorems \ref{thm:mainbis} and \ref{thm:mainbisbis}.

The proof of Theorem \ref{thm:mainbis} relies on a cascading reduction that translates a continuous-time physical flow problem on the two-torus into a discrete-time arithmetic and probabilistic problem on a circle rotation. In Section \ref{sec:proof0} we decompose Theorem \ref{thm:mainbis} into three self-contained propositions: Propositions \ref{prop:main_dist}, \ref{prop:main_red}, and \ref{prop:main_harm}. We then show how these propositions cleanly combine to yield the final Cauchy limit law. The proofs of the three propositions are then carried out in Sections \ref{sec:proof3}, \ref{sec:proof1}, and \ref{sec:proof2}, respectively.

\subsection*{Notation}

\begin{itemize}[leftmargin=*]
    \item Throughout this paper, we identify $\mathbb{T}$ with $[0,1)$.
    
    \item We denote $\mathcal{X}$ the random variable taking values $x$, $\mathcal{W}$ the random variable taking values $W$, $\rv P$ the random variable taking values $p$, and $\mathcal{A}$ the random variable taking values $\alpha$.

    \item Expectations and probabilities are denoted by $\E$ and $\P$, respectively. Whenever several sources of randomness are present, subscripts indicate the corresponding probability measure. For instance, $\P_{\rv  X}(\cdot\,|\,\rv W)$ or $\E_{\rv W}(\cdot)$.
    
    \item Given a random variable $X$, we write
    \begin{itemize}
        \item[] $X\sim\mu$ to indicate that $X$ has distribution $\mu$,
        \item[] $X\,|\,Y=y\sim\mu_y$ meaning that $X$ is sampled conditionally on another random variable $Y$ with distribution $\mu_y$ for $\P_Y-$almost every $y$.
    \end{itemize}
    
    \item For a sequence of random variables $(X_n)$ and a random variable $X$, we write
    \begin{itemize}
        \item[] $X_n \xrightarrow{} X$ to denote almost-surely point-wise convergence,
        \item[] $X_n \xrightarrow{\mathrm{dist}} X$ to denote convergence in distribution,
        \item[] $X_n \xrightarrow{\P} X$ to denote convergence in probability.
    \end{itemize}

    \item For two sequences $a_n$ and $b_n$, we write
    \begin{itemize}
        \item[] $a_n = O(b_n)$ or $a_n\lesssim b_n$ (when $a_n\geq 0$) if there exists a constant $C>0$, independent of the asymptotic parameter, such that $|a_n|\leq C b_n$ for all sufficiently large $n$,
        \item[] $a_n=o(b_n)$ if $\lim_{n\to\infty}\frac{a_n}{b_n}=0$,
        \item[] $a_n\asymp b_n$ if both $a_n\lesssim b_n$ and $b_n\lesssim a_n$ hold.
    \end{itemize}

    \item When random variables are involved, we use the notation
    \begin{itemize}
        \item[] $X_n=O_{\P}(a_n)$ to mean that the family $(X_n/a_n)_n$ is tight,
        \item[] $X_n=o_{\P}(a_n)$ to mean that $\frac{X_n}{a_n}\xrightarrow{\P}0$.
    \end{itemize}
    More generally, we write $O_{\P_{\rv W}}(\cdot)$, etc., to specify the corresponding probability space.

    \item We denote by $\chi_E(x)$ the indicator function of $x\in E$.

    \item Given $x\in \R$, we write $\snorm x$ to denote the signed distance to the nearest integer, that is, $x=\snorm{x}+m$, where $\snorm{x}\in(-\frac{1}{2},\frac{1}{2}]$ and $m\in\Z$, and write $\norm{x}:=\dist(x,\Z)=|\snorm{x}|$ and $\ssnorm{x}=(-1)^m\snorm{x}$. We denote by $\lfloor \cdot\rfloor$, $\lceil \cdot\rceil$ the unique integers such that $\lfloor x\rfloor\leq x\leq\lceil x\rceil$.

    \item Given a function $f:E\to \R$, we denote by $\norm{f}_\infty=\sup_E|f|$, by $\norm{\cdot}_{C^2}=\norm{\cdot}_{\infty}+\norm{D\cdot}_{\infty}+\norm{D^2\cdot}_{\infty}$ the standard $C^2$ norm, and by $\text{Var(f)}$ its variation.

    \item We say that a family of scalar functions $(\rho_t)_{t\in E}$, $E$ any set, is uniformly $C^2(\T)$ if $\rho_t\in C^2(\T)$ for every $t\in E$ and $\sup_E\norm{\rho_t}_{C^2}<\infty$.
\end{itemize}

\section{Arnol'd flows and \textit{Good} families} \label{sec:Af}
\subsection{Roof function and passage constant}\label{sec:passtime}
On their $\phi_\R$-invariant minimal components, smooth locally Hamiltonian flows with finitely many saddles are isomorphic to a suspension flow (recall Sections \ref{sec:susflow} and \ref{sec:flowalpha}) over $\mathbb{T}=[0,1)$, where $\mathbb{T}$ parametrizes a transversal section $\Sigma$ via the Hamiltonian. The roof function $\ph$ is given by the first return time to $\Sigma$ under this identification (see \cite{Ravotti2017,KOZ2026} and the references therein). For genus one surfaces, the base map is a rotation $R_\alpha$.

The rotation parameter $\alpha$ can be computed directly from the $1$-form $\lambda$ corresponding to the flow. Let $z_0$ and $z_0'$ denote the first intersections of the backward and forward separatrices of the saddle $p_0$ with $\Sigma$, respectively. Then
\begin{equation}
\tag{\ref{eq:VF-RotNumb}}   
\alpha = \left(\int_{\Sigma}\lambda\right)^{-1}
\int_{z_0}^{z_0'}\lambda.
\end{equation}
More generally, $\lambda$ induces an identification of $\Sigma$ with $\mathbb{T}$ by the coordinate
\[
\Phi(x,0)=z\longmapsto x=\left(\int_{\Sigma}\lambda\right)^{-1}
\int_{z_0}^z\lambda \quad \operatorname{mod} 1.
\]

Moreover, for every non-degenerate saddle point $p_0$, there exists a neighborhood $p_0\in \cU_0$ where we can introduce local coordinates $(u,v)$ in which $p_0=(0,0)$, the Hamiltonian is $H(u,v)=uv$, and there exists a function $V$ smooth with $V> 0$ such that the area form reads $\omega = V(u,v)du\wedge dv$ and the vector field is given by $W = \frac{x}{V (x,y)}\partial x-\frac{y}{V (x,y)}\partial y$ (see, e.g., \cite[$\S$3.1]{Ravotti2017}).

In particular, for Arnol'd flows with a single saddle, we have the following result.

\begin{lemma}\label{lemma:passtime}
    Let $T(c)$ be the passage time around a fixed neighborhood of the saddle $p_0$ at energy level $c$--close, and let $V(u,v)$ be the potential function arising from the coordinates established above; then the return time to $\Sigma$ of $\Phi(x,0)$ is given by 
    \begin{equation}
        \tag{\ref{eq:def_phi}}
        \ph (x) = \psi (x)+ 2c_0|\log x|+c_0|\log(1-x)|,
    \end{equation}
    where $\psi $ is of bounded variation, and
    \be\tag{\ref{eq:passtime}}
        c_0  = \lim_{c\to 0}\frac{T(c)}{|\log c|}=V(0,0)=\frac{1}{\sqrt{-\text{det}DW(p_0)}}>0.
    \ee
\end{lemma} 
The choice between $c_0$ and $2c_0$ on $\ph$ depends on whether the orbit passes once or twice near $p_0$ before returning to the Poincaré section $\Sigma$, the latter corresponding to the side of the saddle of the loop (see Figure \ref{fig:Af}).

\begin{proof} The proof follows the argument in \cite[$\S$3.1]{Ravotti2017}. We include it for completeness. 

In the coordinates established above, the level set $H=c>0$ is the hyperbola $v=c/u$, and the flow satisfies $\dot{u}=u/V(u,c/u)$. Moreover, assume that, for every $x\in\Sigma$, if the orbit segment $\{\Phi(x,t)\}_{t\in[0,\ph(x)]}$ meets $\cU_{0}$, it does so in at most two connected arcs (the latter possibility occurs precisely when $\Phi(x,t)$ is close to a saddle loop homologous to zero), and that $\cU_0$ corresponds to $(-\eps,\eps)^2$ in these coordinates for some $\eps>0$.

Write
\[
    \ph(x) = \int_{\{t:\Phi(x,t)\notin\cU_0\}\cap[0,\ph(x)]}dt + \int_{\{t:\Phi(x,t)\in \cU_0\}\cap[0,\ph(x)]}dt.
\]
Since the flow is uniformly smooth away from the fixed point and the time spent outside $\cU_0$ is uniformly bounded in $x$ (it depends only on the geometry of the flow away from the saddle), the first integral defines a bounded-variation function of $x$, which we absorb into $\psi(x)$. It therefore suffices to analyze the second term.

Let $c=c(x)<\eps^2$ be the energy of the orbit of $\Phi(x,0)$. The length of a maximal interval in $[0,\ph(x)]\cap \{t:\Phi(x,t)\in \cU_0\}$ becomes
\begin{equation}
    \label{eq:Tc}
    T(c) := \int_0^{T(c)} dt = \int_{c/\eps}^{\eps} \frac{V \left(u,c/u\right)}{u}du,
\end{equation}
where the substitution $u = u(t)$ along the hyperbola absorbs the $\dot u$ factor. We can apply \cite[Lemma B.1]{FraczekUlcigrai2010} to \eqref{eq:Tc} to get
\[
    T(c) = V(0,0)|\log c| + \hat \psi(c),
\]
where $\hat \psi$ is in $L^\infty$ with bounded variation. Since $x\in\T$ parametrizes $\Sigma$ (linearly) via the Hamiltonian, $c(x)=H_0\,\dist(x,\Z)$ for some $H_0>0$. Hence $|\log c(x)| = |\log \dist(x,\Z)| + O(1)$, and
\[
   T(c(x)) = V(0,0)|\log \dist(x,\Z)| + O(1).
\]
where the $O(1)$ remainder lies in $L^\infty$ with bounded variation in $x$. Putting together the bounded-variation terms and accounting for both sides of the saddle completes the proof.
    
\end{proof}

\subsection{A measure class on flows} \label{sec:Kfc} To present examples of \textit{good} pairs $(\cF,\hat\eta)$, we first need to define a measure class on the corresponding space of flows. This measure class, commonly referred to as Katok's fundamental class, was introduced in \cite{Katok1973}. We define it here in the context of Arnol'd flows and recall some of its well-known properties. For a more detailed discussion, we refer the reader to \cite{NikolaevZhuzhoma1999,Ravotti2017}.

Let $\cF_{p_0,q_0}$ be the set of smooth closed 1-forms on $\T^2$ with exactly two singular points: a non-degenerate saddle in $p_0$ and a center on $q_0$. Fix a basis $\{\gamma_1,\gamma_2,\gamma_3\}$ of the first relative homology group $H_1(\T^2,\{p_0,q_0\},\R)$. For every
$\lambda\in \cF_{p_0,q_0}$, define the periodic coordinates
\[\Theta(\lambda)=\left(\int_{\gamma_1}\lambda,\int_{\gamma_2}\lambda,\int_{\gamma_3}\lambda\right)\in\R^3.\]
The function $\Theta$ is well defined on a neighborhood of $\lambda$ in $\cF_{p_0,q_0}$ and it is a complete invariant for isotopy classes. That is, for every $\lambda\in\cF_{p_0,q_0}$, there exists a neighborhood $\cU$ of $\lambda$ such that for any $\lambda_0\in\cU$, $\Theta(\lambda)=\Theta(\lambda_0)$ if and only if there exists a family of 
smooth maps $\{h_t:\T^2\to\T^2\}_{t\in[0,1]}$, such that $h_1^*\lambda=\lambda_0$.
In particular, the orbits of the locally Hamiltonian flows induced by $\lambda$ and $\lambda_0$ are equivalent.

For each $p_0,q_0\in \T^2$ ($p_0\neq q_0$), we equip $\cU\subset \cF_{p_0,q_0}$ with the measure class induced by the pull-back of the Lebesgue measure on $\R^3$ via $\Theta$, namely $\Theta^*(\text{Leb}_{\R^3})$. For each $\lambda\in \cF_{p_0,q_0}$, recall that $\lambda$ generates an Arnol'd flow with a single saddle, and let $\alpha(\lambda)$ denote the rotation number associated to this flow (see Section \ref{sec:flowalpha}). Then, for every $\lambda\in \cF_{p_0,q_0}$ we can choose a basis of the homology group such that $\alpha\int_\Sigma\lambda$ and $\int_\Sigma\lambda$ appear as coordinates of $\Theta(\lambda)$ \cite[$\S$3.1]{Ravotti2017}. 
Thus, every null Lebesgue measure set of $\alpha$ in $\T$ corresponds to a null measure set on $\R^3$, which corresponds to a null set on $\cF_{p_0,q_0}$.

\subsection{Examples of \textit{good} families of vector fields} \label{sec:ex}
We present two examples of \textit{good} families in the sense of Definition \ref{def:good}, one arising from the measure class introduced in Section \ref{sec:Kfc}, and a generic one-parameter analytic family.

\begin{example}
Fix $p_0$ and $q_0$ in $\T^2$ and let $\cF\subset\cF_{p_0,q_0}$ be a positive measure set of Arnol'd flows with a fixed non-degenerate saddle $p_0$ and a fixed center $q_0$, with respect to the measure class introduced in Section \ref{sec:Kfc}. We know $\alpha(\cF)\subset \T$ has positive Lebesgue measure, so let $\eta$ be a probability measure on $\T$, supported on $\alpha(\cF)$ and absolutely continuous with respect to $\text{Leb}_\T$. Picking any probability measure $\hat\eta$ in $\cF$ such that $\eta=\alpha_*\hat\eta$, the pair $(\cF,\hat\eta)$ is \textit{good}.
\end{example}

\begin{example}
Let $I\subset\R$ be a small neighborhood of $0$ and let $\{\lambda(t)=\lambda_0+t\lambda_1\}_{t\in I}$ be a 1-parameter analytic family of closed 1-forms in $\cF_{p_0,q_0}$. As discussed in Section \ref{sec:Kfc}, we can choose a basis $\{\gamma_1,\gamma_2,\gamma_3\}$ of the homology group $H_1(\T^2,\{p_0,q_0\},\R)$ such that the periodic coordinates \[
\alpha(t)=\frac{\int_{\gamma_1}\lambda(t)}{\int_{\gamma_2}\lambda(t)}
=\frac{\int_{\gamma_1}(\lambda_0+t\lambda_1)}{\int_{\gamma_2}(\lambda_0+t\lambda_1)}=\alpha(0)+\frac{t\left (\int_{\gamma_1}\lambda_1-\alpha(0)\int_{\gamma_2}\lambda_1\right)}{\int_{\gamma_2}\lambda_0+t\int_{\gamma_2}\lambda_1}.
\]
Therefore $\alpha(t)$ is an analytic function of $t$ and it is non-constant precisely when \[\int_{\gamma_1}\lambda_1\neq\alpha(0)\int_{\gamma_2}\lambda_1.\] This is a generic condition. Hence, for a generic perturbation $\lambda_1$, if $t$ is distributed according to a bounded density on $I$, the induced distribution of $\alpha(t)$ is absolutely continuous with respect to Lebesgue.
\end{example}

Finally, we note that one can consider time reparametrizations of the flows of a \textit{good} pair $(\cF,\hat\eta)$ to impose $c_0^W$ to be constant. 
\begin{example} Let $a\in\R^{+}$ be fixed and set $\cF' := \{{a\,W}/{\sqrt{-\det DW(p_0)}} : W\in\cF\}$. For $W\in\cF$ and $c(W):=a/\sqrt{-\det DW(p_0)}>0$, we have $\det D(c(W)W)(p_0) = -a^2$, independently of $W$. By \eqref{eq:passtime}, the passage constant for any $W\in\cF'$ is equal to $c_0:=a^{-1}$.

\begin{claim*} Let $\hat\eta'$ be the measure on $\cF'$ given by the push-forward of $\hat\eta$ on $\cF$ under the scaling map $W\mapsto c(W)W$. Then $(\cF',\hat\eta')$ is \textit{good} in the sense of Definition \ref{def:good}.
\end{claim*}
\begin{proof}
Since $c(W)>0$, the flow generated by $c(W)W$ is a constant-speed time reparametri-zation of the flow generated by $W$. 
In particular, the saddle $p_0$, the section $\Sigma$, and the induced base map $R_\alpha$ are unchanged, so $\alpha(\cdot)$ is preserved under this scaling, and $\alpha(\rv W')\sim\alpha(\rv W)$ for $\rv W\sim\hat\eta$ in $\cF$, $\rv W'\sim\hat\eta'$ in $\cF'$. Therefore, the pair is \textit{good}.
\end{proof}
\end{example}

\section{Proof of Theorem \ref{thm:main}}\label{sec:mainproof}
Here, we use Theorems \ref{thm:mainbis} and \ref{thm:mainbisbis} to prove Theorem \ref{thm:main}. The remainder of the paper, Sections \ref{sec:proof0}--\ref{sec:proof2}, is devoted to the proof of Theorem \ref{thm:mainbis}. In Appendix \ref{app:1}, we provide a sketch of the elementary modifications needed to prove Theorem \ref{thm:mainbisbis}.
\begin{proof}
First, we prove the following two claims.
\begin{claim}For each any $A\subset\cF$ with $\hat\eta(A)>0$, set \[d\hat{\eta}_{A}(W):=\hat\eta(A)^{-1} \chi_{A}(W)d\hat \eta(W).\]
    Then $(\cF,\hat\eta_{A})$ is \textit{good} in the sense of Definition \ref{def:good}.
\end{claim}
\begin{proof}
    It follows from the fact that $(\cF,\hat\eta)$ is \textit{good} and that $\hat\eta(A)^{-1} \chi_{A}(W)$ is uniformly bounded for $W\in\cF$.
\end{proof}
\begin{claim} Without loss of generality, we may assume there $\exists c_0>0$ such that $c_0^W\leq c_0$ for $\hat\eta$-almost every $W\in\F$.
\end{claim} 
\begin{proof}
Assume the theorem holds under the hypothesis that $c_0^W$ is uniformly bounded; we show it then holds in full generality.
    
    For every $N\in\N$, let $B_N:=\{W:c_0^W\leq N\}$. Then, since $c_0^W<\infty$ for every $W\in\cF$, we have $\hat\eta(B_N)\to1$ as $N\to\infty$. Let $Y\sim \frac{1}{4\pi}\,\mathfrak{C}$, let
\[F_t:=\frac{\int_0^{t} f \left(\phi^{\rv W}_s\rv P\right)ds - {t}\bar f^{\rv W}}{(\log t)^2},\]
and fix a continuous bounded function $g:\R\to\R$. Then, for $N$ large enough such that $\hat\eta(B_N)>0$, it holds that
\[
|\E\left[g(F_{t})-g(Y)\right]|\leq|\E\left[g(F_{t})-g(Y)|B_N\right]|+2\norm{g}_{\infty}\P(B_N^c).
\]
By Claim 1, $(B_N, \hat{\eta}_{B_N})$ is \textit{good} in the sense of Definition \ref{def:good} and, by definition, $c_0^W$ is uniformly bounded on $B_N$. Hence, by our assumption s$|\E\left[g(F_{t})-g(Y)|B_N\right]|\to 0$ as $t\to\infty$ for every $N$ large enough. Finally, note $2\norm{g}_{\infty}\P(B_N^c)\to 0$ as $N\to\infty$. This concludes the proof.
\end{proof}
Assume there exists $c_0$ as in Claim 2. For every $N\geq 1$, set \[\Delta_Nf(W):=\text{sign}(f(p_0^{W})-\bar f^W)\frac{\left\lceil{N|f(p_0^{W})-\bar f^W|c_0^{W}}\right\rceil}{N}.\]
Note that $|\Delta_Nf(W)|\leq 2c_0\norm{f}_\infty+1$, it takes a finite set of values, and we have uniform convergence $\Delta_Nf(W)\to c_0^W(f(p_0^{W})-\bar f^W)$ as $N\to\infty$. Moreover, note that $\Delta_Nf(W)=0$ if and only if $\bar f^W=f(p_0^W)$.

Let $\rv D=c_0^{\rv W}(f(p_0^{\rv W})-\bar f^{\rv W})$, where $\rv W\sim\hat\eta$, and $\rv D_N=\Delta_Nf({\rv W})$. Denote by \[\{a_{N,1},\ldots,a_{N,r_N}\}=\{a\in \R^*:\P(\rv D_N=a)>0\},\] where $r_N\geq 0$. For $1\leq i\leq r_N$, let $A_{N,i}=\{W\in\cF:\Delta_N(W)=a_{N,i}\}$, which have positive $\hat\eta$-measure by definition of $a_{N,i}$. Then
\[\rv D_N=\sum_{i=1}^{r_N}a_{N,i}\chi_{A_{N,i}}(\rv W).\]

For each $N,i$, denote $\hat\eta_{N,i}=\hat\eta_{A_{N,i}} $. By Claim 1, $(\cF,\hat\eta_{N,i})$ is \textit{good} in the sense of Definition \ref{def:good}.
By definition of $A_{N,i}$, $\hat\eta_{N,i}$-almost surely $f(p_0^{W})-\bar f^W\neq 0$. Hence, by Theorem \ref{thm:mainbis}, taking $\rv W_{N,i}\sim\hat \eta_{N,i}$ in $\cF$, and then $\rv P$ uniform in each $\M^W$, we have
\[
\frac{\int_0^{t} f \left(\phi^{\rv W_{N,i}}_s\rv P\right)ds - {t}\bar f^{\rv W_{N,i}}}{c_0^{\rv W_{N,i}}\left(f(p_0^{\rv W_{N,i}})-\bar f^{\rv W_{N,i}}\right)(\log t)^2}
\xrightarrow[t\to\infty]{\text{dist}}\ \frac{1}{4\pi}\,\mathfrak{C}.
\]
Since the $\{A_{N,i}\}_i$ are pair-wise disjoint for every $N\geq 1$, it follows that 
\[\begin{aligned}
\frac{\rv D_N}{\rv D}\frac{\int_0^{t} f \left(\phi^{\rv W}_s\rv P\right)ds - {t}\bar f^{\rv W}}{(\log t)^2}
&=\sum_{i=1}^{r_N}a_{N,i}\chi_{A_{N,i}}(\rv W)\frac{\int_0^{t} f \left(\phi^{\rv W}_s\rv P\right)ds - {t}\bar f^{\rv W}}{{c_0^{\rv W}}\left(f(p_0^{\rv W})-\bar f^{\rv W}\right)(\log t)^2}\ \xrightarrow[t\to\infty]{\text{dist}}\ \frac{\rv D_ N}{4\pi}\,\mathfrak{C},
\end{aligned}
\]
with the convention that $\frac{\rv D_N}{\rv D}=0$ on $A_0:=\{W\in\cF:f(p_0^{W})-\bar f^W= 0\}$.

Also, by Theorem \ref{thm:mainbisbis}, if $\hat\eta(A_0)>0$, then
\[
\chi_{A_{0}}(\rv W)\frac{\int_0^{t} f \left(\phi^{\rv W}_s\rv P\right)ds - {t}\bar f^{\rv W}}{(\log t)^2}\xrightarrow[t\to\infty]{\text{dist}}\  0.
\]

Now, let $Y\sim \frac{1}{4\pi}\,\mathfrak{C}$ and, for every $t\geq 0$ and $N\geq 1$, let 
\[Z_N:=\frac{\rv D_N}{\rv D}+\chi_{A_0}(\rv W),\quad F_t:=\frac{\int_0^{t} f \left(\phi^{\rv W}_s\rv P\right)ds - {t}\bar f^{\rv W}}{(\log t)^2}.\]
We know that $\rv D_N\to \rv D$ uniformly as $N\to\infty$, so $Z_N\to 1$ almost-surely. Moreover, from the derivations above, it holds that $Z_NF_t\to D_NY$ in distribution as $t\to\infty$, so the family $\{Z_N F_t\}_t$ is tight. Noting $|Z_N|\geq 1$, this implies that $\{F_t\}_t$ is tight.

Fix a continuous bounded function $g:\R\to\R$ and $\eps>0$. Let $M>0$ be such that $\P(2\norm{f}_\infty|Y|>M)<\eps$ and $\P(|F_t|>M)<\eps$ for any $t$ large enough, and let $1>\delta>0$ be such that for any $|x|,|y|\leq M+1$ with $|x-y|<\delta$, $|g(x)-g(y)|<\eps$. Finally, fix $N$ large enough such that $|Z_N-1|<\delta/M$ $\eps$-almost surely and $|\rv D_N-\rv D|<2\delta\norm{f}_\infty/M$ uniformly. Then
\[\begin{aligned}
&|\E\left[g(F_t)-g(\rv DY)\right]|
\\&\qquad \leq\E\left[|g(F_t)-g(Z_NF_t)\right|]+|\E\left[g(Z_NF_t)-g(\rv D_NY)\right]|+\E\left[|g(\rv D_N Y)-g(\rv DY)\right|].
\end{aligned}
\]
The first term
\[\begin{aligned}
\E\left[|g(F_t)-g(Z_NF_t)\right|]\leq\eps+2\norm{g}_\infty\P(|F_t|>M)+2\norm{g}_\infty\P(|F_t-Z_NF_t|\geq\delta, \ |F_t|\leq M)\\\leq\eps+2\norm{g}_\infty\P(|F_t|>M)+2\norm{g}_\infty\P(|1-Z_N|\geq\delta/M)<\eps+4\eps\norm{g}_\infty.
\end{aligned}
\]
For $N$ fixed, for the second term we have the limit
\[
|\E\left[g(Z_NF_t)-g(\rv D_NY)\right]|\xrightarrow{t\to\infty}0.
\]
Finally, for the third term
\[\begin{aligned}
\E\left[|g(\rv D_N Y)-g(\rv DY)\right|]\leq\eps+2\norm{g}_\infty\P(
2\norm{f}_\infty|Y|>M)<\eps+2\eps\norm{g}_\infty.\end{aligned}
\]
Therefore $\limsup_{t\to\infty}|\E\left[g(F_t)-g(\rv DY)\right]|=O(\eps)$ for every $\eps>0$, so it is equal to 0. 

Since $g$ was arbitrary, it follows that $F_t\to \rv DY$ in distribution as $t\to\infty$. This concludes the proof of the theorem.
\end{proof}
\section{Main steps of the proof of Theorem \ref{thm:mainbis}}\label{sec:proof0}

The proof combines three main ingredients, namely Propositions \ref{prop:main_dist}, \ref{prop:main_red}, and \ref{prop:main_harm}. 

Let $W\in\cF$. Recall from Section \ref{sec:susflow} that for $\omega$-almost every $p\in M^W$, there exist $x\in(0,1)$ and $\tau\in[0,\ph^W(x))$ such that $\Phi^W(x,\tau)=p$. 
We define
\[
I(W,x,t):=\int_0^{t} f \left(\phi_{s}^W\left(\Phi^W(x,0)\right)\right)ds - t\bar f^W,\qquad \bar f^W:=\frac{1}{\omega(\M^W)}\int_{\M^W}f(p)\,dp.
\]

\noindent\textbf{Step 1:} The first step, Proposition \ref{prop:main_dist} allows us to consider the function $I(\rv W,\rv X,t)$ and take a distribution on $\rv X\big|\rv W$, instead of $\int_0^{t} f \left(\phi^{\rv W}_s\rv P\right)ds - t\bar f^W$ and $\rv P$ uniform in $\M^{\rv W}$.

\begin{proposition}\label{prop:main_dist}
    Let ($\cF,\hat\eta$) be \textit{good} in the sense of Definition \ref{def:good}. Let $\epsilon>0$ and let $f\in C^{2+\epsilon}(\T^2)$, and suppose $f(p_0^{W})-\bar f^W\neq 0$, for $\bar f^W$ as in \eqref{eq:meanf}, holds $\hat\eta$-almost surely. Suppose that for any family of probability distributions
    $ (\hat\rho_W)_{W\in\cF}$ uniformly $ C^2(\T)$\footnote{Recall we say $(\hat \rho_W)_{W\in \cF}\subset C^2(\T)$ is uniformly $C^2(\T)$ if $\sup_\cF\norm{\hat\rho_W}_{C^2}<\infty$.}, it holds that, under $\rv W\sim\hat\eta$ in $\cF$ and $\rv X\big|\rv W=W\sim\rho_{W}$ we have
\[
\frac{I(\rv W,\rv X,t)}{{c_0^{\rv W}}\left(f(p_0^{\rv W})-\bar f^{\rv W}\right)(\log t)^2}\ \xrightarrow[t\to\infty]{\text{dist}}\ \frac{1}{4\pi}\,\mathfrak{C},
\]
where $\mathfrak{C}$ is the standard Cauchy distribution. Then Theorem \ref{thm:mainbis} holds.
\end{proposition}

For each $W\in\cF$, let $\mu^W:=\int_{\T}\ph^W(x)dx$. We know from the definition of $\ph$ as \eqref{eq:def_phi} that $0<\mu^W<\infty$. The natural distribution induced by the invariant measure of the suspension flow has density proportional to the roof function (see Lemma \ref{lemma:distx}),
\[
\rv X\big|\rv W=W\sim\frac{\ph^{W}(X)}{\mu^W},
\]
which inherits the logarithmic singularities of $\ph^W$ for each $W$ and is therefore not uniformly $C^2$. On the other hand, our limit theorem will be established first for arbitrary families of uniformly $C^2$ probability densities, since this regularity is convenient in the harmonic analysis carried out in Section \ref{sec:proof2}. 

The role of Proposition \ref{prop:main_dist} is showing that this restriction is inessential: convergence for all uniformly $C^2$ densities automatically extends to every strictly positive $L^1$ density (see Lemma \ref{lemma:relax}), in particular to the densities $(\ph^W/\mu^W)_W$, and ultimately to the sampling of initial conditions according to $\omega$ on $\M^W$.

\noindent\textbf{Step 2:} We define the (mean zero) singular profile of $\ph $ by 
\be
\tag{\ref{eq:sing_p}}
h(x)= 2|\log x|+|\log(1-x)|-3.
\ee
Then \eqref{eq:def_phi} yields $\ph^W(x)=c_0^Wh(x)+3c_0^W+\psi^W(x)$. 

Given $N\geq 0$, $\alpha\in\R$, $x\in \T$, and $g\in L^2(\T)$, define the Birkhoff sum of $g$
\[
    S(\alpha,g,x,N):=\sum_{m=0}^{N-1}g(R_\alpha^mx).
\]

The return times of the induced suspension flow contain logarithmic singularities created by the saddle, suggesting that this singular part should encode all large fluctuations of Birkhoff integrals of a function $f$. The goal of Proposition \ref{prop:main_red} is to make this heuristic precise by showing that, after the normalization $(\log t)^2$, every regular contribution becomes negligible, and only the logarithmic profile $h$ survives.

\begin{proposition}
    \label{prop:main_red} Let ($\cF,\hat\eta$) be \textit{good} in the sense of Definition \ref{def:good}. Let $\epsilon>0$ and let $f\in C^{2+\epsilon}(\T^2)$, and suppose $f(p_0^{W})-\bar f^W\neq 0$, for $\bar f^W$ as in \eqref{eq:meanf}, holds $\hat\eta$-almost surely. 
    Let $\rv W\sim\hat\eta$ and $\rv X\big|\rv W=W\sim \hat\rho_{W}$ uniformly $C^2(\T)$. If there exists a continuous random variable $Y$ such that
    \be \label{eq:convY}
    \frac{S\left(\alpha(\rv W),h,\rv X,N\right)}{(\log {N})^2}\ \xrightarrow[N\to\infty]{\text{dist}}\ Y,
    \ee
    for $h$ defined in \eqref{eq:sing_p}, then 
        \be\label{eq:distrib_sum}
    \frac{I(\rv W,\rv X,t)}{{c_0^{\rv W}}\left(f(p_0^{\rv W})-\bar f^{\rv W}\right)(\log t)^2}\ \xrightarrow[t\to\infty]{\text{dist}}\ Y.
    \ee
\end{proposition}

\noindent\textbf{Step 3:} Finally, Proposition \ref{prop:main_harm} establishes the convergence of $S(\rv A,h,\rv X,N)$ to a Cauchy distribution, for $\rv A$ with smooth density on $\T$, after normalizing by $(\log N)^2$. Since the function on \eqref{eq:convY} is measurable with respect to $\alpha(\rv W)$, taking $\rv A=\alpha(\rv W)$ in Proposition \ref{prop:main_harm} will give us the limit law we need to apply Proposition \ref{prop:main_red}.

The Proofs of Propositions \ref{prop:main_dist}, \ref{prop:main_red}, and \ref{prop:main_harm} are carried out in Sections \ref{sec:proof3}, \ref{sec:proof1}, and \ref{sec:proof2}, respectively. Here, we use them to deduce Theorem \ref{thm:mainbis}.

\begin{proof}[Proof of Theorem \ref{thm:mainbis}]
    Given $\rv W\sim\hat\eta$ in $\cF$, let $\eta$ be the distribution of its rotation parameter $\rv A=\alpha(\rv W)$, which is absolutely continuous with respect to Lebesgue$(\T)$ by assumption.
    
    Fix a family of probability distributions $(\hat \rho_W)_{W\in\cF}$ uniformly $ C^2(\T)$ and let $( \rho_\alpha)_{\alpha\in\text{Supp}(\eta)}$ be the equivalent distributions with respect to $\alpha$. That is, letting $\E_{\rv W}$ denote the expectation with respect to $\rv W$, $\rho_\alpha(x)= \E_{\rv W}[\hat\rho_{\rv W}(x)|\alpha(\rv W)=\alpha]$ for $\eta-$almost every $\alpha$. Note that, by the linearity of $\E_{\rv W}$ and the dominated convergence theorem, $( \rho_\alpha)_{\alpha\in\text{Supp}(\eta)}$ are also uniformly $C^2(\T)$.

    By assumption $\eta\in L^1(\T)$, so there exists a sequence of smooth probability densities $\eta_n\to\eta$ in ${L^1}$. Let $\rv A_n\sim \eta_n$ for every $n$, coupled with $\rv A\sim\eta$ such that $\P(\rv A_n\neq \rv A)\leq\int_\T|\eta_n(x)-\eta(x)|dx\to0$, as $n\to\infty$. Define for $n\geq 1$ random variables $\rv X_n$: $\rv X_n\big|\rv A_n=\alpha\sim\rho_\alpha$.
    
    Let $F_N(x,\alpha)=\frac{S(\alpha,h,x,N)}{(\log N)^2}$ for $N\geq 1$, $x,\alpha\in\T$. By Proposition \ref{prop:main_harm}, for every $n$ it holds 
    \[F_N(\rv A_n ,\rv X_n)\xrightarrow[N\to\infty]{\dist} \tfrac{1}{4\pi} \mathfrak{C}.\] 

    Fix a continuous bounded function $g:\T\times\T\to\R$ and let $\rv X$ be a random variable with $\rv X|\rv A=\alpha\sim\rho_\alpha$. Then for every $n$ and $\rv A$, $\rv A_n$ and $\rv X_n$ as above,
    \[\begin{aligned}
    \left|\E[g(F_N(\rv A,\rv X))]-\E[g(Y)]\right|&\leq\left|\E[g(F_N(\rv A_n,\rv X_n))]-\E[g(Y)]\right|+2\norm{g}_\infty\P(\rv A\neq \rv A_n)
    \\&\xrightarrow[N\to\infty]{}2\norm{g}_\infty\P(\rv A\neq \rv A_n) \xrightarrow[n\to\infty]{}0.
    \end{aligned}
    \]

    Therefore, under $\rv A\sim \eta$ and $\rv X|\rv A=\alpha\sim\rho_\alpha$, we have the convergence result:
    \[
    \frac{S\left(\rv A,h,\rv X,N\right)}{(\log {N})^2}\ \xrightarrow[N\to\infty]{\text{dist}}\  Y,
    \] 
    where $Y\sim (4\pi)^{-1}\mathfrak{C}$ and $\mathfrak{C}$ is the standard Cauchy distribution. Moreover, by the definitions of $\eta$ and $(\rho_\alpha)_\alpha$, this is equivalent to convergence under $\rv W\sim\hat\eta$ and $\rv X\big|\rv W=W\sim\hat\rho_W$ of 
 \[
    \frac{S\left(\alpha(\rv W),h,\rv X,N\right)}{(\log {N})^2}\ \xrightarrow[N\to\infty]{\text{dist}}\ Y.
\]
Therefore, the hypothesis of Proposition \ref{prop:main_red} are fulfilled and we can conclude
    \[
    \frac{I(\rv W,\rv X,t)}{c_0^{\rv W}\left(f(p_0^{\rv W})-\bar f^{\rv W}\right)(\log {t})^2}\ \xrightarrow[t\to\infty]{\text{dist}}\ \frac{1}{4\pi}\mathfrak{C}.
    \]
Applying Proposition \ref{prop:main_dist} concludes the proof.
\end{proof}

\section{Proof of Proposition \ref{prop:main_dist}}\label{sec:proof3}

Let $\cF $ be a family of vector fields in $\T^2$, let $\hat\eta$ be a probability measure on $\cF$, let $\epsilon>0$, and let $f\in C^{2+\epsilon}(\T^2)$. Throughout this Section, assume $(\cF,\hat\eta)$ is \textit{good} in the sense of Definition \ref{def:good}, and that $f(p_0^{W})-\bar f^W\neq 0$ holds $\hat\eta$-almost surely. .

Recall that we define, for each $x\in\T$, $W\in\cF$, and $t\geq 0$,
\[
I(W,x,t):=\int_0^{t} f \left(\phi_{s}^W\left(\Phi^W(x,0)\right)\right)ds -t\bar f^W,\qquad \bar f^W:=\frac{1}{\omega(\M^W)}\int_{\M^W}f(p)\,dp,
\]
with the identification $(x,\ph^W(x))\sim$\footnote{Here $\sim$ denotes an identification, not a probability distribution.}$(R_{\alpha(W)}x,0)$ at the points where $\ph^W(x)<\infty$, and
\[\mu^W:=\int_T\ph^W(x)dx.\]
Moreover, recall that for $\omega$-almost every $p\in M^W$, there exist $x\in(0,1)$ and $\tau\in[0,\ph^W(x))$ such that $\Phi^W(x,\tau)=p$, so
$
\int_0^{t} f \left(\phi_{t}^Wp\right)dt =\int_0^{t} f \left(\phi_{t}^W\left(\Phi^W(x,\tau)\right)\right)dt.
$

The proof of Proposition \ref{prop:main_dist} proceeds in two steps. First, we show we can reduce the problem to considering the function $I(\rv W,\rv X,t)$ and taking the conditional distribution on $\rv X\big|\rv W=W\sim\ph^W/\mu^W$, instead of $\int_0^{t} f \left(\phi^{\rv W}_s\rv P\right)ds - t\bar f^W$ and $\rv P$ uniform in each $\M^W$. 

\begin{lemma}\label{lemma:distx}
    Suppose that there exists a random variable $Y$ such that, under $\rv W\sim\hat\eta$ and $\rv X\big|\rv W=W\sim\ph^W/\mu^W$
     \[
     \frac{I(\rv W,\rv X,t)}{{c_0^{\rv W}}\left(f(p_0^{\rv W})-\bar f^{\rv W}\right)(\log t)^2}\ \xrightarrow[t\to\infty]{\text{dist}}\ Y,
     \]
     then taking $\rv W\sim\hat \eta$ in $\cF$ and $\rv P|\rv W=W\sim\text{Uniform}(\M^W)$ we have
    \be\label{eq:convfin}
    {\frac{\int_0^{t} f \left(\phi^{\rv W}_s\rv P\right)ds - t\bar f^W}{{c_0^{\rv W}}\left(f(p_0^{\rv W})-\bar f^{\rv W}\right)(\log t)^2}\ \xrightarrow[t\to\infty]{\text{dist}}\ Y.}
    \ee
\end{lemma}

Then, we use a standard argument to relax the (singular) profiles $(\ph^W/\mu^W)_W$ to a family of probability distributions $(\hat \rho_W)_W$ uniformly $C^2(\T)$. 

\begin{lemma} \label{lemma:relax}
    Let $\cF$ be a set, let $F_t:\mathcal{F\times\T}\to\R$ be a family of scalar functions for $t\geq 0$, and let $\rv W$ be a fixed $\cF$-valued random variable. Suppose there exists a random variable $Y$ such that \[F_t(\rv W,Z)\xrightarrow[t\to\infty]{\dist} Y\qquad \text{ whenever}\quad Z\big|\rv W=W\sim\rho_W,\] for any family of probability distributions $(\rho_W)_W$ uniformly $C^2(\T)$. Then for any family of probability distributions $(\ph_W)_W\subset L^1(\T)$, the same distributional limit holds taking $Z\big|\rv W=W\sim \ph_W$.
\end{lemma}

We now combine both results to conclude the Proposition. Then, we prove the lemmas.

\begin{proof}[Proof of Proposition \ref{prop:main_dist}]
For each $W\in\cF$, let $\mu^W=\int_\T\ph^W(x)\,dx\geq \inf_\T\ph^W(x)>0$ and $\ph_W:=\ph^W/\mu^W$. Since $\ph^W\in L^1(\T)$, each $\ph_W$ is a probability density in $L^1(\T)$. Let \[F_t(W,x):=\frac{I(W,x,t)}{c_0^{\rv W}\left(f(p_0^{W})-\bar f^W\right) (\log t)^2}.\] 

By hypothesis, we can apply Lemma \ref{lemma:relax} with this family $(\ph_W)_{W\in\cF}$ and $F_t$, and conclude that, when $\rv W\sim\hat\eta$ in $\cF$ and $\rv X\big|\rv W=W\sim\ph^W/\mu^W$,
\[
\frac{I(\rv W,\rv X,t)}{c_0^{\rv W}\left(f(p_0^{\rv W})-\bar f^{\rv W}\right) (\log t)^2}\ \xrightarrow[t\to\infty]{\text{dist}}\ \frac{1}{4\pi}\,\mathfrak{C}.
\]

This is precisely the hypothesis of Lemma \ref{lemma:distx} (with $Y=\frac{1}{4\pi}\mathfrak{C}$), so applying that 
lemma yields, for $\rv W\sim\hat\eta$ in $\cF$ and $\rv P\big|\rv W=W\sim\text{Uniform}(\M^W)$,
\[
\frac{\int_0^{t} f \left(\phi^{\rv W}_s\rv P\right)ds - t\bar f^W}{c_0^{\rv W}\left(f(p_0^{W})-\bar f^W\right) (\log t)^2}\ \xrightarrow[t\to\infty]{\text{dist}}\ \frac{1}{4\pi}\,\mathfrak{C}.
\]
This is Theorem \ref{thm:mainbis}.
\end{proof}

\begin{proof}[Proof of Lemma \ref{lemma:distx}]
For every $W\in\cF$ and almost every $p=\ph^W(x,\tau)\in\M^W$, with $x=x(p)\in(0,1)$ and $\tau=\tau(p)\in[0,\ph^W(x))$, we have the almost-surely point-wise limit
\[\begin{aligned}
&\frac{\left|\int_0^{t} f \left(\phi_{s}^Wp\right)ds-\int_0^{t} f \left(\phi_{s}^W\left(\Phi^W(x,0)\right)\right)ds\right|}{c_0^{W}\left|f(p_0^{W})-\bar f^W\right| (\log t)^2}
\\
&\qquad =\frac{\left|\int_{t}^{t+\tau} f \left(\phi_{s}^W\left(\Phi^W(x,0)\right)\right)ds-\int_0^{\tau} f \left(\phi_{s}^W\left(\Phi^W(x,0)\right)\right)ds\right|}{c_0^{W}\left|f(p_0^{W})-\bar f^W\right|(\log t)^2}
\\&\qquad\leq\frac{2\tau\norm{f}_{\infty}}{c_0^{W}\left|f(p_0^{W})-\bar f^W\right|(\log t)^2}\xrightarrow[t\to\infty]{}0. 
\end{aligned}
\]
Therefore, under $\rv W\sim\hat \eta$ in $\cF$ and $\rv P|\rv W=W\sim\text{Uniform}(\M^W)$,
\[
\frac{\left|\int_0^{t} f \left(\phi_{t}^{\rv W}\rv P\right)dt-t\bar f^W-I(\rv W,x(\rv P),t)\right|}{c_0^{\rv W}\left|f(p_0^{\rv W})-\bar f^{\rv W}\right|(\log t)^2}\xrightarrow[t\to\infty]{} 0\quad.
\]

Now, the random variable $I(\rv W,x(\rv P),t)$ is measurable with respect to $(\rv W,x(\rv P))$, and $x(\rv P)|\rv W=W\sim\omega|_{\M^W}/\omega(\M^W)$ has the same law as $\rv X\big|\rv W=W\sim\ph^W/\mu^W$.
To see this, recall that $\Phi^W : X^W \to \M^W$, where $X^W = \{(x,s) : x \in \T,\ 0 \le s < \varphi^W(x)\}$, is the composition of a linear rescaling of the base coordinate by a constant $c^W>0$ with a measure-preserving map. Hence $\Phi^W$ carries the Lebesgue measure $dx\,ds$ on $X^W$ to $\omega|_{\M^W}/c^W$. Integrating for $s \in [0,\varphi^W(x))$ at fixed $x$ gives $\varphi^W(x)\,dx$ as the marginal of $dx\,ds$ on the base. Finally, normalizing by its total mass $\mu^W := \int_\T \varphi^W(x)\,dx = \omega(\M^W)/c^W$ cancels $c^W$, so the marginal law of the base coordinate under the $\phi_t$-invariant probability measure on $\M^W$ is exactly $\varphi^W/\mu^W$.

Since 
\[\frac{I(\rv W,\rv X,t)}{{c_0^{\rv W}}\left(f(p_0^{\rv W})-\bar f^{\rv W}\right)(\log t)^2}\ \xrightarrow[t\to\infty]{\text{dist}}\ Y,\] by hypothesis and the difference above vanishes almost surely as $t\to\infty$, Slutsky's theorem yields \eqref{eq:convfin}.
\end{proof}
\begin{proof}[Proof of Lemma \ref{lemma:relax}]
    $C^2(\T)$ is dense in $L^1(\T)$, thus for any $\eps>0$ there exists $(\hat\rho_W^\eps)_W\subset C^2(\T)$ such that 
    \be\label{eq:cond1}\sup_W\int_\T|\ph_W(x)-\hat\rho_W^\eps(x)|dx<\eps\ee 
    Moreover, since $\ph_W$ are probability distributions, we can further assume $\hat\rho_W^\eps\geq 0 $ and $ \int_\T \hat\rho_W^\eps(x)dx=1$, so that $(\hat\rho_W^\eps)_W$ are probability distributions.
    Also, by boundedness almost surely of $\norm{\hat\rho_W^\eps}_{C^2}$, there exists $\cF_\eps$ with probability $>1-\eps$ such that $(\hat\rho_W^\eps)$ restricted to ${W\in\cF_\eps}$ is uniformly $C^2(\T)$. Define for every $x\in\T$
    \be\label{eq:cond2}\rho_W^\eps(x)=
    \begin{cases}
        \hat\rho_W^\eps(x)&\text{for }W\in \cF_\eps,\\
        1&\text{for }W\notin \cF_\eps.
    \end{cases}
    \ee
    $(\rho_W^\eps)_W$ is a family of probability distributions uniformly $C^2(\T)$, so  \[F_t(\textbf{}W, Z)\xrightarrow[t\to\infty]{\dist} Y\qquad \text{when}\quad Z\big|\rv W=W\sim\rho_W^\eps.\]

    Let $X$ and $Z_\eps$ be two random variables on $\T$ satisfying $X\big|\rv W=W\sim\ph_W$ and $Z_\eps|\rv W=W\sim\rho_W^\eps$, and couple them so that $\P(X\neq Z|\rv W=W)\leq\eps$ provided $W\in {\mathcal{F_\eps}}$. This can be done because of \eqref{eq:cond1}--\eqref{eq:cond2}. Then
    \[
        \P(X\neq Z_\eps)\leq\P(X\neq Z_\eps|\rv W\in{\mathcal{F_\eps}})+\P(\cF_\eps^c)<2\eps \quad\text{and}\quad F_t(\rv W,Z_\eps)\xrightarrow[t\to\infty]{\dist} Y.\]
    Hence for any bounded continuous $g:\R\to\R$,
    \[
    |\E[g(F_t(\rv W,X))]-\E[g(Y)]|\leq 2\|g\|_\infty\,\P(X\neq Z_\eps)+|\E[g(F_t(\rv W,Z_\eps))]-\E[g(Y)]|.
    \]
    Letting $t\to\infty$ (for fixed $\eps$) kills the second term, giving
    \[
    \limsup_{t\to\infty}|\E[g(F_t(\rv W,X))]-\E[g(Y)]|\leq 4\eps\|g\|_\infty,
    \]
    and letting $\eps\to0$ yields $\E[g(F_t(\rv W,X))]\to\E[g(Y)]$ for every such $g$, i.e. $F_t(\rv W,X)\xrightarrow[t\to\infty]{\dist} Y.$

\end{proof}

\section{Proof of Proposition \ref{prop:main_red}}\label{sec:proof1}
\subsection{ Key steps in the proof.}
Let $\cF $ be a family of vector fields in $\T^2$, let $\hat\eta$ be a probability measure on $\cF$, let $\epsilon>0$, and let $f\in C^{2+\epsilon}(\T^2)$. Throughout this Section, assume $(\cF,\hat\eta)$ is \textit{good} in the sense of Definition \ref{def:good}, and that $f(p_0^{W})\neq \bar f^{W}$, $\hat\eta$-almost surely. Let $\rv W\sim\hat\eta$ and $\rv X\big|\rv W=W\sim \hat\rho_{W}$ a fixed family uniformly $C^2(\T)$.

Recall that we define, for each $x\in\T$, $W\in\cF$, and $t\geq 0$,
\[
I(W,x,t):=\int_0^{t} f \left(\phi_{s}^W\left(\Phi^W(x,0)\right)\right)ds - t\bar f^W,\qquad \bar f^W:=\frac{1}{\omega(\M^W)}\int_{\M^W}f(p)\,dp.
\]
and, given $N\geq 0$, $\alpha\in\R$, $x\in \T$, and $g\in L^2(\T)$, the Birkhoff sum of $g$
\be\label{eq:defBS}
    S(\alpha,g,x,N):=\begin{cases}
        0&\text{if }N=0,\\
        \sum_{m=0}^{N-1}g(R_\alpha^mx)&\text{if }N\geq 1.
    \end{cases}
\ee
The goal of this section is to prove Proposition \ref{prop:main_red}, which reduces the study of Birkhoff integrals of $f$ to that of the Birkhoff sums of the singular function $h$ in \eqref{eq:sing_p}. 

The proof proceeds in three stages. 

\subsection*{Step 1:} Define the centered {Birkhoff integral before a return}
\begin{equation}
    \label{eq:def_fhat}
    \hat{f}^W(x) := \int_0^{\ph^W(x)} f \left(\Phi^W(x,t)\right)dt-\ph^W(x)\bar f^W,
\end{equation}
where $\Phi^W(x,\cdot)$ denotes the orbit under $(\phi_t^W)$ of $\Phi^W(x,0)\in\Sigma^W$. Note $\int_\T\hat{f}^W(x)\,dx=0$.

We perform an asymptotic expansion of $\hat f$, separating its singular profile from a bounded-variation remainder. 

\begin{lemma}
    \label{lemma:fhat}
    For every $W\in\cF$, there exists $g^W\in L^\infty(\T)$ of bounded variation such that
    \begin{equation}
        \label{eq:fhat_expansion}
        \hat{f}^W(x) = g^W(x) + c_0^W\left(f(p_0^{W})-\bar f^W\right)h(x),
    \end{equation}
    where $p_0^W$ is the saddle, $c_0^W$ is the passage constant defined as \eqref{eq:passtime}, $h$ is the mean-zero, scaled, singular profile of the roof function defined as \eqref{eq:sing_p}, and $\bar f^W$ is the mean of $f$ on $\M^W$ defined as \eqref{eq:meanf}. 
\end{lemma}

Lemma \ref{lemma:fhat} is proven in Section \ref{sec:fhat}.

\subsection*{Step 2:} We establish precise probabilistic control over the fluctuations of some random hitting times $N(\rv W,\rv X,t)$ and present two Lemmas that simplify the computations in the final step.

Given a time $t>0$, $N\geq 0$, $W\in\cF$, and $x\in \T$, define
\[
    N(W,x,t):=\max\{k\geq 0 :\;t-S(\alpha(W),\ph^W,x,k)\geq 0\},\quad
    N_*(W,t):=\sup_x N(W,x,t).
\]

\begin{lemma}
    \label{cor:bounds} Suppose Proposition \ref{prop:main_red} holds under the following extra assumptions
    \be\gamma:=\inf_\cF\inf_\T\ph^W>0,\quad \mu:=\sup_\cF \mu^W<\infty,\label{eq:def_gamma}\ee
    \be V:=\sup_\cF\left(\text{Var}(g^W)+\text{Var}(\psi^W)+\norm{g^W}_\infty+\norm{\psi^W}_\infty\right)<\infty,\label{eq:def_V}\ee
    \be f_p:=\inf_\cF c_0^W\left|f(p_0^{W})-\bar f^W\right|>0,\quad c_0:=\sup_\cF c_0^W<\infty
    ,\label{eq:def_f0}\ee
    and, uniformly on $W$,
    \be\label{eq:boundsN}
    N_*(W,t)\asymp t.
    \ee
    Then the result holds in full generality as stated in Proposition \ref{prop:main_red}.
\end{lemma} 
Recall $\mu^W := \int_\T \varphi^W(x)\,dx$ and $\psi^W$, $g^W$ are the bounded-variation components of $\ph^W$ and $\hat{f}^W$, respectively. We also show that these bounded-variation remainders cannot contribute at the scale $(\log t)^2$.

\begin{lemma} \label{lemma:red3}
    Let $\rv W\sim\hat\eta$ and $\rv X\big|\rv W=W\sim \hat\rho_{W}$ any probability distribution for $\hat\eta-$almost every ${W\in\cF}$. Let $(g_0^W)_{W\in\cF}$ be a family of mean-zero functions with uniformly bounded variation in $\T$ (i.e. $\sup_\cF\text{Var}(g_0^W)<\infty$). It holds that 
    \[
    \frac{S\left(\alpha(\rv W),g_0^{\rv W},\rv X,k\right)}{(\log k)^2}\xrightarrow{\P } 0 \quad \text{as } k\to\infty.
    \]
\end{lemma}
Lemmas \ref{cor:bounds} and \ref{lemma:red3} are proven in Section \ref{sec:controlN}.

\subsection*{Step 3:} We transition from the continuous integral over $t$ to the Birkhoff sums of $\hat f$ up to $N(\rv W,\rv X,t)$, and then up to an $\rv X-$independent variable, showing at each step that the error terms generated by the reduction are asymptotically negligible under the $(\log t)^2$ scaling.

\begin{lemma} \label{lemma:red2}
    Under the hypotheses of Proposition \ref{prop:main_red} and \eqref{eq:def_gamma}--\eqref{eq:boundsN}, it holds that 
    \[
    \frac{|I(\rv W,\rv X,t)-S(\alpha(\rv W),\hat{f}^{\rv W},\rv X,\floor{t/\mu^{\rv W}})|}{c_0^{\rv W}\left|f(p_0^{\rv W})-\bar f^{\rv W}\right|(\log t)^2}\xrightarrow{\P } 0 \quad \text{as } t\to\infty.
\]
\end{lemma}
We prove Lemma \ref{lemma:red2} in Section \ref{sec:fits}. The proof builds on the results in Section \ref{sec:controlN}.

\subsection*{Proof of Proposition  \ref{prop:main_red}} 

\begin{proof}[Proof of Proposition  \ref{prop:main_red}]
Recall from Lemma \ref{lemma:fhat} that the induced function $\hat{f}^W$ can be written as
$
\hat{f}^W(x) = c_0^W\left(f(p_0^{W})-\bar f^W\right)h(x)+ g^W(x).
$
Taking the Birkhoff sum under $R_{\alpha(W)}$ of both sides up to $\floor{t/\mu^W}$ yields
\[
S(\alpha, \hat{f}^W, x, \floor{t/\mu^W}) = c_0^W\left(f(p_0^{W})-\bar f^W\right) S(\alpha, h, x, \floor{t/\mu^W}) + S(\alpha, g^W, x, \floor{t/\mu^W}).
\]
From Lemma \ref{cor:bounds}, we may assume \eqref{eq:def_gamma}--\eqref{eq:boundsN}. Thus, Lemma \ref{lemma:red2} applies and the continuous Birkhoff integral $I(W,x,t)$ can be approximated by this sum up to a negligible error in probability
\[
\frac{I(\rv W,\rv X,t)}{{c_0^{\rv W}}\left(f(p_0^{\rv W})-\bar f^{\rv W}\right)(\log t)^2} = \frac{S(\alpha(\rv W), \hat{f}^{\rv W}, \rv X, \floor{t/\mu^{\rv W}})}{{c_0^{\rv W}}\left(f(p_0^{\rv W})-\bar f^{\rv W}\right)(\log t)^2} + o_{\P }(1) \quad \text{as } t \to \infty.
\]
Combining these two expressions, we expand the normalized integral as
\[
\frac{I(\rv W,\rv X,t)}{c_0^{\rv W}\left(f(p_0^{\rv W})-\bar f^{\rv W}\right) (\log t)^2} = \frac{S(\alpha(\rv W), h, \rv X, \floor{t/\mu^{\rv W}})}{ (\log t)^2} + \frac{S(\alpha(\rv W), g^{\rv W}, \rv X, \floor{t/\mu^{\rv W}})}{c_0^{\rv W}\left(f(p_0^{\rv W})-\bar f^{\rv W}\right) (\log t)^2} + o_{\P }(1).
\]

Recall $\int_\T \hat{f}^W(x)\,dx=0$ so, for every $W\in\cF$, \[\int_\T g^W(x)\,dx=\int_\T \hat{f}^W(x)\, dx-c_0^W\left(f(p_0^{W})-\bar f^W\right)\int_\T h(x)\,dx=0.\]
Since $(g^W|f(p_0^{W})-\bar f^W|^{-1})_W $ are mean zero and, assuming \eqref{eq:def_V}--\eqref{eq:def_f0}, of uniformly bounded variation, Lemma \ref{lemma:red3} applies, giving
\[
\frac{S(\alpha(\rv W), g^{\rv W}, \rv X, \floor{t/\mu^{\rv W}})}{c_0^{\rv W}\left(f(p_0^{\rv W})-\bar f^{\rv W}\right) (\log t)^2} \xrightarrow{\P } 0 \quad \text{as }t\to\infty,
\]
where we have used the point-wise limit for each $W$: $\lim_{t\to\infty} \frac{\log \floor{t/\mu^{W}}}{\log t} = 1$. 

We apply the main hypothesis of Proposition \ref{prop:main_red} regarding the convergence to the random variable $Y$: \eqref{eq:convY}, to conclude 
\[
\frac{I(\rv W,\rv X,t)}{c_0^{\rv W}\left(f(p_0^{\rv W})-\bar f^{\rv W}\right) (\log t)^2} = \frac{S(\alpha(\rv W), h, \rv X, \floor{t/\mu^{\rv W}})}{(\log t)^2}+ o_{\P }(1)\ \xrightarrow[t\to\infty]{\text{dist}}\ Y
,\]
again by $\frac{\log\floor{t/\mu^W}}{\log  t}\to1$. This completes the proof of \eqref{eq:distrib_sum}. 
\end{proof}

\subsection{Asymptotic expansion of $\hat{f}^W$}\label{sec:fhat} The induced observable $\hat f^W$ measures the contribution of one complete excursion between successive returns to the transversal, so we expect $\hat f^W$ to inherit the same logarithmic singularity as the roof function $h$ in \eqref{eq:sing_p}. This is the content of Lemma \ref{lemma:fhat}.

\begin{proof}[Proof of Lemma \ref{lemma:fhat}] The proof is analogous to the proof of Lemma \ref{lemma:passtime}.

We fix $W$ and drop the notation dependence on $W$. Fix an open neighborhood $\cU_{0}$ of the saddle $p_0$ such that, for every $x\in\Sigma$, if the orbit segment $\{\Phi(x,t)\}_{t\in[0,\ph(x)]}$ meets $\cU_{0}$, it does so in at most two connected arcs (see \cite[$\S$3.1]{Ravotti2017}).

Write
\[
    \hat{f}(x) = \int_{\{t:\Phi(x,t)\notin\cU_0\}\cap[0,\ph(x)]} f \left(\Phi(x,t)\right)dt + \int_{\{t:\Phi(x,t)\in \cU_0\}\cap[0,\ph(x)]} f \left(\Phi(x,t)\right)dt-\ph(x)\bar f.
\]
Since $f$ is bounded, the flow is uniformly smooth away from the fixed point, and the time spent outside $\cU_0$ is uniformly bounded in $x$, the first integral defines a bounded-variation function of $x$, which we absorb into $g(x)$. It therefore suffices to analyze the second term, which we denote by $\hat F(x)$. 

We introduce local coordinates $(u,v)$ on $\cU_0$ in which $p_0=(0,0)$, the Hamiltonian is $H(u,v)=uv$, and the area form reads $\omega = V(u,v)du\wedge dv$ with $V$ smooth and $V> 0$ (\cite[$\S$3.1]{Ravotti2017}). In these coordinates the level set $H=c>0$ is the hyperbola $v=c/u$, and the flow satisfies $\dot{u}=u/V(u,c/u)$. Moreover, assume $\cU_0$ corresponds to $(-\eps,\eps)^2$ in these coordinates for some $\eps>0$.

Let $c=c(x)<\eps^2$ be the energy of the orbit of $\Phi(x,0)$. The integral along a maximal connected arc $\{\gamma(x,t)\}_{t\in[0,T]}$ of $\{t:\Phi(x,t)\in\cU_0\}\cap[0,\ph(x)]$ becomes
\begin{equation}
    \label{eq:local_integral}
    F(c) := \int_0^{T} f \left(\gamma(x,t)\right)dt = 
    { \int_{c/\eps}^{\eps}} f \left(u,{c}/{u}\right) \frac{V \left(u,c/u\right)}{u}du.
\end{equation}

By \cite[Lemma B.1]{FraczekUlcigrai2010}, applied to \eqref{eq:local_integral},
\begin{equation}
    \label{eq:Fexpansion}
    F(c) = f(0,0)V(0,0)|\log c| + \hat g(c),
\end{equation}
where $\hat g$ is in $L^\infty$ with bounded variation.
In the original coordinates $f(0,0)=f(p_0)$.

Note that, since $x$ parametrizes $\Sigma$ via the Hamiltonian, $c(x)=H_0\,\dist(x,\Z)$ for some $H_0> 0$. Hence, as $x\to 0^\pm$
\[
    |\log c(x)| = |\log \dist(x,\Z)| + O(1),
\]
so that, substituting into \eqref{eq:Fexpansion},
\[
   F(c(x)) = f(p_0)V(0,0)|\log \dist(x,\Z)| + O(1).
\]

By Lemma \ref{lemma:passtime}, the passage constant satisfies $c_0 = V(0,0) $, yielding the coefficient $2c_0$ for a double-passage orbit. Hence, after accounting for the number of passages by distinguishing the two sides, we get
\[
    \hat F(x) = 2c_0f(p_0)|\log x| + c_0f(p_0)|\log(1-x)| + O(1),
\]
where the $O(1)$ remainder lies in $L^\infty$ with bounded variation in $x$.

Finally, we need to consider the singular and bounded variation components coming from the term $-\ph(x)\bar f$.
Using \eqref{eq:def_phi} and collecting all bounded-variation remainders into a single function $g\in L^\infty(\T)$, we obtain expression \eqref{eq:fhat_expansion}.
\end{proof}

\subsection{Control on $N(\rv W,\rv X,t)$}\label{sec:controlN}
Given a time $t>0$, $N\geq 0$, $W\in\cF$, and $x\in \T$, define
\[
    N(W,x,t):=\max\{k\geq 0 :\;t-S(\alpha(W),\ph^W,x,k)\geq 0\},\quad
    N_*(W,t):=\sup_x N(W,x,t),
\]
where $\ph^W$ is the roof function given by \eqref{eq:def_phi} and $S(\alpha,g,x,k)$ is the Birkhoff sum of $g$ as defined in \eqref{eq:defBS}. 

Let $\norm{x}:=\dist(x,\Z)$,  $\gamma^W:=\inf_\T\ph^W>0$, and $\mu^W:=\int_\T\ph^W(x)\,dx \in[\gamma^W,\infty)$. Recall $\psi^W$ and $g^W$ are the bounded variation components of $\ph^W$ and $\hat f^W$, respectively. Also, recall $p_0^W$ is the saddle of each flow, and $c_0^W>0$ is the passage constant at $p_0^W$.

To establish Lemmas \ref{cor:bounds} and \ref{lemma:red2}, we need a deterministic control on the hitting time $N(W,x,t)$ and its uniform maximum $N_*(W,t)$. The following lemma is standard and shows $N_*$ grows linearly with respect to $t$.

\begin{lemma} \label{lemma:maxN} 
    For any $W\in\cF$ and $t>0$ it holds
    $\floor{t/\mu^W}\leq N_*(W,t)\leq \floor{t/\gamma^W}
    $.
\end{lemma}

\begin{proof}
    Fix $W$. The upper bound follows from the definition of $\gamma^W$: $S(\alpha,\ph^W,x,k)\geq\gamma^Wk$ for very $W,k$. For the lower bound, note that 
    \[\begin{aligned}
    \sup_x(t-S(\alpha,\ph^W,x,\floor{t/\mu^W}))\geq \int_\T (t-S(\alpha,\ph^W,x,\floor{t/\mu^W}))\,dx
    \\=t-\int_\T \sum_{m=0}^{\floor{t/\mu^W}-1}\ph^W(R_\alpha^mx)=(t-\mu^W\floor{t/\mu^W})\geq 0,\end{aligned}
    \]
    so $\exists x\in \T$ with $N(W,x,t)\geq\floor{t/\mu^W}$.
\end{proof}

Recall $(\cF,\hat\eta)$ is \textit{good} in the sense of Definition \ref{def:good}, $\eta$ is the distribution of its rotation parameter $\alpha(\rv W)$ and it is absolutely continuous with respect to Lebesgue$(\T)$, $f(p_0^{W})\neq \bar f^{W}$, $\hat\eta$-almost surely, and $(\hat \rho_W)_{W\in\cF}$ is a fixed family of probability distributions uniformly $ C^2(\T)$. 

Fix a probability space with the random variables $\rv W\sim \hat \eta$ and $\rv X\big|\rv W=W\sim\hat \rho_W$.

We present a purely probabilistic Lemma, from which Lemma \ref{cor:bounds} follows.

\begin{lemma}\label{lemma:bound}
    For each $n\geq 1$, let 
    \[
    \begin{aligned}\cF_n=\{&W\in\cF:\gamma^W>n^{-1},\quad c_0^W\left|f(p_0^{W})-\bar f^W\right|>n^{-1},\\&\quad c_0^W+\mu^W+\text{Var}(g^W)+\text{Var}(\psi^W)+\norm{g^W}_\infty+\norm{\psi^W}_\infty<n \},
    \end{aligned}
    \]
    and set $\rv W_n=\E_{\rv W}[\rv W\big|\cF_n]$, $\rv X_n|\rv W_n=W\sim\hat\rho_W$, whenever $\P_{\rv W}(\cF_n)>0$.
    
    Assume that there exists a random variable $Y$ such that, for $n$ large enough
    \[
    \frac{I(\rv W_n,\rv X_n,t)}{c_0^{\rv W_n}\left(f(p_0^{\rv W_n})-\bar f^{\rv W_n}\right)(\log t)^2}\ \xrightarrow[t\to\infty]{\text{dist}}\ Y,
    \]
    then 
    \[
    \frac{I(\rv W,\rv X,t)}{{c_0^{\rv W}}\left(f(p_0^{\rv W})-\bar f^{\rv W}\right)(\log t)^2}\ \xrightarrow[t\to\infty]{\text{dist}}\ Y.
    \]
\end{lemma}
\begin{proof}
    From the point-wise bounds for each $W$, we have $\P_{\rv W}(\cF_n)\to 1$ as $n\to\infty$. 
    
    Let $Z_N(W,x)=\frac{I(W,x,t)}{c_0^W\left(f(p_0^{W})-\bar f^W\right)(\log t)^2}$.
    Let $n$ be large enough so that $\P_{\rv W}(\cF_n)>0$. Then, for every bounded continuous function $g$,
    \[ 
    \E [g(Z_N(\rv W,\rv X))] =\P_{\rv W}(\cF_n)\E[g(Z_N(\rv W,\rv X)\big|\cF_n)]+\P_{\rv W}(\cF_n^c)\E [g(Z_N(\rv W,\rv X)\mid\cF_n^c)].
    \]
    Note that $\E[g(Z_N(\rv W,\rv X)\big|\cF_n)]=\E[g(Z_N(\rv W_n,\rv X_n)]$. Hence
    \[
    \left|\E [g(Z_N(\rv W,\rv X))]-\E[g(Z_N(\rv W_n,\rv X_n))]\right|\le 2\|g\|_\infty\,\P_{\rv W}(\cF_n^c).
    \]
    Taking $N\to\infty$ we get
    \[
    \limsup_{N\to\infty} \left|\E [g(Z_N(\rv W,\rv X))] -\E[g(Y)] \right| \le 2\|g\|_\infty\,\P_{\rv W}(\cF_n^c)\xrightarrow[N\to\infty]{}0.
    \]
    Therefore $\E [g(Z_N(\rv W,\rv X))]\xrightarrow[N\to\infty]{} \E[g(Y)],$ so $Z_N(\rv W,\rv X)\ \xrightarrow[N\to\infty]{\text{dist}}\ Y$, as claimed.
\end{proof}

\begin{proof}[Proof of Lemma \ref{cor:bounds}]
    Let $\cF_n$, $\rv W_n$ be as in Lemma \ref{lemma:bound} and denote by $\hat \eta_n$ the law of $\rv W_n$. Note that, provided $\P_{\rv W}(\cF_n)>0$, the distribution of $\alpha(\rv W_n)\in L^1$, and $(\hat\rho_W)_{W\in\cF_\eps}$ is uniformly $C^2(\T)$. Thus, $(\cF_n,\hat\eta_n)$ is \textit{good} in the sense of Definition \ref{def:good}, and satisfies both the hypothesis of Proposition \ref{prop:main_red} and \eqref{eq:def_gamma}--\eqref{eq:def_f0}. By Lemma \ref{lemma:bound}, showing the result under these extra conditions implies the general result in Proposition \ref{prop:main_red}.

    Combining Lemma \ref{lemma:maxN} and \eqref{eq:def_gamma} gives the bound \eqref{eq:boundsN} uniformly on $W$.
\end{proof}

To eventually bound the fluctuations of $N(\rv W,\rv X,t)$ and prove Lemma \ref{lemma:red2}, we need to show that the bounded-variation remainders of $\ph^W$ and $\hat{f}^W$ cannot contribute at the scale $(\log t)^2$, namely Lemma \ref{lemma:red3}.
To do so, we use the arithmetic rigidity of irrational rotations.

\begin{proof}[Proof of Lemma \ref{lemma:red3}]
For each $W$ with $\alpha=\alpha(W)\in(0,1)\setminus\Q$, we express $k$ in its Ostrowski representation with respect to the continued fraction expansion of $\alpha=[0;a_1,a_2,\dots]$ (see e.g. \cite{khinchin1964})
\[
k = \sum_{j=0}^{n} m_j q_j, \quad \text{where } 0 \leq m_j \leq a_{j+1},\quad q_{n+1}=a_{n+1}q_n+q_{n-1},\quad q_0=1.
\]
Recall that from the standard properties of continued fractions, $n = O(\log k)$.

Since $\int_\mathbb{T} g_0^W \,dx = 0$, we can directly apply the standard Denjoy--Koksma inequality \cite{Herman1979} to each denominator block, bounding $|S(\alpha, g_0^W, x, q_j)|$ uniformly by $\mathrm{Var}(g_0^W)$. Combining this with the Ostrowski representation via the triangle inequality gives, uniformly on $x\in\T$:
\[
|S(\alpha, g_0^W, x, k)| \leq \mathrm{Var}(g_0^W) \sum_{j=0}^{n} m_j \leq \sup_\cF\mathrm{Var}(g_0^W) \sum_{j=0}^{n} a_{j+1}.
\]

To bound the sum of the partial quotients, we use a classical metric result by Khinchin \cite{khinchin1964}. Under the uniform Lebesgue measure on $(0,1)$, the sum of the partial quotients satisfies a weak law of large numbers
\[
\frac{1}{n \log n} \sum_{j=1}^{n} a_j =O_{\P_{\text{Leb}}}(1) \quad \text{as } n \to \infty.
\]
In terms of asymptotic growth, this and $\rv A=\alpha(\rv W)\sim\eta$ absolutely continuous with respect to Lebesgue$(\T)$ implies that $\sum_{j=0}^{n} a_{j+1} = O_{\P_{\rv A}}(n \log n)$, and consequently $ O_{\P_{\rv W}}(n \log n)$. Substituting our step bound $n = O(\log k)$, we find that
\[
|S(\alpha, g_0^{\rv W}, \rv X, k)|\leq\sup_{\cF}\mathrm{Var}(g_0^W)\sum_{j=0}^{n} a_{j+1} = O_{\P }(\log k \log\log k).
\]

Finally, normalizing by $(\log k)^2$ yields
\[
\frac{|S\left(\alpha(\rv W),g_0^{\rv W},\rv X,k\right)|}{(\log k)^2} = O_{\P }\left(\frac{\log k \log\log k}{(\log k)^2}\right) = O_{\P }\left(\frac{\log\log k}{\log k}\right)\xrightarrow[k\to\infty]{\P }0.
\]
\end{proof}

Lemma \ref{lemma:red3} is a key ingredient for the next Lemma. 

For each $W$, recall $\mu^W=\int_\T\ph^W(x)\,dx\geq\gamma^W=\inf_\T\ph^W>0$ for every $W\in\cF$.

\begin{lemma} \label{lemma:logDN}
    Let $\rv W\sim\hat\eta$ in $\cF$ and $\rv X\big|\rv W=W\sim \hat\rho_{W}$ uniformly $C^2(\T)$ and suppose there exists a continuous random variable $Y$ such that
    \[
    \frac{S\left(\alpha(\rv W),h,\rv X,N\right)}{(\log {N})^2}\ \xrightarrow[N\to\infty]{\text{dist}}\ Y.
    \]
    Moreover, assume $\cF$ is such that \eqref{eq:def_gamma}--\eqref{eq:def_f0} hold. Then for any $r\in\R$
    \[\begin{aligned}
    \lim_{t\to\infty}&
    \P \left(\frac{\floor{t/\mu^{\rv W}}-N(\rv W,\rv X,t)}{c_0^{\rv W}(\log t)^2}\leq \frac{r}{\mu^{\rv W}}\right)=\P\left(Y\leq {r}\right).
    \end{aligned}
    \]
\end{lemma}
\begin{proof}
For simplicity, we drop the dependence on $\rv W$ on the notation. First, we claim that we can assume
\[
    \lim_{k\to \infty}\P \left(\frac{S\left(\alpha(\rv W),\ph ,\rv X,k\right)-{k}\mu }{c_0(\log {k})^2}\leq r\right)=\P\left(Y\leq r\right).
    \] 
Indeed, by \eqref{eq:def_V}--\eqref{eq:def_f0} and the boundedness of $f$, the correction terms \[g_0^W(x)=g_0(x):=c_0^{-1}\ph (x)-c_0^{-1}\mu -h(x)=c_0^{-1}\psi (x)+3-c_0^{-1}\mu,\] satisfy the hypotheses of Lemma \ref{lemma:red3}, so
\[
\frac{S(\alpha(\rv W), g_0,\rv X, k)}{(\log k)^2} \xrightarrow{\P } 0 \quad \text{as }k \to \infty.
\]
Hence
\[
\frac{S\left(\alpha(\rv W),\ph-\mu,\rv X,k\right)}{c_0(\log {k})^2}=\frac{S\left(\alpha(\rv W),h,\rv X,k\right)}{(\log {k})^2}+\frac{S\left(\alpha(\rv W),g_0,\rv X,k\right)}{(\log {k})^2}\ \xrightarrow[k\to\infty]{\text{dist}}\ Y.
\]

Now, for each $r\in \R$,
\[
\begin{aligned}
\P &\left(\frac{\floor{t/\mu }-N(\rv W,\rv X,t)}{c_0(\log t)^2}\leq \frac{r}{\mu }\right)
=\P \left({\floor{\mu \floor{t/\mu }-r{c_0(\log t)^2}}\leq \mu N(\rv W,\rv X,t)} \right)
\\&=\P \left(t\geq S\left(\alpha(\rv W),\ph ,\rv X,\floor{\floor{t/\mu }-rc_0{(\log t)^2}/\mu }\right)\right)\\
&= \P \Bigg( 
    \frac{\floor{t} - \floor{\floor{t}-rc_0(\log t)^2}}{\big(\log \floor{\floor{t}-rc_0(\log t)^2}\big)^2} 
     \\ &\hspace{2cm} \geq 
    \frac{S\big(\alpha(\rv W),\ph ,\rv X,\floor{\floor{t/\mu }-rc_0(\log t)^2/\mu }\big) - \floor{\floor{t}-rc_0(\log t)^2}}{\big(\log \floor{\floor{t}-rc_0(\log t)^2}\big)^2} 
\Bigg).
\end{aligned}
\]
Now, taking limit $t\to\infty$ in the equalities above and letting $k(t,r):=\floor{\floor{t}-rc_0{(\log t)^2}}/\mu $ and $k'(t,r):=\floor{\floor{t/\mu }-rc_0{(\log t)^2/\mu }}$, we have
\[\begin{aligned}
&\limsup_{t\to\infty}
\P \left(\frac{\floor{t/\mu }-N(\rv W,\rv X,t)}{c_0(\log t)^2}\leq \frac{r}{\mu }\right)\\
&\qquad =\limsup_{t\to \infty}\P \left(\frac{S\left(\alpha(\rv W),\ph ,\rv X,k'(t,r)\right)-{\mu k(t,r)}}{(\log {k(t,r)}+\log \mu )^2}\leq \frac{\floor{t}-\floor{\floor{t}-rc_0{(\log t)^2}}}{(\log \floor{\floor{t}-rc_0{(\log t)^2}})^2}\right),
\end{aligned}
\]
and similarly for $\liminf_{t\to\infty}$. 

Note the following point-wise limits as $t\to\infty$: for any $r\in\R$ and $W\in \cF$
\[
k(t,r),k'(t,r)\to\infty,\qquad \frac{k(t,r)}{k'(t,r)}\to 1, \qquad \frac{\floor{t}-\floor{\floor{t}-rc_0{(\log t)^2}}}{(\log \floor{\floor{t}-rc_0{(\log t)^2}})^2}\to r.
\]
Also, by \eqref{eq:def_gamma} and $\mu^W\geq \gamma^W$,  $|k(t,r)-k'(t,r)|$ is uniformly bounded on $t$, $r$ and $W$. Thus, we have for every $\eps>0$ that
\[\begin{aligned}
\limsup_{t\to\infty}
\P &\left(\frac{\floor{t/\mu }-N(\rv W,\rv X,t)}{c_0(\log t)^2}\leq \frac{r}{\mu }\right)\\&\leq\limsup_{k\to \infty}\P \left(\frac{S\left(\alpha(\rv W),\ph ,\rv X,k\right)-{\mu k}}{c_0(\log {k})^2}\leq r+\eps \right)=\P\left(Y\leq {r+\eps}\right),
\end{aligned}\]
and similarly
\[\begin{aligned}
\liminf_{t\to\infty}
\P &\left(\frac{\floor{t/\mu }-N(\rv W,\rv X,t)}{c_0(\log t)^2}\leq \frac{r}{\mu }\right)\geq\P\left(Y\leq {r-\eps}\right),
\end{aligned}\]
so taking $\eps\to 0$, the continuity of $Y$ yields \[\begin{aligned}
\P\left(Y\leq {r}\right)&\leq \liminf_{t\to\infty}
\P \left(\frac{\floor{t/\mu }-N(\rv W,\rv X,t)}{c_0(\log t)^2}\leq \frac{r}{\mu }\right)\\&\leq \limsup_{t\to\infty}
\P \left(\frac{\floor{t/\mu }-N(\rv W,\rv X,t)}{c_0(\log t)^2}\leq \frac{r}{\mu }\right)\leq \P\left(Y\leq {r}\right).\end{aligned}
\]
\end{proof}

\subsection{From Birkhoff integrals to Birkhoff sums} \label{sec:fits}
With these linear and distributional controls on the hitting times in hand, we can begin the process of discretizing the Birkhoff integral. The first reduction step consists of replacing $I(W,x,t)$ by the Birkhoff sum of the induced function $\hat{f}^W$ evaluated up to $N(W,x,t)$. Then, we will complete the proof of Lemma \ref{lemma:red2} by showing that the discrepancy between the Birkhoff sums evaluated at $N(\rv W,\rv X,t)$ and the $\rv X$-independent random variable $\floor{t/\mu^{\rv W}}$ is asymptotically negligible, up to our scaling factor.
The following lemma establishes this initial substitution.

\begin{lemma}\label{lemma:red1}
{ Let $\alpha(\mathcal{W})$ denote the rotation number corresponding to the vector field $\mathcal{W}.$}
Assuming \eqref{eq:def_gamma}--\eqref{eq:boundsN}, it holds that
\[
    \frac{|I(\rv W,\rv X,t)-S(\alpha(\rv W),\hat{f}^{\rv W},\rv X,N(\rv W,\rv X,t))|}{c_0^{\rv W}\left|f(p_0^{\rv W})-\bar f^{\rv W}\right|(\log t)^2}\xrightarrow{\P } 0 \quad \text{as } t\to\infty.
\]
\end{lemma}

\begin{proof}
We will show that there exists { a constant} $C=C(\cF,f)>0$ and for every $W\in \cF$  { there exists a set} 
$ E_1(W,t)\subset\T$ 
such that for $t$ large enough, and $x\in\T\setminus E_1(W,t)^c$ 
    \[
       \frac{|I(W,x,t)-S({\alpha(\mathcal{W})},\hat{f}^W,x,N(W,x,t))|}{c_0^{W}|f(p_0^{W})-\bar f^W|}\leq C \log t,
    \]
    and
    \[
         \lim_{t\to\infty}\P_{\rv X}(E_1(W,t)|\rv W=W)=0 \quad\text{uniformly on } W.
    \]
    The result immediately follows.
    
Let $N_*=N_*(W,t)$ and, for each $W$, let 
\[
    E_1(W,t):=\bigcup_{m=0}^{N_*} R_\alpha^{-m}\left(-1/N_*^2,1/N_*^2\right).
\]
By \eqref{eq:boundsN} and the $C^2$ uniformity of $(\hat\rho_W)$ it holds, uniformly on $W$,
\[
    \P_{\rv X}(E_1(W,t)|\rv W=W)\leq 2\norm{\hat\rho_W}_{\infty}(N_*+1)/N_*^2\lesssim 1/t\xrightarrow{t\to \infty}0.
\] 

Fix $W$ and set $\alpha=\alpha(W)$. We have \[
    \begin{aligned}
    I(W&,f,x,t)= \int_0^{t} f \left(\phi_s^W\Phi^W(x,0)\right)ds-t\bar f^W\\
   & = S(\alpha,\hat{f}^W,x,N(W,x,t))+\int_{S(\alpha,\ph^W,x,N(W,x,t))}^{t} \left(f \left(\phi_s^W\Phi^W(x,0)\right)-\bar f^W\right)ds.
    \end{aligned}
\]
Note that we can bound the last term using \eqref{eq:def_phi} as follows
\[
    \begin{aligned}\left|\int_{S(\alpha,\ph^W,x,N(W,x,t))}^{t}  \left(f \left(\phi_s^W\Phi^W(x,0)\right)-\bar f^W\right)ds\right|\leq2\norm{f}_\infty|t-S(\alpha,\ph^W,x,N(W,x,t))|\\\leq 2\norm{f}_\infty \ph^W(R_\alpha^{N(W,x,t)}x)\leq2\norm{f}_\infty\left(-3c_0\log\|R_\alpha^{N(W,x,t)}x\|+\norm{\psi^W}_\infty\right),
    \end{aligned}
\]
where we have used $ |\bar f^W|\leq \norm{f}_\infty.$

From \eqref{eq:def_V}, $\norm{\psi^W}_\infty\leq V$, from \eqref{eq:def_f0}, $c_0^{W}\left|f(p_0^{W})-\bar f^W\right|\geq f_p>0$, and from \eqref{eq:boundsN}, $N_*\asymp t$ uniformly on $W$. Therefore there exists $C=C(\cF,f)>0$ such that for $t$ large enough, for any $W\in \cF$, and any $x\in E_1(W,t)^c$
\[
    \frac{|I(W,x,t)-S(\alpha,\hat{f}^W,x,N(W,x,t))|}{c_0^{W}\left|f(p_0^{W})-\bar f^W\right|}\leq C  \log t,
\]
where we have used $N(W,x,t)\leq N_*$, $\forall x$, so $\|R_\alpha^{N(W,x,t)}x\|\geq N_*^{-2}$, $\forall x\in E_1(W,t)^c$.
\end{proof}

While Lemma \ref{lemma:red1} reduces the continuous integral to a discrete sum, the stopping time $N(\rv W,\rv X,t)$ remains a random variable that depends heavily on the geometry of each roof function $\ph^W$. To replace $N(\rv W,\rv X,t)$ with the $\rv X-$independent random variable $\floor{t/\mu^{\rv W}}$, we use the distribution of its fluctuations obtained in Lemma \ref{lemma:logDN}. 
\begin{proof}[Proof of Lemma \ref{lemma:red2}]
We drop the notation dependence on $\rv W$. From Lemma \ref{lemma:red1}, it is enough to show 
\[
    \frac{|S(\alpha(\rv W),\hat{f} ,\rv X,N(\rv W,\rv X,t))-S(\alpha(\rv W),\hat{f} ,\rv X,\floor{t/\mu })|}{c_0\left|f(p_0)-\bar f\right|(\log t)^2}\xrightarrow{\P } 0 \quad \text{as } t\to\infty.
\]

Note that for every $W\in\cF$, $\alpha=\alpha(W)$, and $x\in\T$ 
\[\begin{aligned}
|S(\alpha,\hat{f} ,x&,N(W,x,t))-S(\alpha,\hat{f} ,x,\floor{t/\mu })|\\&=
\begin{cases}
    |S(\alpha,\hat{f} ,R_\alpha^{\floor{t/\mu}}x,N(W,x,t)-\floor{t/\mu })| &\text{if } N(W,x,t)\geq\floor{t/\mu },\\
    |S(\alpha,\hat{f} ,R_\alpha^{N(W,x,t)}x,\floor{t/\mu }-N(W,x,t))| &\text{if } N(W,x,t)<\floor{t/\mu }.
\end{cases}
\end{aligned}
\]

Similarly to the proof of Lemma \ref{lemma:red1}, set $N_*=N_*(W,t)$ and let, for $t>1$, 
\[
    E_1(W,t)=\bigcup_{m=0}^{N_*} R_\alpha^{-m}\left(-1/N_*^2,1/N_*^2\right).
\]
Recall $\lim_{t\to\infty}\P_{\rv X}(E_1(W,t)|\rv W=W)=0$ uniformly on $W$.
Moreover, define \[
E_2(t)=\left\{(W,x):{|\floor{t/\mu }-N(W,x,t)|}> (\log t)^3\right\},
\]
It follows from Lemma \ref{lemma:logDN}, \eqref{eq:def_gamma}, \eqref{eq:def_f0}, and $\mu^W\geq\gamma^W$ that $\lim_{t\to\infty}\P (E_2(t))=0$.

Consider the family of functions given by $\hat g^W =\hat g= \hat{f} \chi_{\{\norm{x}\geq N_*^{-2}\}}$. Note that $\hat g $ does not have exactly zero mean; instead, $\int_\T \hat g  \,dx = -\int_{\{\norm{x}<   N_*^{-2}\}} \hat{f} \,dx$. By the expression of $\hat{f}$ in Lemma \ref{lemma:fhat} and \eqref{eq:def_gamma}--\eqref{eq:def_V}, there exist constants $C,C'>0$ depending only on $\cF$ and $f$ such that for $t$ large enough and every $W\in \cF$ the mean satisfies \[\left|\int_\T \hat g  \,dx\right| \leq C\frac{\log N_*}{N_*^2}+\frac{2\norm{g^W}_\infty}{N_*^2}\leq C'\frac{\log t}{t^2}.\]

For almost-every $W$, $\alpha(W)\in\RQ$. Let $(q_n)$ be the sequence of best denominators of $\alpha$. By the Denjoy--Koksma inequality, for all $n$ and all $x \in \T$,
\[
\left|S(\alpha, \hat g , x, q_n) - q_n \int_\T \hat g  \,dx \right| \leq \mathrm{Var}(\hat g ).
\]
Moreover, from \eqref{eq:def_V} and \eqref{eq:boundsN}, there exists $\hat C,\hat C'>0$ depending only on $\cF$ and $f$ such that \[
\mathrm{Var}(\hat g )\leq 2c_0\norm{f}_\infty\text{Var}\left(h\chi_{\{\norm{x}\geq N_*^{-2}\}}\right)+2\norm{g}_\infty+\text{Var}(g) \leq \hat C\log N_* +3V\leq \hat C'\log t,\]
for $t$ large enough, yielding
\[
|S(\alpha, \hat g , x, q_n)| \leq \mathrm{Var}(\hat g ) + q_n \left|\int_\T \hat g  \,dx \right| \leq \hat C'\log t + q_n C'\frac{\log t}{t^2}.
\]

Now, assume $(W, x) \in E_2(t)^c$. Writing $\Delta N(W,x,t) = |\floor{t/\mu }-N(W,x,t)|$ we have $\Delta N \leq (\log t)^3$, so for any block index $j$ appearing in its Ostrowski expansion, $q_j \leq (\log t)^3$. Thus, the mean term satisfies $q_j \left|\int_\T \hat g \,dx \right| \leq C'(\log t)^4/t^2\to 0$, and $|S(\alpha, \hat g , x, q_j)| \leq 2\hat C' \log t$ for all blocks, for $t$ large enough.

Writing $\Delta N = \sum_{j=0}^{n} m_j q_j$ and applying the block-wise Denjoy--Koksma bound gives, for any $(W,x)\in E_2(t)^c$ with $x\in E_1(W,t)^c$,
\[
|S(\alpha, \hat{f} , x, N(W,x,t))-S(\alpha,\hat{f} ,x,\floor{t/\mu})| \leq \sup_{x \in E_1^c} |S(\alpha, \hat g , x, \Delta N)| \leq 2\hat C' \log t \sum_{j=0}^{n} m_j,
\]
where we have used that the $N_*-$orbit of $x\in E_1(W,t)^c$ under $R_\alpha$ does not visit $(-N_*^{-2}N_*^{-2})$.

Finally, note that, by \eqref{eq:def_f0}, $c_0\left|f(p_0)-\bar f\right|\geq f_p>0$ uniformly on $W$ and, using the same classical metric result as in the proof of Lemma \ref{lemma:red3}, \[\sum_{j=0}^{n} m_j=O_{\P_{\rv A}}(n\log n)=O_{\P_{\rv W}}(\log \Delta N\log\log \Delta N),\] with $\Delta N\leq (\log t)^3$ so $\sum_{j=0}^{n} m_j=o_{\P_{\rv W}}((\log\log t)^2)$, and thus
\[
\frac{|S(\alpha(\rv W), \hat{f} , \rv X, N(\rv W,\rv X,t))-S(\alpha(\rv W),\hat{f} ,\rv X,\floor{t/\mu})|}{c_0\left|f(p_0)-\bar f\right|}\chi_{E(t)}(\rv W,\rv X)=o_{\P }((\log t)^2).
\]
for $E(t)=\{(W,x)\in E_2(t)^c:\ x\in E_1(W,t)^c\}$. The result follows from $\P (E(t))\xrightarrow {t\to\infty}1$.
\end{proof}

\section{Proof of Proposition \ref{prop:main_harm}}\label{sec:proof2}
At this point, the problem has been reduced to understanding Birkhoff sums
\[
    S(\alpha,h,x,N)=\sum_{m=0}^{N-1}h(R_\alpha^mx)=\sum_{m=0}^{N-1}h(x+m\alpha\bmod1),
\]
of the explicit logarithmic profile $h=2|\log x|+|\log(1-x)|-3$, for $\alpha$ and $x$ varying according to some distribution.

The philosophy, following Dolgopyat--Sarig \cite{DolgopyatSarig2019} or Dolgopyat--Fayad \cite{DolgopyatFayad2020}, is that only a sparse collection of resonant Fourier modes contributes at the scale $(\log N)^2$, while all remaining frequencies produce negligible fluctuations.

We first introduce some probability preliminaries in $\S$\ref{sec:prelim}. Then, the proof of Proposition \ref{prop:main_harm} ($\S$\ref{sec:s5}) proceeds in four steps. Step 1 ($\S$\ref{sec:s1}) identifies the Fourier coefficients of $h$. Steps 2--3 ($\S$\ref{sec:s2}--$\S$\ref{sec:s3}) reduce the ergodic sum, up to negligible errors, to a finite random series indexed by the prime resonant harmonics (see Definition \ref{def}). And Step 4 ($\S$\ref{sec:s4}) establishes the joint limiting law of the data attached to the resonant harmonics of each irrational $\alpha$. Assembling these inputs and the characterizations in $\S$\ref{sec:prelim} identifies the limit as a Cauchy distribution.

\subsection{Probability preliminaries}\label{sec:prelim}

Recall (see e.g. \cite[$\S$2.1]{DolgopyatFayad2020} or \cite[Appendix B]{DolgopyatSarig2019}) the following characterization of the Cauchy random variable.
\begin{lemma}(See e.g. \cite[Lemma 2.2]{DolgopyatFayad2020})\label{lemma:cauchy}
Let $\Xi = \{\xi_n\}$ be a Poisson process on $\mathbb{R}$ with intensity $c$ and let
\[ Y_R = \sum_{|\xi_n| < R} \frac{\Phi_n}{\xi_n} \]
where $\Phi_n$ are i.i.d. bounded random variables with zero mean independent of $\Xi$. Define
\[
\hat{c} = c\pi \E(|\Phi|),
\]
and denote by $\mathfrak{C}$ the standard Cauchy distribution with density $\frac{1}{\pi(1+x^2)}$, and by $\mathfrak{C}_{R}$ the law of $Y_R/{\hat c}$. Then 
\[
\mathfrak{C}_{R}\ \xrightarrow[R \to \infty]{\text{dist}}\ \,\mathfrak{C}.
\] 
\end{lemma}

For $x,y\in\T$, define the function
\begin{equation}
  \label{DefTheta} 
\Theta(x,y):=\sum_{k\geq 1}\frac{\sin(2\pi k x) \sin(2\pi k y)}{\pi^2k^2}.
\end{equation}
\begin{lemma} \label{lemma:rv}
The { series \eqref{DefTheta}}  converges at every point to a bounded function. Moreover, letting $X,Y$ denote two independent uniform random variables on $\T$, we have
\be\label{eq:mean_theta}
\E_{X,Y}(\Theta(X,Y))=0,\quad \E_{X,Y}(|\Theta(X,Y)|)= \tfrac{1}{24}.
\ee
\end{lemma}
\begin{proof}
Convergence and boundedness of $\Theta$ follows from summability of $k^{-2}$.

Now, for each $k\geq 1$, it holds $\E_{X}[\sin(2\pi k X)]=0$, so by linearity of the expectation, $\E_{X,Y}(\Theta(X,Y))=0$.

We denote by $\norm x=\dist(x,\Z)$. We have $\sum_{k\geq 1}\frac{\cos(2\pi kx)}{2\pi^2k^2}=\frac{x(x-1)}{2}+\frac{1}{12}$, for $x\in [0,1)$, thus
\[\begin{aligned}
\Theta(x,y)&:=\sum_{k\geq 1}\frac{\sin(2\pi k x) \sin(2\pi k y)}{\pi^2k^2}= \sum_{k\geq 1} \frac{\cos(2\pi k (x-y))}{2\pi^2k^2}-\sum_{k\geq 1} \frac{\cos(2\pi k (x+y))}{2\pi^2k^2},
\\&\ =\frac{\norm{x+y}(1-\norm{x+y})}{2}-\frac{\norm{x-y}(1-\norm{x-y})}{2}.
\end{aligned}
\]
Noting that $X-Y \ (\bmod 1)$ and $X+Y\ (\bmod 1)$ are also independent uniform r.v. on $\T$, the second equality in \eqref{eq:mean_theta} follows from a standard computation.
\end{proof}

\subsection{Main steps of the proof of Proposition \ref{prop:main_harm}} \label{sec:s5} Let us state the main steps of the proof and show how they yield the result.

We first present some definitions as stated in \cite[\S 3]{DolgopyatSarig2019}.

We denote by $\snorm\cdot:\R\to(-\tfrac{1}{2},\tfrac{1}{2}]$ the signed distance to the nearest integer and by $\norm{\cdot}=|\snorm{\cdot}|$. Fix $N$ and $\alpha\in \RQ$. An integer $1\leq k\leq N$ is called a \textit{prime harmonic}, if $k\alpha = \snorm{k\alpha}+ m$ where $\gcd(k, m) = 1$. If $k$ and $m$ are not coprime, that is $r := \gcd(k, m)$, then we call $k/r$ the \textit{prime harmonic associated to $k$}.

\begin{definition} \label{def}
Fix $\delta > 0$, and $N \gg 1$. Let $T = N (\log N)^2$.

(1) $p \in \mathbb{N}$ is called a \emph{prime resonant harmonic}, if $p \leq N$, $p$ is a prime harmonic, and
\[
p\|p\alpha\| \leq (\delta \log N)^{-1}.
\]

(2) $k \in \mathbb{N}$ is called a \emph{resonant harmonic}, if $k \leq T$, and the prime harmonic associated to $k$ is a prime resonant harmonic.
    
\end{definition} 

Let $\mathcal{R}(\delta, N,\alpha)$ denote the set of resonant harmonics, $\mathcal{P}(\delta, N,\alpha)$ the set of prime resonant harmonics, and $\mathcal{O}(\delta, N,\alpha)$ the set of non-resonant harmonics less than $N$. 

\subsection*{Steps 1 and 2:} Let $K\subseteq\Z^+$ and define, for $\alpha\in\RQ$, $x\in\T$, and $N\geq 1$
\[
\mathcal{G}(\alpha,x,N,K):=\sum_{k\in K}\frac{\log k}{k}\frac{\sin(\pi N k \alpha)f_{N,k}(x,\alpha)}{\pi\sin(\pi k\alpha)},
\]
where, for each $k\geq 1$, we define
\[
f_{N,k}(x,\alpha):=\sin\big(2\pi kx+\pi (N-1)k\alpha\big).
\]

We identify the Fourier coefficients of $S(\alpha,h,x,N)$ (Step 1) and then reduce the problem to the study of $\mathcal{G}(\alpha,x,N,\mathcal{R})$ (Step 2). The asymmetry of $h$ is crucial for this step.

Fix a probability space with the random variables $\rv A\sim \eta$ and $\rv X\big|\rv A=\alpha\sim \rho_\alpha$. 

\begin{proposition} \label{prop:harmonics} For any $\epsilon>0$ $\exists\delta_0>0$ such that, for any $\delta_0>\delta>0$ and $N$ large enough,
\be
    \P \left(\frac{|S(\rv A,h,\rv X,N)-\mathcal{G}(\rv A,\rv X,N,\mathcal{R}(\delta,N,\rv A)) |}{(\log N)^2}>\eps\right)<\eps.
\ee
\end{proposition}
Proposition \ref{prop:harmonics} is proven in Sections \ref{sec:s1}--\ref{sec:s2}.

\subsection*{Step 3:} We decompose the sum $\mathcal{G}(\alpha,x,N,\mathcal{R}(\delta,N,\alpha))$ prime-by-prime.

Let
\[
 \Theta(x,y):=\sum_{k\geq 1}\frac{\sin(2\pi kx)\sin(2\pi k y)}{\pi^2k^2},
\]
which we know is well-defined for any $x,y$ by Lemma \ref{lemma:rv}, and let 
\be
\nu_p:=\frac{\log p}{\log N},\qquad
\xi_p:=p\ssnorm{p\alpha}\log N, \qquad y_p:=Np\alpha\bmod 2.
\label{eq:def_rv1}\ee
\be
\Theta_p(x,\alpha):=\begin{cases} \Theta \left(\frac{p(N-1)\alpha+2px\bmod 2}{2},\frac{y_p}{2}\right),& p\alpha-((p\alpha))\in 2\Z,\\
\Theta \left(\frac{p(N-1)\alpha+2px\bmod 2}{2},\frac{1-y_p}{2}\right),& p\alpha-((p\alpha))\in 2\Z+1.
\end{cases}
\label{eq:def_rv2}\ee
\begin{lemma} \label{lemma:sumrv}
For any $\alpha\in\RQ$, $x\in\T$, and $\delta>0$ sufficiently small,
\be\begin{aligned}\label{eq:asymp}
\frac{\mathcal{G}(\alpha,x,N,\mathcal{R}(\delta,N,\alpha))}{(\log N)^2}
&=\sum_{p\in\mathcal{P}(\delta,N,\alpha)}\frac{1}{\xi_p} \Bigl(\nu_p\Theta_p(x,\alpha)+O \bigl(\tfrac{1}{\log N}\bigr)\Bigr)
,
\end{aligned}
\ee
where the constant in the $O(\cdot)$ term depends on $\delta$ but is uniform in $\alpha$, $x$ and $p$.
\end{lemma}

Lemma \ref{lemma:sumrv} is proven in Section \ref{sec:s3}. The proof is analogous to the proof of \cite[Proposition 3.5]{DolgopyatSarig2019}.

\subsection*{Step 4:} The sum \eqref{eq:asymp} is now a finite series of the form $\sum_p \frac{\nu_p}{\xi_p}\Theta_p + \text{lower order}$. The remaining task is to identify the joint law of the triple $(\nu_p,\xi_p,\Theta_p)$ as $N\to\infty$.

We use the following result from \cite{DolgopyatFayad2020}.

\begin{proposition}(\cite[Theorem 5.1]{DolgopyatFayad2020})\label{prop:PPP}
Suppose that $\rv A$ is distributed according to a smooth density on $\mathbb{T}$. For each $\delta$ the point process
\[
\left\{
\left(
\nu_p,\ \xi_p,\ y_p 
\right)\right\}_{p \in \mathcal{P}(\delta,N,\rv A)}
\]
converges in distribution as $N \to \infty$ to the Poisson Point Process on $[0,1] \times \left[-\frac{1}{\delta}, \frac{1}{\delta}\right]\times[0,2], $ with constant intensity $c={6}/{\pi^2}.$
\end{proposition}

We use the randomness of $x$ to find the limiting distribution of the vectors \[\bigl(p_1[\alpha(N-1)+2x],\dots, p_r[\alpha(N-1)+2x]\bigr)\bmod 2,\] for $\{p_1,\dots,p_r\}\in\cP(\delta,N,\alpha).$

\begin{proposition} \label{prop:equidx}
Suppose that $\rv A$ is distributed according to a bounded density on $\mathbb{T}$. For every $r>1$ and $\varepsilon>0$, if $\delta$ is small enough there are $N_0$ and $ A(N,\delta,r)\subset \mathbb{T}$ such that $ \P_{\rv A}(A(N,\delta,r))>1-\eps$ and
\begin{enumerate}
\item If $\alpha\in A(N,\delta,r)$ then $|\mathcal{P}(\delta,N,\alpha)|\geq r$ for all $N\geq N_0$.
\item For each neighborhood $\mathcal{V}$ of the uniform distribution on $[0,2]^{r}$ there exists $N_\mathcal{V}$ such that for $N\geq N_\mathcal{V}$ the following holds. Let $\alpha\in A(N,\delta,r)$ and let $p_1,\ldots,p_r$ be an enumeration of the prime resonant harmonics in $\mathcal P(\delta,N,\alpha)$ so that $\norm{p_k\alpha}\, p_k\leq\norm{p_{k'}\alpha}\, p_{k'}$ for $k\leq k'$. Then the distribution of the random vectors
\[
\bigl(p_1[\alpha(N-1)+2\mathcal{X}_\alpha],\dots, p_r[\alpha(N-1)+2\mathcal{X}_\alpha]\bigr)\bmod 2,
\]
where $\rv X_\alpha\sim \rho_\alpha$ for some $(\rho_\alpha)_{\alpha}$ uniformly $C^{2}(\T)$, belongs to $\mathcal{V}$.
\end{enumerate}
\end{proposition}

Proposition \ref{prop:PPP} gives the joint law of $(\nu_p,\xi_p,y_p)$, and Proposition \ref{prop:equidx} complements this by showing that, as $\rv X\big|\rv A=\alpha$ varies according to some well-behaved distribution $\rho_\alpha$, the phases $(p(N-1)\rv A+2p\rv X)$ are asymptotically i.i.d., independently of $(\nu_p,\xi_p,y_p)$. Together they identify the limit of the sum \eqref{eq:asymp} as a Cauchy distribution.

Proposition \ref{prop:equidx} is proven in Section \ref{sec:s4}. The proof does not use the chosen ordering of the elements of $\mathcal{P}(\delta,N,\alpha)$ in any essential way. However, this ordering is relevant for describing the behavior of the dominant terms in the sum expression established in Lemma \ref{lemma:sumrv} in the proof of Proposition \ref{prop:main_harm}. 
\subsection*{Proof of Proposition  \ref{prop:main_harm}} 

\begin{proof}[Proof of Proposition \ref{prop:main_harm}]
We know from Proposition \ref{prop:harmonics} that
\[
S(\alpha,h,x,N)=\mathcal{G}(\alpha,x,N,\mathcal{R}(\delta,N,\alpha))+\Delta(\alpha,x,N,\delta),
\]
with
\[
\lim_{\delta\to 0}\lim_{N_\to\infty}\P \left(\frac{|\Delta(\rv A,\rv X,N,\delta)|}{(\log N)^2}>\eps\right)=0,
\]
for any $\eps>0$.

We prove the Proposition by showing that \[\frac1{\hat c}\frac{\mathcal{G}(\rv A,\rv X,N,\mathcal{R}(\delta,N,\rv A))
}{(\log N)^2}\ \xrightarrow[N\to\infty]{\text{dist}}\ \mathfrak{C}_{1/\delta}\ ,\]
when $\rv A\sim\eta$ and $\rv X\big|\rv A=\alpha\sim\rho_\alpha$, and using \[ \mathfrak{C}_{1/\delta}\ \xrightarrow[\delta\to0]{\text{dist}}\ \mathfrak{C},\]
where $\hat c:=(4\pi)^{-1}$, and $\mathfrak{C}_{1/\delta}$ and $\mathfrak{C}$ are as in Lemma \ref{lemma:cauchy}.

Recall the definitions \eqref{eq:def_rv1}--\eqref{eq:def_rv2}. Let $p_k$ be an enumeration of the elements of $\mathcal{P}(\delta,N,\alpha)$ from smaller to larger $|\xi_{p_k}|$. This ordering is what makes Proposition \ref{prop:equidx} apply to the dominant terms of the sum \eqref{eq:asymp}: for every fixed $r$, the proposition gives joint convergence of the first $r$ terms $(p_1,\dots,p_r)$ in this ordering, i.e. precisely the $r$ terms contributing the most to $\mathcal{G}(\rv A,\rv X,N,\mathcal{R}(\delta,N,\rv A))$. Letting $r\to\infty$ ($\delta\to 0$) after $N\to\infty$ then identifies the limit of the full sum with the series $\sum_k \Phi_k/\xi_k$ in Lemma \ref{lemma:cauchy}.

By Propositions \ref{prop:PPP} and \ref{prop:equidx}, the law of 
$\{\xi_{p_k}\}$ converges to a Poisson Point Process $\{\xi_k\}$ in $\left[-\frac{1}{\delta},\frac{1}{\delta}\right]$ with intensity $c=2\frac{6}{\pi^2}$, while
\[\begin{aligned}
\nu_{p_k}\Theta_{p_k}(\rv X,\rv A)+O\left((\log N)^{-1}\right)
\xrightarrow[N\to\infty]{\text{dist}}  Z_k\,\Theta(X_k,Y_k):=\Phi_k,
\end{aligned}
\]
where $X_k,Y_k,Z_k$ are i.i.d. random variables $\sim\text{Uniform}(\T)$, independent of $\xi_{k}$. In particular, by Lemma \ref{lemma:rv}, $\E(\Phi_k)=0$ and $\E(|\Phi_k|)=\frac{1}{2}\frac{1}{24}=\frac{1}{48}$.

Hence, by Lemma \ref{lemma:cauchy},
\[
\sum_{k=1}^{|\mathcal{P}|}\frac{1}{\xi_{p_k}} \Bigl(\nu_{p_k}\Theta_{p_k}(\rv X,\rv A)+O \bigl(\tfrac{1}{\log N}\bigr)\Bigr)\xrightarrow[N\to\infty]{\text{dist}}\sum_k\frac{\Phi_k}{\xi_k}=c\pi\E(|\Phi_k|)\mathfrak{C}_{1/\delta}=\frac{1}{4\pi}\mathfrak{C}_{1/\delta}.
\]
In particular, by Lemma \ref{lemma:sumrv},
\[
4\pi\frac{\mathcal{G}(\rv A,\rv X,N,\mathcal{R}(\delta,N,\rv A))
}{(\log N)^2}\ \xrightarrow[N\to\infty]{\text{dist}}\ \mathfrak{C}_{1/\delta},
\]
which completes the proof.
\end{proof}
\color{black}

\begin{remark}[Multiple singularities]
The result in Proposition \ref{prop:main_harm} can be stated and proven in more generality for a function with multiple (possibly asymmetric) logarithmic singularities following the same steps.

Let $\{\cdot\}:\R\to[0,1)$ denote the decimal part of a real number. For $y_0=0$, $d\geq 0$, and $(C_0^+,C_0^-,\dots,C_d^+,C_d^-)\in\R^{2(d+1)}$, define
\be
\begin{aligned}
\label{eq:sing_pC}     \hat h\left(y_1,\dots,y_d\right)\left(x\right) &= \sum_{i=0}^d C_i^+\,\left[|\log(\{x-y_i\})|-1\right]+ \sum_{i=0}^d C_i^-\,\left[|\log(\{y_i-x\})|-1\right].
\end{aligned}
\ee
Clearly, \eqref{eq:sing_pC} generalizes the logarithmic profile \eqref{eq:sing_p} to multiple singularities.

For $\rv A\sim\eta\in C^\infty(\T)$, and  $(\rv X,\rv Y_1,\dots,\rv Y_d)|\rv A=\alpha\sim\rho_{\alpha}$ uniformly $C^2(\T^{d+1})$, the following result analogous to Proposition \ref{prop:main_harm} holds: 
\[
\frac{S(\rv A,\hat h(\rv Y_1,\dots,\rv Y_d),\rv X,N)}{(\log N)^2}\ \xrightarrow[N\to\infty]{\text{dist}}\ \frac{6c}{\pi} \mathfrak{C},
\]
where $\mathfrak{C}$ is the standard Cauchy distribution and, letting $C_i:=C_i^+-C_i^-$ for $i=0,\dots,d$,
\be\label{eq:cttlimit}
c:=\int_{[0,1]^{d+2}}\left|\sum_{i=0}^d C_i\frac{\norm{x_i+z}(1-\norm{x_i+z})-\norm{x_i-z}(1-\norm{x_i-z})}{2}\right|dx_0\dots dx_ddz.
\ee    

Note that, in particular, $c= 0$ if and only if $\hat h$ is symmetric around each singularity (i.e. $C_i^+=C_i^-$ for every $i$). That is, $(\log N)^2$ becomes the appropriate scaling factor only in the presence of an asymmetric singularity.

We don't write a full proof of this result, but note the key ideas behind it. 

Regarding the Fourier decomposition, the extra singularities contribute with a phase shift factor $\sum C_j e^{-2\pi i k y_j}$, so steps 1--3 hold, with the final expression equivalent to \eqref{eq:asymp} being
\[\begin{aligned}
\frac{\mathcal{G}(\alpha,y,x,N,\mathcal{R}(\delta,N,\alpha))}{(\log N)^2}
&=\sum_{p\in\mathcal{P}(\delta,N,\alpha)}\frac{1}{\xi_p} \Bigl(\nu_p\sum_{i=0}^dC_i\Theta_p(x-y_i,y_p)+O \bigl(\tfrac{1}{\log N}\bigr)\Bigr).
\end{aligned}
\]
Then, given our randomization, it is easy to see that an analogue of Proposition \ref{prop:equidx} holds for the $\R^{2(d+1)}$-valued vectors
\[\bigl(p_1[\alpha(N-1)+2x],\dots,p_r[\alpha(N-1)+2(x-y_d)]\bigr)\bmod 2.\]
The expression \eqref{eq:cttlimit} is therefore a consequence of the results in Section \ref{sec:prelim}.
\end{remark}

\subsection{Proof of Step 1 (Fourier coefficients)} \label{sec:s1}
Let $R_\alpha: \T \to \T$ be a circle rotation defined by $R_\alpha x = x + \alpha \bmod 1$, where $\alpha \in \T \setminus \mathbb{Q}$ is the rotation number. We consider the singular profile $h \in L^2(\T)$ given by
\[
\begin{aligned}
    h(x) &= 2|\log x|+|\log(1-x)|-3.
\end{aligned}
\]
Recall we denote by $S(\alpha,h,x,N):=\sum_{l=0}^{N-1} h(R_\alpha ^l x)$ the $N$-th Birkhoff sum. 
To analyze $S(\alpha,h,x,N)$, we decompose the system into the sum of two single-singularity components.

Let 
\[
h_0(x)=-\log x-1,\qquad x\in (0,1),
\]
define the zero-centered asymmetric log profile on the circle, extended periodically in $\T$. Note $\int_\T h_0(x)\, d x=0$.
\begin{lemma} \label{lemma:ank_vs_bnk}
For $k\in\Z$, let $a_k := \int_0^1 h_0(x) e^{-2\pi i k x} \, d x$ be the Fourier coefficients of $h_0$ and $b_k := \int_0^1 h(x) e^{-2\pi i k x} \, d x$ be the Fourier coefficients of the global function $h$. Then $b_k$ decompose into
\begin{equation}\label{eq:fourier_decomp_sum}
    b_k = [2a_k+a_{-k}].
\end{equation}
\end{lemma}
\begin{remark}
    Since $h_0,h\in L^2(\mathbb T)$, their Fourier coefficients are well-defined. Moreover, $a_0=b_0=0$.
\end{remark}
\begin{proof}
We observe that $h(x)=2h_0(x) +h_0(1-x)$. Applying the Fourier transform to both sides, linearity allows us to transform each component individually. By the spatial shift and reflection properties of the Fourier transform, we have
\[ \int_0^1 h_0(1-x) e^{-2\pi i k x} \, d x = \int_0^1 h_0(-z) e^{-2\pi i k z}\, d z =a_{-k}. \]
Together with the first term, we get \eqref{eq:fourier_decomp_sum}.
\end{proof}

Lemma \ref{lemma:ank_vs_bnk} reduces the global Fourier analysis of $h$ to the study of the single prototype $h_0$. We therefore turn to an asymptotic analysis of the coefficients $a_k$.

\begin{lemma} \label{lemma:asymp}
For $k\geq 1$, the Fourier coefficients $a_k$ of $h_0(x)$ satisfy the asymptotic expansion
\begin{equation}\label{eq:ank_asymptotic}
    a_k = \frac{\log k}{2\pi i k} + O\left(\frac{1}{k}\right), \qquad a_{-k} = -\frac{\log k}{2\pi i k} + O\left(\frac{1}{k}\right).
\end{equation}
Therefore, for $k\neq 0$, it holds
\begin{equation}\label{eq:bnk_asymptotic}
    b_k = \frac{\log |k|}{2\pi i k}+ O\left(\frac{1}{|k|}\right).
\end{equation}
\end{lemma}

\begin{proof}
Clearly, \eqref{eq:bnk_asymptotic} follows from \eqref{eq:ank_asymptotic} and Lemma \ref{lemma:ank_vs_bnk}.

Let us now prove \eqref{eq:ank_asymptotic}. For $k\neq 0$, the constant $-1$ is orthogonal to $e^{2\pi kx}$,
so
\[
a_k = -\int_0^1 \log x\, e^{-2\pi i k x}\,dx.
\]
Integration by parts with $u=\log x$ and $dv=e^{-2\pi i k x}dx$ gives
\[
a_k=\lim_{\eps\to 0}
\left[\frac{\log \eps}{2\pi i k}(e^{-2\pi i k\eps}-1)+\frac{1}{2\pi i k}{ \int_\eps^1} \frac{1-e^{-2\pi i k x}}{x}\,dx\right]=\frac{1}{2\pi i k}\int_0^1 \frac{1-e^{-2\pi i k x}}{x}\,dx.
\]
Making the change of variables $y=2\pi k x$ and splitting the real and imaginary parts yields
\be
a_k=\frac{1}{2\pi i k}\int_0^{2\pi k}\frac{1-e^{-iy}}{y}\,dy=\frac{1}{2\pi i k}\left(\int_0^{2\pi k}\frac{1-\cos y}{y}\,dy+i\int_0^{2\pi k}\frac{\sin y}{y}\,dy\right).\label{eq:aux_fourier}
\ee

Using the classical asymptotics (for $x\to \infty$)
\[
\int_0^x \frac{1-\cos y}{y}\,dy = \log x + O(1),
\qquad
\int_0^x\frac{\sin y}{y}\,dy = O(1),
\]
with $x=2\pi k$, \eqref{eq:aux_fourier} becomes
\[
a_k = \frac{\log k}{2\pi i k} + O \left(\frac{1}{k}\right).
\]

Since $h_0$ is real-valued,
$\displaystyle
a_{-k}=\overline{a_k} = -\frac{\log k}{2\pi i k} + O \left(\frac{1}{k}\right).
$
\end{proof}

Combining Lemmas \ref{lemma:ank_vs_bnk} and \ref{lemma:asymp}, we obtain a sharp description of the Fourier tail of $h$. Note that the leading order of $b_k$ is $\asymp ({2\pi i k})^{-1}{\log|k|}$ because of the asymmetry of the logarithmic singularities of $h$: the contributions of $2a_k$ and $a_{-k}$ do not cancel out. This implies $\sum_{k\leq N}\frac{\log k}{k}\asymp\frac{1}{2}(\log N)^2$, which explains the normalization scale.
The task is now to identify which frequencies contribute at this scale.

\subsection{Proof of Step 2 (Proposition \ref{prop:harmonics})} \label{sec:s2}Following the work of Dolgopyat and Sarig \cite{DolgopyatSarig2019}, we will show that the contribution of large Fourier modes (with respect to $N$) is negligible, and decompose the rest into resonant and non-resonant harmonics. The idea is that although every Fourier mode appears in the expansion of the ergodic sum, the sum behavior is determined only by the frequencies for which $\dist{(k\alpha},\Z)$ is exceptionally small. 

Recall $\mathcal{R}(\delta, N,\alpha)$ denotes the set of resonant harmonics, $\mathcal{P}(\delta, N,\alpha)$ the set of prime resonant harmonics, and $\mathcal{O}(\delta, N,\alpha)$ the set of non-resonant harmonics less than $N$. 

Let $K\subseteq\Z^+$ and recall we define, for $\alpha\in\RQ$, $x\in\T$, and $N\geq 1$
\be\label{eq:defG}
\mathcal{G}(\alpha,x,N,K):=\sum_{k\in K}\frac{\log k}{k}\frac{\sin(\pi N k \alpha)f_{N,k}(x,\alpha)}{\pi\sin(\pi k\alpha)},
\ee
where, for each $k\geq 1$,
\begin{align} \label{eq:deff}
f_{N,k}(x,\alpha)&:=\sin\big(2\pi kx+\pi (N-1)k\alpha\big).
\end{align}

Fix a probability space with the random variables $\rv A\sim \eta$ and $\rv X\big|\rv A=\alpha\sim \rho_\alpha$. Recall we assume $\eta\in C^\infty(\T)$ and $(\rho_\alpha)_{\alpha\in\text{Supp}(\eta)}$ uniformly $C^{2}(\T)$. In particular, $\eta\in L^\infty(\T)$ and the Fourier coefficients of $\rho_\alpha$ satisfy $\hat \rho_{\alpha,k}=O(|k|^{-2})$ uniformly for $k\neq 0$.

The goal of this subsection is to prove Proposition \ref{prop:harmonics}.
First, we introduce two lemmas that provide intuition for the definitions \eqref{eq:defG}--\eqref{eq:deff}. Then, we show on various steps that $S(\alpha,h,x,N)$ and $\mathcal{G}(\alpha,x,N,K)$ are close for a suitable choice of $K$. 

\begin{lemma} For each $K\subset \N$, $\alpha\in\RQ$ and $N\geq 1$, $\mathcal{G}(\alpha,x,N,K)$ as in \eqref{eq:defG} is well-defined and finite $\rho_\alpha-$almost every $x$. \label{lemma:Gwd}  
\end{lemma}
\begin{proof} Let us define 
\[
D_N(\beta) := \sum_{\ell=0}^{N-1} e^{2\pi i \ell \beta}=\frac{\sin(\pi N\beta)}{\sin(\pi \beta)}e^{i\pi(N-1)\beta}.
\]
By the orthonormality of the Fourier basis in $L^2(\mathbb T)$, for any mean-zero function $g\in L^2(\mathbb T)$ with Fourier coefficients $\hat g_k$, for every $\alpha\in\mathbb T$, and Lebesgue--almost every $x\in\T$ we have the exact identity
\[
S(\alpha,g,x,N) = \sum_{k\neq 0} \hat g_k\,D_N(k\alpha)\,e^{2\pi i kx}.
\]

Note that for $\alpha\in\RQ$
\be\label{eq:aux_lemma59}\begin{aligned}
D_N(k\alpha)e^{2\pi i kx}-D_N(-k\alpha)e^{-2\pi i kx}&=2i\frac{\sin(\pi Nk\alpha)}{\sin(\pi k\alpha)}\sin(2\pi k x+\pi(N-1)k\alpha)\\&=2i\frac{\sin(\pi Nk\alpha)}{\sin(\pi k\alpha)}f_{N,k}(x,\alpha),
\\
D_N(k\alpha)e^{2\pi i kx}+D_N(-k\alpha)e^{-2\pi i kx}&=2\frac{\sin(\pi Nk\alpha)}{\sin(\pi k\alpha)}\cos(2\pi k x+\pi(N-1)k\alpha)
.\end{aligned}
\ee
Thus, $\mathcal{G}(\alpha,x,N,K)$ is the projection onto the Fourier modes in $K\cup\{-k:\ k\in K\}$ of the $N-$th Birkhoff sum of a mean-zero function $g$ with Fourier coefficients $\hat g_k=\frac{\log |k|}{2\pi i k}$, for $k\neq 0$, which is in $L^2$. Therefore $\mathcal{G}(\alpha,x,N,K)$ is well-defined Lebesgue--almost surely on $\T$, which implies it holds $\rho_\alpha-$almost surely.
\end{proof}

Analogously to \eqref{eq:deff}, Define
\begin{align} \label{eq:deffh}
\hat f_{N,k}(x,\alpha)&:=\cos\big(2\pi kx+\pi (N-1)k\alpha\big).
\end{align}

\begin{lemma} \label{lemma:fourier} There exists two sequences $(r_k^+),(r_k^-)\subset\R,$ $r_k^\pm=O((\log k)^{-1})$, such that for every $\alpha\in\RQ$ and $x\in\T$ it holds that
\begin{equation}\label{eq:SN-fourier}
\begin{aligned}
S(\alpha,h,x,N)&=\sum_{k\geq 1}\frac{\log k}{k}\frac{\sin(\pi N k \alpha)\left[(1+r_k^-)f_{N,k}(x,\alpha)+r_k^+\hat f_{N,k}(x,\alpha)\right]}{\pi\sin(\pi k\alpha)},
\end{aligned}
\end{equation}
where $f_{N,k}$ and $\hat f_{N,k}$ are as in \eqref{eq:deff}--\eqref{eq:deffh}.
\end{lemma}

\begin{proof} Recall from the proof of Lemma \ref{lemma:Gwd} that for every $\alpha,x\in\mathbb T$ we have the exact identity
$S(\alpha,h,x,N) = \sum_{k\neq 0} b_k\,D_N(k\alpha)\,e^{2\pi i kx},$ where $D_N(\beta) := \sum_{\ell=0}^{N-1} e^{2\pi i \ell \beta}.$

Therefore, it follows from \eqref{eq:bnk_asymptotic} that
\be\label{eq:aux_lemma45}
S(\alpha,h,x,N) = \sum_{k\neq 0}\left(\frac{\log |k|}{2\pi i k}+\hat r_{k}\right)\,D_N(k\alpha)e^{2\pi i kx},
\ee
where $\hat r_{k}=O(|k|^{-1})$. Substituting \eqref{eq:aux_lemma59} into \eqref{eq:aux_lemma45} and using 
\[\begin{aligned}
\hat r_{k}D_N(k\alpha)e^{2\pi i kx}+\hat r_{-k}D_N(-k\alpha)e^{-2\pi i kx}&=\frac{\hat r_{k}+\hat r_{-k}}{2}(D_N(k\alpha)e^{2\pi i kx}+D_N(-k\alpha)e^{-2\pi i kx})\\&+\frac{\hat r_{k}-\hat r_{-k}}{2}(D_N(k\alpha)e^{2\pi i kx}-D_N(-k\alpha)e^{-2\pi i kx}),
\end{aligned}
\] 
gives the result.
\end{proof}

Lemma \ref{lemma:fourier} provides a Fourier representation of the ergodic sum in terms of the coefficients $f_{N,k}$ and $\hat f_{N,k}$.
Since the Fourier coefficients decay like $ k^{-1}\log k$, one expects the contribution of very large frequencies to be negligible after normalization at the scale $(\log N)^2$.

Let $K\subseteq\Z^+$ and define (cf. \eqref{eq:defG} and Lemma \ref{lemma:fourier}) 
\[\begin{aligned}
\Sc(\alpha,x,N,K):&=\sum_{k\in K}\frac{\log k}{k}\frac{\sin(\pi N k \alpha)\left[(1+r_k^-)f_{N,k}(x,\alpha)+r_k^+\hat f_{N,k}(x,\alpha)\right]}{\pi\sin(\pi k\alpha)},
\end{aligned}
\]
for $f_{N,k}$ and $\hat f_{N,k} $ as in \eqref{eq:deff}--\eqref{eq:deffh}, and $r_k^\pm$ as in Lemma \ref{lemma:fourier}.
Note that, from Lemma \ref{lemma:fourier} $\Sc(\alpha,x,N,\Z^+)=S(\alpha,h,x,N)$.

\begin{lemma} \label{lemma:truncation} Fix $N$ and let $K_N:=\{k:1\leq k\leq N\}$. Then, for any $K_N\subseteq K\subseteq \Z^+$ and any $\eps>0$,
\[
\lim_{N\to\infty}\P \left(\frac{|\Sc(\rv A,\rv X,N,K)-\Sc(\rv A,\rv X,N,\Z^+)|}{(\log N)^2}>\eps\right)=0.
\]    
\end{lemma}
\begin{proof}
Let $F_{N,k}(x,\alpha):=\frac{\log k}{k}(1+r_k^-)f_{N,k}(x,\alpha)+r_k^+\hat f_{N,k}(x,\alpha)$. By Lemma \ref{lemma:fourier} and Markov's inequality, we have
    \[\begin{aligned}
    \P \left(\frac{|\Sc(\rv A,\rv X,N,K)-\Sc(\rv A,\rv X,N,\Z^+)|}{(\log N)^2}>\eps\right)
    =\P \left(\frac{|\sum_{k\in \N\setminus K}\frac{\sin(\pi N k \rv A)F_{N,k}(\rv X,\rv A)}{\pi\sin(\pi k\rv A)}|}{(\log N)^2}>\eps\right)
    \\\leq\frac{1}{\eps^2}\E \left(\left[\frac{\sum_{k\in \N\setminus K}\frac{\sin(\pi N k \rv A)F_{N,k}(\rv X,\rv A)}{\pi\sin(\pi k\rv A)}}{(\log N)^2}\right]^2\right).
    \end{aligned}
    \]

Recall $\rv A\sim\eta\in C^\infty(\T)$ and $\rv X\big|\rv A=\alpha\sim \rho_\alpha$, with Fourier coefficients $\hat\rho_{\alpha,k}=O(|k|^{-2})$ uniformly on $\alpha$, for $k\neq 0$ (and $\hat\rho_{\alpha,0}=1$). To compute $\E_{\rv X}(\cdot|\rv A)$, we expand the square of the sum and split it into diagonal and off-diagonal components.

\textbf{Diagonal.} For each $k$, we have a uniform bound $|F_{N,k}(x,\alpha)|=O(\frac{\log k}{k})$, so \[\E_{\rv X}\left[|F_{N,k}(\rv X,\rv A)|^2\big|\rv A=\alpha\right]=O\left(\frac{(\log k)^2}{k^2}\right).\]
We use $\N\setminus K\subset\{k\geq N\}$ and the $L^2$-norm identity for every $k$ (see e.g. \cite[$\S$1]{Katznelson2004})
\[
\left|\E_{\rv A}\left[\frac{\sin^2(\pi Nk\rv A)}{\sin^2(\pi k\rv A)}\right]\right|\leq\norm{\eta}_{\infty}\int_\T\frac{\sin^2(\pi Nk\alpha)}{\sin^2(\pi k\alpha)}d\alpha=\norm{\eta}_{\infty}N,
\]
to conclude that there exists $C>0$ independent of $k$ and $N$ such that, for $N$ large enough,
    \[\begin{aligned}
    \text{diagonal} &\leq \frac{C}{\eps^2(\log N)^4}{N}\sum_{k\geq N}\left(\frac{\log k}{k}\right)^2
   \lesssim \frac{N}{N(\log N)^2}=O((\log N)^{-2})=o(1).
    \end{aligned}
    \]

\textbf{Off-diagonal ($k\neq k'$).}
By the hypothesis $\hat\rho_{\alpha,k}=O(|k|^{-2})$, the definitions of $f_{N,k}$ as \eqref{eq:deff} and $\hat f_{N,k}$ as \eqref{eq:deffh}, $r_k^\pm=O((\log k)^{-1})$, and the definition of $F_{N,k}$, we have
\[\begin{aligned}
\left|\E_{\rv X}[F_{N,k}(\rv X,\rv A)\,F_{N,k'}(\rv X,\rv A)\big|\rv A=\alpha]\right|&\lesssim \frac{\log k}{k}\frac{\log k'}{k'}(\hat\rho_{k-k'}+\hat\rho_{k+k'})\\&= O\left({\frac{\log k}{k}\frac{\log k'}{k'}\frac{1}{|k-k'|^{2}}}\right),\end{aligned}
\]
for every $k,k'\geq 1$, $k\neq k'$, uniformly on $\alpha$.
By Cauchy--Schwarz and the bounds above,
\[
\left|\E_{\rv A}\left[\frac{\sin(\pi Nk\rv A)}{\sin(\pi k\rv A)}\cdot
\frac{\sin(\pi Nk'\rv A)}{\sin(\pi k'\rv A)}\right]\right|
\leq O(N).
\]
Hence, for some $C>0$ independent of $k,k',N$, the off-diagonal contribution is bounded by
\[
\text{off-diagonal}\leq\frac{CN}{\eps^2(\log N)^4}
\sum_{\substack{k\neq k'\\k,k'\geq N}}\frac{\log k}{k}\frac{\log k'}{k'}\frac{1}{|k-k'|^2}.
\]

Let $m=k-k'$, $m\neq 0$, and estimate the $k$-sum for each $m$.
Cauchy--Schwarz gives
\[
\sum_{k\geq N+|m|}\frac{\log k}{k}\frac{\log(k-m)}{k-m}
\leq\sqrt{\sum_{k\geq N}\frac{(\log k)^2}{k^2}}
   \cdot\sqrt{\sum_{k'\geq N}\frac{(\log k')^2}{k^2}}
\lesssim\frac{(\log N)^2}{N},
\]
so summing $\frac{(\log N)^2}{N}\frac{1}{|m|^2}$ over $m\neq 0$ yields 
\[
\text{off-diagonal}=\frac{N}{\eps^2(\log N)^4}
O\left(\frac{(\log N)^2}{N}\right)=O((\log N)^{-2})=o(1).
\]

Combining both contributions concludes the proof.
\end{proof}

Lemma \ref{lemma:truncation} shows that, given our normalization by $(\log N)^2$, the ergodic sum is determined by its truncation to frequencies below $N$. 

The next lemma shows that the error introduced by discarding the coefficients $r_{k}^\pm$ is negligible up the scale $k\leq N(\log N)^2$, thus replacing $\Sc$ by the simpler model sum $\mathcal G$, defined in \eqref{eq:defG}. The advantage of working with $\mathcal G$ is that only the principal asymptotic term of the Fourier coefficients remains.

\begin{lemma} \label{lemma:large_coef} Let $T=N\log ^2 N$ and $K_T:=\{k:1\leq k\leq T\}$. Then, for any $K\subseteq K_T$ and any $\eps>0$,
\[
\lim_{N\to\infty}\P \left(\frac{|\mathcal{G}(\rv A,\rv X,N,K)-\Sc(\rv A,\rv X,N,K)|}{(\log N)^2}>\eps\right)=0.
\]    
\end{lemma}
\begin{proof}
    For $1\leq k \leq T$, define  (c.f. \cite[$\S$3]{DolgopyatFayad2020})
    \[E_1(k):=\left\{\alpha\in \T: k\norm{k\alpha}\leq(\log N)^{-2}\right\},\qquad E_1:=\bigcup_{k=1}^TE_1(k),\]
    and, for $n\geq 1$, define the disjoint Diophantine sorting sets
    \[
    E_2(k,n) = \left\{\alpha\in\T : \frac{n}{({\log N})^{2}}<k\|k\alpha\|\leq\frac{n+1}{(\log N)^{2}} \right\},
    \]
By definition, $E_1(k)^c= \bigsqcup_{n\geq1}E_2(k,n)$ and $E_1^c=\cap_kE_1(k)^c$. 

Since $\rv A\sim\eta\in C^\infty(\T)$, we have for any $0\leq a\leq b$
\be\label{eq:Palpha}\P_{\rv A}(\{\rv A:\norm{k\rv A}\in(a,b)\})\leq 2\norm{\eta}_{\infty}(b-a).\ee 
Thus,
\[
\P_{\rv A} (E_1)\leq\sum_{k=1}^T\norm{\eta}_{\infty}\frac{2}{k(\log N)^{2}}=O\left(\log T/(\log N)^{2}\right)=o(1).
\]

Now, let $\hat F_{N,k}(x,\alpha):=\frac{\log k}{k}\sin(\pi N k \alpha)\left[r_k^-f_{N,k}(x,\alpha)+r_k^+\hat f_{N,k}(x,\alpha)\right]$ and consider the difference of the sums restricted to $E_1^c$:
\[\begin{aligned}
\E \left[|\mathcal{G}(\rv A,\rv X,N,K)-\Sc(\rv A,\rv X,N,K)|^2\chi_{E_1^c}(\rv A)\right]
=\E \left[{\left|\sum_{k\in K}\frac{\hat F_{N,k}(\rv X,\rv A)}{\pi\sin(\pi k\rv A)}\right|^2}\chi_{E_1^c}(\rv A)\right].
\end{aligned}
\]

We proceed as in the proof of Lemma \ref{lemma:truncation}. That is, we expand the square and compute the expectation of the diagonal and off-diagonal summands separately.

\textbf{Diagonal.}
Recall $r_k^\pm=O((\log k)^{-1})$. For each $k$, \[
\E_{\rv X}\left[|\hat F_{N,k}(\rv X,\rv A)|^2\big|\rv A=\alpha\right]\lesssim\frac{(\log k)^2}{k^2}(|r_k^+|+|r_k^-|)^2=O(k^{-2}).\] 
Using the standard inequality $2|\sin(\pi x)|\geq \norm{x}$ and \eqref{eq:Palpha}, we have there exists $C>0$ such that
\[
\begin{aligned}
\text{diagonal}
&\leq{C}\sum_{k=1}^T\frac{1}{k^2}\E_{\rv A}\left[\frac{1}{\norm{k\rv A}^2}\chi_{E_1(k)^c}(\rv A)\right]\leq{C}\sum_{k=1}^T\frac{1}{k^2}\E_{\rv A}\left[\frac{1}{\norm{k\rv A}^2}\sum_{n\geq 0}\chi_{E_2(k,n)}(\rv A)\right]
\\
&\leq{C}\sum_{k=1}^T\frac{1}{k^2}\sum_{n\geq 1}\P_{\rv A}({E_2(k,n)})\frac{k^2(\log N)^4}{n^2} \leq{C}2\norm{\eta}_{\infty}\sum_{k=1}^T\frac{1}{k^2} \sum_{n\geq 1}\frac{k(\log N)^2}{n^2}
\\
&=O((\log N)^2\log T)=O((\log N)^3)=o((\log N)^4).
\end{aligned}
\]

\textbf{Off-diagonal.}
Similarly to the off-diagonal case in the proof of Lemma \ref{lemma:truncation}, we have \[\begin{aligned}\left|\E_{\rv X}\left[\hat F_{N,k}(\rv X,\rv A)\,\hat F_{N,k}(\rv X,\rv A)\big|\rv A=\alpha\right]\right|&\lesssim \frac{\log k}{k}\frac{\log k'}{k'}(|r_k^+|+|r_k^-|)(|r_{k'}^+|+|r_{k'}^-|){|k-k'|^{-2}}\\&=O\left(k^{-1}k'^{-1}|k-k'|^{-2}\right),\end{aligned}\] uniformly on $\alpha$, and for each $k,k'\geq 1$ with $k\neq k'$. By Cauchy-- Schwarz it holds
\[\begin{aligned}
    \E_{\rv A}&\left[\chi_{E_1^c}(\rv A)\,\left|\frac{1}{\sin(\pi k\rv A)}\frac{1}{\sin(\pi k'\rv A)}\right|\right]
\\
&\leq 4\left(\E_{\rv A}\left[\chi_{E_1(k)^c}(\rv A)\,\frac{1}{\norm{k\rv A}^2}\right]\right)^{1/2}\left(\E_{\rv A}\left[\chi_{E_1(k')^c}(\rv A)\,\frac{1}{\norm{k'\rv A}^2}\right]\right)^{1/2}.
\end{aligned}\]
which, by our computations for the diagonal case, is $O\left(\sqrt{kk'}(\log N)^2\right)$.

Hence, the off-diagonal contribution is bounded by
 \[
    \begin{aligned} \text{off-diagonal}&\leq{C'}(\log N)^2\sum_{1\leq k'<k\leq T}\frac{1}{{kk'}}\frac{1}{(k-k')^2}\sqrt{kk'},
    \end{aligned}
\]
for some constant $C'>0$ independent of $k,k',N$.
Let $m=k-k'$ and consider two cases: $m\leq k/2$ and $m> k/2$. In the first case, $k'\geq k/2$, so
 \[
    \begin{aligned} \sum_{1\leq k\leq T}\sum_{m=1}^{k/2}\frac{1}{\sqrt{k(k-m)}}\frac{1}{m^2}\lesssim \sum_{1\leq k\leq T}\frac{1}{k}\sum_{m=1}^{k/2}\frac{1}{m^2}=O(\log T)=O(\log N),
    \end{aligned}
\]
and in the second case, $m^{-2}\lesssim k^{-2}$, and $k'\leq k/2$, so 
\[
    \begin{aligned} \sum_{1\leq k\leq T}\sum_{k'=1}^{k/2}\frac{1}{\sqrt{kk'}}\frac{1}{(k-k')^2}\lesssim \sum_{1\leq k\leq T}\frac{1}{k^{5/2}}\sum_{m=1}^{k/2}\frac{1}{\sqrt{k'}}\lesssim \sum_{1\leq k\leq T}\frac{1}{k^{2}}=O(1),
    \end{aligned}
\]
so the off-diagonal sum (multiplying by $C'(\log N)^2$) is also $O((\log N)^3)=o((\log N)^4)$.

Together, the diagonal and off-diagonal cases imply, by Markov,
\[
\lim_{N\to\infty}\P \left(\frac{|\mathcal{G}(\rv A,\rv X,N,K)-\Sc(\rv A,\rv X,N,K)|}{(\log N)^2}\chi_{E_1^c}(\rv A)>\eps\right)=0,
\]
which, together with $\lim_{N\to\infty}\P_{\rv A} (E_1^c)=1$, concludes the proof.
\end{proof}

Combining Lemmas \ref{lemma:truncation} and \ref{lemma:large_coef}, we see that the study of $S(\alpha,h,x,\Z^+)$ reduces to the study of the simplified sum $\mathcal G$. Large values of $({\sin(\pi k\alpha)})^{-1}$ can occur only when $k\alpha$ is exceptionally close to an integer, and such frequencies are precisely the resonant harmonics. Recall $\mathcal{O}(\delta, N,\alpha)=[1,N]\setminus \mathcal{R}(\delta, N,\alpha)$ is the set of non-resonant harmonics less than $N$. The next lemma shows that the contribution of $\sum_{k\in \mathcal{O}(\delta, N,\alpha)}$ is asymptotically negligible.

\begin{lemma} \label{lemma:nonres} For any $\epsilon>0$ $\exists\delta_0>0$ such that, for any $\delta_0>\delta>0$ and $N$ large enough,
\[
\P \left(\frac{|\mathcal{G}(\rv A,\rv X,N,\mathcal{O}(\delta,N,\rv A))|}{(\log N)^2}>\eps\right)<\frac{\eps}{2}.
\]    
\end{lemma}
\begin{proof}
    We will proceed analogously to the proof of Lemma \ref{lemma:large_coef}. For $n\geq 1$, define the Diophantine sorting sets
\begin{equation}\label{eq:B_sets}
    B(k,n,\delta) = \left\{\alpha\in\T : \frac{n}{{\delta\log N}}<k\|k\alpha\|\leq\frac{n+1}{\delta\log N} \right\}, \quad B(k,\delta) = \bigsqcup_{n\geq 1}B(k,n,\delta).
\end{equation}
Let $g_{N,k}(x,\alpha):=\sin(\pi N k \alpha)f_{N,k}(x,\alpha)$. By Markov and \eqref{eq:defG},
    \[
    \begin{aligned} \P \left(\frac{|\mathcal{G}(\rv A,\rv X,N,\mathcal{O}(\delta,N,\rv A))|}{(\log N)^2}>\eps\right)\leq\frac{1}{\eps^2(\log N)^4}\E \left({|\mathcal{G}(\rv A,\rv X,N,\mathcal{O}(\delta,N,\rv A))|}^2\right)
    \\
    \leq\frac{1}{\eps^2(\log N)^4}\E \left(\left[{\sum_{k\leq N}\chi_{\{k\in \mathcal{O}(\delta,N,\rv A)\}}(\rv A)\frac{\log k}{k}}\frac{g_{N,k}(\rv X,\rv A)}{\pi \sin(\pi k \rv A)}\right]^2\right)
    \end{aligned}
    \]
Again, compute the diagonal and off-diagonal terms.

    \textbf{Diagonal.} $\E_{\rv X}[|g_{N,k}(\rv X,\rv A)|^2\big|\rv A=\alpha]=O(1)$, so we have for some $C>0$
    \[
    \begin{aligned} 
    \text{diagonal}&\leq\frac{C
    }{\eps^2(\log N)^4}\sum_{k\leq N}\frac{(\log k)^2}{k^2}\E_{\rv A}\left(\frac{\chi_{\{k\in \mathcal{O}(\delta,N,\rv A)\}}(\rv A)}{\sin^2(\pi k \rv A)}\right).
    \end{aligned}
    \]
    If the harmonic $k$ is not prime and it is associated to the prime harmonic $p$, $k=rp$, then $\snorm{k\alpha}=r\snorm{p\alpha}$ (see \cite[Lemma 4.1]{DolgopyatSarig2019}). Therefore, for each $\alpha$, $\mathcal{O}(\delta, N,\alpha)\subseteq[1,N]\cap\{k:\alpha\in B(k,\delta)\}$ and, in particular, $\chi_{\{k\in \mathcal{O}(\delta,N,\alpha)\}}(\alpha)\leq \chi_{B(k,\delta)}(\alpha)=\sum_{n\geq 1}\chi_{B(k,n,\delta)}(\alpha)$. Together with the definition \eqref{eq:B_sets} of $\{\chi_{B(k,n,\delta)}\}_n$ and its orthogonality for fixed $k$, this implies
    \[
    \begin{aligned} \text{diagonal}&\leq\frac{C'}{\eps^2(\log N)^4}\sum_{k\leq N}\frac{(\log k)^2}{k^2}\E_{\rv A}\left({\sum_{n\geq 1}\chi_{B(k,n,\delta)}(\rv A)\frac{1}{\norm{k\rv A}^2}}\right)
    \\ &\leq\frac{C'\norm{\eta}_\infty}{\eps^2(\log N)^4}\sum_{k\leq N}\frac{(\log k)^2}{k^2}\sum_{n\geq 1}2\delta\frac{k(\log N)}{n^2}\\
    &\leq\frac{2\delta C'\norm{\eta}_\infty}{\eps^2(\log N)^3}\sum_{k\leq N}\frac{(\log k)^2}{k}\sum_{n\geq 1}\frac{1}{n^2}=\frac{C''\delta}{\eps^2}\leq\frac{\eps}{4},
    \end{aligned}
\]
for $\delta$ small enough, where $C',C''>0$ are constants independent of $k$ and $N$, and where we have used \eqref{eq:Palpha}. 

\textbf{Off-diagonal.} Again, we have $|\E_{\rv X}\left[g_{N,k}(\rv X,\rv A)\,g_{N,k'}(\rv X,\rv A)|\rv A=\alpha\right]|\lesssim {|k-k'|^{-2}}$ for each $k,k'\geq 1$ with $k\neq k'$ (cf. Lemmas \ref{lemma:truncation} and \ref{lemma:large_coef}) and, by Cauchy--Schwarz and the computations above, it holds for some $\hat C>0$
\[\begin{aligned}
    \E_{\rv A}&\left[\chi_{\{k\in \mathcal{O}(\delta,N,\rv A)\}\cap \{k'\in \mathcal{O}(\delta,N,\rv A)\}}(\rv A)\,\left|\frac{1}{\sin(\pi k\rv A)}\frac{1}{\sin(\pi k'\rv A)}\right|\right]
\\
&\leq \left(\E_{\rv A}\left[\frac{\chi_{\{k\in \mathcal{O}(\delta,N,\rv A)\}}(\rv A)}{\sin^2(\pi k\rv A)}\right]\right)^{1/2}\left(\E_{\rv A}\left[\frac{\chi_{\{k'\in \mathcal{O}(\delta,N,\rv A)\}}(\rv A)}{\sin^2(\pi k'\rv A)}\right]\right)^{1/2}
\leq \hat C \delta\log N \sqrt{kk'}.
\end{aligned}
\]
Hence, there exists $\hat C'>0$ independent of $k,k'$, and $N$ such that the off-diagonal contribution is bounded by
 \[
    \begin{aligned} \text{off-diagonal}&\leq\frac{\delta \hat C'}{\eps^2(\log N)^3}\sum_{1\leq k'<k\leq N}\frac{\log k\log k'}{\sqrt{kk'}}\frac{1}{(k-k')^2}.
    \end{aligned}
\]
Let $m=k-k'$ and consider two cases: $m\leq k/2$ and $m> k/2$ (i.e. $k'\leq k/2$). We have
\[
\begin{aligned} \sum_{1\leq k\leq N}\sum_{m=1}^{k/2}\frac{\log k\log (k-m)}{\sqrt{k(k-m)}}\frac{1}{m^2}\lesssim \sum_{1\leq k\leq N}\frac{(\log k)^2}{k}\sum_{m=1}^{k/2}\frac{1}{m^2}=O((\log N)^3),
\end{aligned}
\]
and
\[
    \begin{aligned} \sum_{1\leq k\leq N}\sum_{k'=1}^{k/2}\frac{\log k\log k'}{\sqrt{kk'}}\frac{1}{(k-k')^2}\lesssim \sum_{1\leq k\leq N}\frac{\log k}{k^{5/2}}\sum_{m=1}^{k/2}\frac{\log k'}{\sqrt{k'}}\lesssim \sum_{1\leq k\leq N}\frac{(\log k)^2}{k^{2}}=O(1),
    \end{aligned}
\]
so the off-diagonal sum also satisfies, for some $\hat C''>0$ and provided $\delta$ small enough and $N$ large enough, 
\[
\text{off-diagonal}< \frac{\hat C''\delta}{\eps^2}<\frac{\eps}{4}.
\]

This concludes the proof.
\end{proof}

Lemma \ref{lemma:nonres} establishes that every contribution except the resonant ones are negligible after normalization by $(\log N)^2$. The proof of Proposition \ref{prop:harmonics} follows by combining the successive approximation steps obtained above.
\begin{proof}[Proof of Proposition \ref{prop:harmonics}]
Fix $\varepsilon>0$. By the Lemma \ref{lemma:truncation}, for $N$ sufficiently large,
\[
\P  \left(\frac{|S(\rv A,h,\rv X,\Z^+)-\Sc(\rv A,\rv X,N,K)|}     {(\log N)^2}>\varepsilon\right)=o(1),
\]
for any $K\subseteq \Z^+$ with $K_N=\{1,\dots,N\}\subseteq K$.

Next, letting $T=N(\log N)^2$ Lemma \ref{lemma:large_coef} shows that, for any $K\subseteq K_T=\{1,\dots,T\}$,
\[
\P  \left(\frac{|\Sc(\rv A,\rv X,N,K)-\mathcal G(\rv A,\rv X,N,K)|}     {(\log N)^2}>\varepsilon\right)=o(1).
\]
Therefore, for any $K_N\subseteq K\subseteq K_T$,
\be\label{eq:auxproph}
\frac{|S(\rv A,h,\rv X,\Z^+)-\mathcal G(\rv A,\rv X,N,K)|}{(\log N)^2}\xrightarrow{\P }0.
\ee

Define $K:=\mathcal R(\delta,N,\alpha) \sqcup \mathcal O(\delta,N,\alpha)$. Then
\[
\mathcal G(\alpha,x,N,K)
=
\mathcal G(\alpha,x,N,\mathcal R(\delta,N,\alpha))
+
\mathcal G(\alpha,x,N,\mathcal O(\delta,N,\alpha)),
\]
and, by definition, $K_N\subset K\subseteq K_T,$ so \eqref{eq:auxproph} applies.

By Lemma \ref{lemma:nonres}, for every $\eps>0$ there exists $\delta_0>0$ such that whenever $0<\delta<\delta_0$ and $N$ is large,
\[
\P  \left(\frac{|\mathcal G(\rv A,\rv X,N,\mathcal O(\delta,N,\rv A))|} {(\log N)^2} >\eps \right) <\frac{\eps}{2}.
\]

Combining the above estimates and applying the triangle inequality yields
\[
\P  \left( \frac{ |S(\rv A,h,\rv X,N) -\mathcal G(\rv A,\rv X,N,\mathcal R(\delta,N,\rv A))| } {(\log N)^2} >\eps \right)<\eps,
\]
for all sufficiently large $N$.

Since $\eps>0$ is arbitrary, the conclusion follows.
\end{proof}

\subsection{Proof of Step 3 (Lemma \ref{lemma:sumrv})} \label{sec:s3}We show the asymptotic behavior of the term of the sum associated with each prime resonant harmonic.

Recall we write $\snorm\cdot$ to denote the signed distance to the nearest integer, that is, $x=\snorm{x}+m$, where $\snorm{x}\in(-\frac{1}{2},\frac{1}{2}]$ and $m\in\Z$, and write $\norm{x}:=\dist(x,\Z)=|\snorm{x}|$. For $x$ and $m$ as above, we write $\ssnorm{x}=(-1)^m\snorm{x}$.

\begin{proof}[Proof of Lemma \ref{lemma:sumrv}]

Throughout this proof, we follow the strategy of \cite[Proposition 3.5]{DolgopyatSarig2019}, with a minor adaptation to account for the extra weight $\log k$.
Fix $\delta>0$ small, $\alpha\in\T\setminus\Q$, and $N\geq 1$, and set $ \mathcal{R}=\mathcal{R}(\delta, N,\alpha)$, $\mathcal{P}=\mathcal{P}(\delta, N,\alpha)$, and $R_p := \{$resonant harmonics in $\mathcal{R}$, associated to $p\}$.  The following three claims are proved in Section 4.3 of \cite{DolgopyatSarig2019}\footnote{In fact, Claim 1 (b) in \cite{DolgopyatSarig2019} has a minor inaccuracy: the term $(-1)^{(k-1)m}$ is missing in their formula. However, the proof in \cite{DolgopyatSarig2019} gives the formulation we use.}:

\textsc{Claim 1:} For $p \in \mathcal{P}$, set
$ L_p := \min \left( \left\lfloor \frac{1}{2}\,\delta p \log N \right\rfloor, \left\lfloor \frac{N (\log N)^2}{p} \right\rfloor \right). $
For every $1 \leq k \leq L_p$,

\begin{enumerate}
\item[(a)] $kp \in \mathcal{R}_p(\delta)$ and $((kp\alpha)) = k((p\alpha)),$
\item[(b)] Let $m=p\alpha-((p\alpha))$, then $\frac{1}{2\pi kp\,\sin(\pi kp\alpha)}=\frac{(-1)^{(k-1)m}}{p\langle\langle p\alpha\rangle\rangle}\left(\frac{1}{2\pi^2 k^2}+O \bigl(\|p\alpha\|^2\bigr)\right).$
\end{enumerate}

\textsc{Claim 2:} If  $p \in \mathcal{P}$ and  $kp \in \mathcal{R}_p(\delta)$, then $\frac{1}{2\pi kp \sin(\pi kp\alpha)} = O \left( \frac{1}{k^2 p \|p\alpha\|} \right)$.

\textsc{Claim 3:} If $p \in \mathcal{P}$, then
$\{kp : k = 1, \dots, L_p\} \subset \mathcal{R}_p(\delta) \subset \left\{kp : k = 1, \dots, \frac{N (\log N)^2}{p}\right\}.
$

Every resonant harmonic $j\in\mathcal{R}$ is associated to a unique prime resonant harmonic $p\in\mathcal{P}$ via $j=kp$ for some integer $k\geq 1$. More precisely, by \eqref{eq:defG}, we may write
\[
\mathcal{G}(\alpha,x,N,\mathcal{R}) =\sum_{p\in\mathcal{P}}\;\sum_{kp\,\in\,\mathcal{R}_p(\delta)}\log(kp)\frac{\sin(Nkp\alpha)f_{N,kp}(x,\alpha)}{\pi kp\sin(\pi kp\alpha)},
\]
where $f_{N,k}(x,\alpha):= \sin(\pi(N-1)k\alpha + 2\pi kx)$.

We have the point-wise bound $|f_{N,kp}(x,\alpha)|\leq 1 $, so Claim 2 implies that, for every $p\in\mathcal{P}$,
\[\begin{aligned}
\sum_{kp\in\mathcal{R}_p(\delta)}
\log k
\frac{\sin(Nkp\alpha)f_{N,kp}(x,\alpha)}
{\pi kp\sin(\pi kp\alpha)}=O\left(\frac{1}{p\norm{p\alpha}}\sum_{k\geq 1}\frac{\log k}{k^2}\right)=\left(\frac{1}{p\norm{p\alpha}}\right),
\end{aligned}
\]
by the summability of $k^{-2}\log k$, and Claims 2 and 3 imply that, for every $p\in\mathcal{P}$,
\[\begin{aligned}
\sum_{k=L_p+1}^{{N(\log N)^2}/{p}}
\frac{\sin(Nkp\alpha)f_{N,kp}(x,\alpha)}
{\pi kp\sin(\pi kp\alpha)}=O\left(\frac{1}{p\norm{p\alpha}}\sum_{k=L_p+1}^{{N(\log N)^2}/{p}}\frac{1}{k^2}\right).
\end{aligned}
\]
If $L_p = \floor{\tfrac{1}{2}\delta p\log N}$, the error term equals 
\[
O\left(\sum_{k=L_p+1}^{{N(\log N)^2}/{p}}\frac{1}{k^2}\right)=O\left(\frac{1}{\delta p\log N}\right),
\]
where we have used that $\sum_{k=1}^K \tfrac{1}{k^2}=\tfrac{\pi^2}{6}-O(\tfrac{1}{K})$. Otherwise $L_p= \floor{ N (\log N)^2/p}$ and the error term is zero.

Applying Claim 1 and the previous computations, we obtain
\[
\begin{aligned}
\mathcal{G}&(\alpha,x,N,\mathcal{R})\\& = \sum_{p \in \mathcal{P}} \frac{1}{p \langle\langle p\alpha \rangle\rangle} \left( \log p\sum_{k=1}^{L_p}(-1)^{(k-1)m} \left( \frac{\sin(Nkp\alpha)f_{N,kp}(x,\alpha)}{\pi^2 k^2} + O(\|p\alpha\|^2) \right) + O\left(1\right) \right) 
\\&
= \sum_{p \in \mathcal{P}} \frac{1}{p \langle\langle p\alpha \rangle\rangle} \left( \log p\sum_{k=1}^{L_p} (-1)^{(k-1)m}\frac{\sin(Nkp\alpha)f_{N,kp}(x,\alpha)}{\pi^2 k^2} + O(\log p\,\|p\alpha\|^2L_p) + O\left(1\right) \right) 
\end{aligned}
\]
where the constants in the $O(\cdot)$ terms depend on $\delta$ but are uniform in $p$ and $N$.

For $p \in \mathcal{P}$, $\|p\alpha\| < \frac{1}{\delta p \log N}$, so \[\log p\,\|p\alpha\|^2 L_p \leq \frac{\frac{1}{2} \delta p \log N}{\delta^2 p^2 (\log N)^2} = O \left( \frac{\log p}{\delta p\log N} \right)=O(1).\]

Also, for $p \in \mathcal{P}$, $p \leq N$, so $L_p \geq \min(\lfloor \frac{1}{2}\delta p \log N \rfloor, \lfloor (\log N)^2 \rfloor)$, hence
\[
\sum_{k=L_p+1}^{\infty} (-1)^{(k-1)m}\frac{\sin(Nkp\alpha)f_{N,kp}(x,\alpha)}{\pi^2 k^2} = O \left( \frac{1}{\delta p \log N} +\frac{1}{(\log N)^2}\right).
\]

Putting all together yields
\[
\begin{aligned}
\mathcal{G}(\alpha,x,N,\mathcal{R}) &= \sum_{p \in \mathcal{P}} \frac{1}{p \langle\langle p\alpha \rangle\rangle} \left( \log p\sum_{k\geq 1} (-1)^{(k-1)m}\frac{\sin(Nkp\alpha)f_{N,kp}(x,\alpha)}{\pi^2 k^2}+O(1)\right)
\end{aligned}
\]
Dividing by $(\log N)^2$ gives \eqref{eq:asymp}.
\end{proof}

\subsection{Proof of Step 4 (Proposition \ref{prop:equidx})} \label{sec:s4}Recall we defined the random variables
\[\{(\nu_p:=\frac{\log p}{\log N},\quad \xi_p:=p\ssnorm{p\rv A}{\log N},\quad y_p:= Np\rv A\bmod 2)\}_{p\in\mathcal{P}(\delta,N,\rv A)},\] 
which, by Propostion \ref{prop:PPP}, converge as $N\to\infty$ to a Poisson point process on $[0,1] \times \left[-\frac{1}{\delta}, \frac{1}{\delta}\right]\times[0,2], $ with constant intensity $c={6}/{\pi^2}$. We show convergence as $N\to\infty$ of the variables 
\[\{(N-1)p\rv A+2p\rv X\bmod 2\}_{p\in\mathcal{P}(\delta,N,\rv A)}.\]

\begin{proof}[Proof of Proposition \ref{prop:equidx}]
Similarly to Lemma \ref{lemma:sumrv}, we follow the strategy of \cite[Proposition 3.7]{DolgopyatSarig2019}.

First, we claim that for any $r\geq 1$ and $\eps>0$ there exists $\delta_0>0$ so that for any $\delta_0>\delta>0$ there exists $N_0\geq 1$ satisfying $\P_{\rv A}(\{\rv A:\ |\mathcal{P}(\delta,N,\rv A)|\geq r\})>1-\eps/2,$ $ \forall N\geq N_0.$

Indeed, by Lemma \ref{prop:PPP} $|\mathcal{P}(\delta,N,\rv A)|$ converges in distribution as $N\to\infty$ to a Poisson process with expectation $\frac{24}{\delta \pi^2}$, so taking $\delta\ll\frac{\pi^2}{24r}$ yields the existence of such $N_0$. 

Let $\mathcal{V}$ be an arbitrary neighborhood of the uniform distribution on $[0,2]^{r}$. 
Without loss of generality, we can assume that $\mathcal{V}$ is given by a basic open set of the form
\[
\mathcal{V} = \bigcap_{m=1}^{m_0} \left\{ \mu : \left| \int_{[0,2]^{r}} f_m \, d\mu \right| \leq \hat\varepsilon \right\},
\]
where $\hat \varepsilon > 0$, $m_0 \in \mathbb{N}$, and the test functions are non-trivial characters
\[
f_m(x_1,\ldots,x_{r}) = e^{i \pi \sum_{k=1}^{r} l_k^{(m)} x_k}
\]
associated with non-zero frequency vectors $l^{(m)} := \bigl(l_1^{(m)},\ldots,l_{r}^{(m)}\bigr) \in \mathbb{Z}^{r} \setminus \{0\}$ for $1 \leq m\leq m_0$.

Given $\eps$, $r$, $\delta$, and $N_0$ as above, fix $N\geq N_0$ and $ \alpha\in \hat A(\delta,N,r):=\{\alpha:\ |\mathcal{P}(\delta,N,\alpha)|\geq r\}.$
Recall $\P_{\rv A}(\hat A(\delta,N,r))>1-\eps/2$. Order the elements of $\mathcal{P}(\delta,N,\alpha)$ as in the statement of the Proposition, i.e. increasing with respect to $\xi_p$, and let $\mu_N(\alpha)$ denote the distribution of the random vector
\[
Z(\mathcal X_\alpha) = \bigl(p_1[\alpha(N-1)+2\mathcal X_\alpha], \, \dots, \, p_r[\alpha(N-1)+2\mathcal X_\alpha] \bigr) \bmod 2
\]
where $\mathcal X_\alpha\sim \rho_\alpha\in C^{2}(\T)$, with decaying Fourier coefficients $\hat\rho_{\alpha,k}\to 0$ as $|k|\to\infty$ (uniformly on $\alpha$). To establish that $\mu_N \in \mathcal{V}$ for all sufficiently large $N$, it suffices to show that for each fixed character $f_m$ ($1 \leq m \leq m_0$) and $\hat \eps>0$, the expectation satisfies
\[
\lim_{N \to \infty} \left |\E_{\mathcal X_\alpha} [f_m(Z)]\right| <\hat\eps
\]
uniformly for $\alpha \in \hat A(N,\delta,r)$ up to a set of small measure. We can write this expectation explicitly as
\[\begin{aligned}
\E_{\mathcal X_\alpha} [f_m(Z(\mathcal X_\alpha))] &= \E_{\mathcal X_\alpha} \left[ \exp\left( i \pi \sum_{k=1}^r  l_k^{(m)} p_k[\alpha(N-1) + 2\mathcal X_\alpha] \right) \right]
\\&=\mathcal{Z}^{(m)}
(\alpha,N)\,\E_{\mathcal X_\alpha} \left[ \exp\left( 2 \pi i \mathcal X_\alpha \sum_{k=1}^r l_{k}^{(m)} p_k \right) \right]=\mathcal{Z}^{(m)}(\alpha,N)\,\hat\rho_{\alpha,\ell},
\end{aligned}
\]
for $\ell={-\sum_{k=1}^r l_{k}^{(m)} p_k}$, and where $|\mathcal Z^{(m)}(\alpha, N)| = 1$. 

Let $L(\hat\eps)\geq 1$ be so that $|\hat\rho_{\alpha,k}|<\hat\eps$ for any $\alpha\in\T$ and $|k|\geq L$, which exists by the uniform decaying of $\hat\rho_{\alpha,k}$. Then, it is enough to show that there exists $\{A(\delta,N,r)\}_{N\geq N_0}$, $A(\delta,N,r)\subset\hat A(\delta,N,r)$ with $P_{\rv A}(A(\delta,N,r))>1-\eps$, satisfying that for any given $l\in\Z^r\setminus\{0\}$ and $\hat\eps>0$ there exists $N_{l,\hat\eps}\geq N_0$ so that $|\sum_{k=1}^r l_{k} p_k|\geq L(\hat\eps)$ for every $\alpha\in A(\delta,N,r)$, provided $N\geq N_{l,\hat\eps}$.

Now, by the convergence of $\{\nu_p\}$ given by Proposition \ref{prop:PPP}, we have, for any $R>0$,
\[
\P_{\rv A}\left(\left\{\alpha:\min_{1\leq k< k'\leq r}\max\left\{\tfrac{p_k}{p_{k'}},\tfrac{p_{k'}}{p_{k}}\right\}>R\right\}\ \vert \ \hat A(\delta,N,r)\right)\xrightarrow[N\to\infty]{}1, \]
and, for any $\alpha\in\RQ$, by definition of $\mathcal{P}(\delta,N,\alpha)$ it holds that $\lim _{N\to\infty}\min\{p:p\in\mathcal{P}(\delta,N,\alpha)\}=\infty$.
Thus, for each $n\geq 1$ $\exists N_n$ so that, setting \[
E_{N,n}:=\left\{\alpha\in \hat A(\delta,N,r):\min_{1\leq k< k'\leq r}\max\left\{\tfrac{p_k}{p_{k'}},\tfrac{p_{k'}}{p_{k}}\right\}>n,\ \min_{1\leq k\leq r} p_k>n\right\},\]
$\P_{\rv A}(E_{N,n})>\P_{\rv A}(\hat A(\delta,N,r))-\eps 2^{-n-1}$ for every $N\geq N_n$. For each $N$, let
\[
A(\delta,N,r):=\bigcap_{n:\ N_n\leq N}E_{N,n}.
\]
Then $\P_{\rv A}(A(\delta,N,r))>1-\eps$ for every $N$.

Now, given a non-zero vector $l=(l_1,\dots,l_r)$ and $\hat\eps>0$, set $k_*=\text{argmax}\{p_k:l_k\neq 0\}$ and let $R:=2\max\{\sum_{k=1}^r|l_k|,\ L(\hat\eps)\}$. Then, for every $N\geq N_{R_l}$, for every $\alpha \in A(\delta,N,r)\subseteq E_{N,R}$, it holds that $\sum_{k\neq k_*}p_k|l_k|<\frac{p_{k_*}}{R}\sum_{k\neq k_*}|l_k|<\frac{p_{k_*}}{2}$, so $|\sum_{k=1}^r l_kp_k|> R/2\geq L$. This concludes the proof.

\end{proof}

\appendix

\section{Sketch of the proof of Theorem \ref{thm:mainbisbis}}\label{app:1}
It is only left to show the case $a=0$, since the rest follow from Theorem \ref{thm:mainbis}. As in the proof of Theorem \ref{thm:mainbis}, the theorem follows from the two following propositions.

\begin{proposition}\label{prop:main_distbis}
    Let ($\cF,\hat\eta$) be \textit{good} in the sense of Definition \ref{def:good}. Suppose $f(p_0^{W})-\bar f^W=0 $, $\hat\eta$-almost surely. Moreover, suppose that for any family of probability distributions $(\hat\rho_W)_{W\in\cF}$ uniformly $ C^2(\T)$, it holds that, under $\rv W\sim\hat\eta$ in $\cF$ and $\rv X\big|\rv W=W\sim\rho_{W}$ we have
\[
\frac{I(\rv W,\rv X,t)}{(\log t)^2}\ \xrightarrow[t\to\infty]{\text{dist}}\ 0.
\]
Then Theorem \ref{thm:mainbisbis} holds for $a=0$.
\end{proposition}

\begin{proposition}
    \label{prop:main_redbis} Let ($\cF,\hat\eta$) be \textit{good} in the sense of Definition \ref{def:good}. Suppose $f(p_0^{W})-\bar f^W=0 $, $\hat\eta$-almost surely.
    Let $\rv W\sim\hat\eta$ and $\rv X\big|\rv W=W\sim \hat\rho_{W}$ uniformly $C^2(\T)$. Then 
    \[
    \frac{I(\rv W,\rv X,t)}{(\log t)^2}\ \xrightarrow[t\to\infty]{\text{dist}}\ 0.
    \]
\end{proposition}

Theorem \ref{thm:mainbisbis} clearly follows from combining Propositions \ref{prop:main_distbis} and \ref{prop:main_redbis}.
We give a sketch of the proofs of the propositions.

\subsection{Sketch of the proof of Proposition \ref{prop:main_distbis}.} We proceed exactly as in Section \ref{sec:proof3}. Analogously to Lemma \ref{lemma:distx}, we have the following lemma.

\begin{lemma}\label{lemma:distxbis}
    Suppose that, under $\rv W\sim\hat\eta$ and $\rv X\big|\rv W=W\sim\ph^W/\mu^W$
     \[
     \frac{I(\rv W,\rv X,t)}{(\log t)^2}\ \xrightarrow[t\to\infty]{\text{dist}}\ 0,
     \]
     then taking $\rv W\sim\hat \eta$ in $\cF$ and $\rv P|\rv W=W\sim\text{Uniform}(\M^W)$ we have
    \[
    {\frac{\int_0^{t} f \left(\phi^{\rv W}_s\rv P\right)ds - t\bar f^W}{(\log t)^2}\ \xrightarrow[t\to\infty]{\text{dist}}\ 0.}
    \]
\end{lemma}

The result then follows from combining Lemma \ref{lemma:distxbis} and Lemma \ref{lemma:relax}.

\subsection{Sketch of the proof of Proposition \ref{prop:main_redbis}.}
From Lemma \ref{lemma:fhat}, we have that $\hat f^W=g^W$ is zero mean and of bounded variation $\hat\eta$-almost surely. 

From Lemma \ref{lemma:maxN} and an analogous to Lemma \ref{lemma:bound} replacing $c_0 Y$ by 0 and removing the $f(p_0^{W})-\bar f^W$ term, we get the following analogous to Lemma \ref{cor:bounds}.
\begin{lemma}
    \label{cor:boundsbis} To prove Proposition \ref{prop:main_redbis} we may assume without lost of generality
    \be\gamma:=\inf_\cF\inf_\T\ph^W>0,\quad \mu:=\sup_\cF \mu^W<\infty,\label{eq:def_gammabis}\ee
    \be V:=\sup_\cF\left(\text{Var}(g^W)+\text{Var}(\psi^W)+\norm{g^W}_\infty+\norm{\psi^W}_\infty\right)<\infty,\label{eq:def_Vbis}\ee
    and, uniformly on $W$,
    \be\label{eq:boundsNbis}
    N_*(W,t)\asymp t.
    \ee
\end{lemma} 
Hence, by Lemma \ref{lemma:red3}, if follows that
\[
\frac{S\left(\alpha(\rv W),\hat f^{\rv W},\rv X,k\right)}{(\log k)^2}\xrightarrow{\P } 0 \quad \text{as } k\to\infty,
\]
as well as the following analogous to Lemma \ref{lemma:red1}.

\begin{lemma}\label{lemma:red1bis}
Assuming \eqref{eq:def_gammabis}--\eqref{eq:boundsNbis}, it holds that
\[
    \frac{|I(\rv W,\rv X,t)-S(\alpha(\rv W),\hat{f}^{\rv W},\rv X,N(\rv W,\rv X,t))|}{(\log t)^2}\xrightarrow{\P } 0 \quad \text{as } t\to\infty.
\]
\end{lemma}

The result then follows from noting $\frac{\log N(W,x,t)}{\log t}\to 1$ point-wise as $t\to\infty$.

\section{Proof of Lemma \ref{claim:boundsEsums}}\label{sec:boundsEsums}
We divide the proof into two sections, corresponding to the proofs of \eqref{eq:small} and \eqref{eq:large}.

\subsection{Small Birkhoff sums}
\begin{proof}[Proof of \eqref{eq:small} in Lemma \ref{claim:boundsEsums}]
For $\alpha\in (0,1)\setminus \Q$, let $\alpha=[0;a_1,a_2,\dots]$ be its continued fraction expansion and let $p_n/q_n$ be the $n$th convergent. The denominators $q_n$ provide a sequence of times at which such sums can be controlled by a truncation argument and the Denjoy--Koksma inequality \cite{Herman1979}. Indeed, truncating $h$ in a $q_n^{-2}$-neighbourhood of its singularities gives a function 
\[h_n=h\chi_{\norm{x}>q_n^{-2}},\] of bounded variation, satisfying
\[
\operatorname{Var}(h_n)\lesssim\log q_n,\quad \int_\T h_n=-3q_n^{-2}\log q_n^2-3(1-q_n^{-2})\log(1-q_n^{-2})\lesssim q_n^{-2}\log q_n^2.
\]

for $n$ large enough. Thus, by the Denjoy--Koksma inequality,
\[
\left|S(\alpha,h_n,x,q_n)\right|\leq\left|q_n\int_{\T}h_n\right|+\left|S(\alpha,h_n,x,q_n)-q_n\int_{\T}h_n\right|\lesssim\log q_n,
\]
uniformly in $x$. On the other hand, the set of $x$ for which one of the first $q_n$ iterates enters the truncation neighborhood $\{\norm{x}\leq q_n^{-2}\}$ has measure $O(q_n^{-1})$. Since the denominators $q_n$ of any irrational rotation grow at least exponentially, these exceptional sets are summable. Hence, by the Borel--Cantelli lemma, for almost every $x$ 
\[
|S(\alpha,h,x,q_n)|\lesssim\log q_n,
\]
for all sufficiently large $n$.
\end{proof}
%Note that, using Ostrowski expansion, the typical growth $\sum a_n=O(n^{1+\epsilon})$ for any $\eps>0$ (see e.g. \cite{DiamondVaaler1986}), and the results above, we obtain that for almost every $(\alpha,x)\in\T^2$\[\limsup_{N\to\infty}\frac{|S(\alpha,h,x,N)|}{(\log N)^{2+\eps}}=0,\quad\forall \eps>0.\]

\subsection{Large Birkhoff sums}

  For $\alpha\in (0,1)\setminus \Q$, let $\alpha=[0;a_1,a_2,\dots]$ be its continued fraction expansion and let $p_n/q_n$ be the $n$th convergent.
   Letting $\norm{x}=\text{dist}(x,\Z)$, recall the following basic continued-fraction facts (see e.g. \cite{khinchin1964}):
    \[
    \max\left\{\frac{1}{2q_{n+1}},a_{n+2}\norm{q_{n+1}\alpha}\right\}<\|q_n\alpha\|=q_n\left|\alpha-\frac{p_n}{q_n}\right|<\frac{1}{q_{n+1}},
    \quad q_{n+1}=a_{n+1}q_n+q_{n-1}.
    \]
    \setcounter{claim}{0}

Let $\psi(n)=\log n\log\log n\log\log \log n$. By Khinchin's theorem \cite{khinchin1964}, for almost-every $\alpha\in(0,1)$ there exists $N\geq 1$ and an increasing sequence $n_k$ such that
    \be\label{eq:dioph1}
    q_{n+1}\leq q_n(\log q_n)^2, \quad \forall n\geq N,
    \ee
        \be\label{eq:dioph2}
    a_{n_k+1}\geq \psi(q_{n_k}), \quad \forall k\geq 1.
    \ee

We will show that certain $q_n-$orbits contribute a factor $\log q_n$ to the Birkhoff sums of $h$, and we will use the resonance at times $q_{n_k}$, given by \eqref{eq:dioph2} to improve the lower bound to $\log q_{n_k}(q_{n_k}\psi(q_{n_k}))^{-1}$.
Then, we will use the Diophantie condition \eqref{eq:dioph1} to show that the waiting time before encountering a good piece of orbit is not too large.

Fix $\alpha$ satisfying \eqref{eq:dioph1}--\eqref{eq:dioph2} and $x\in\T$ not in the $R_\alpha$-orbit of $0$. For any $I\subset\T$, let $|I|$ denote its Lebesgue measure. First, we prove the following lemma.
\begin{lemma}\label{lemma:logEsum}
    For every $n$ large enough there exists a segment $I_n\subset \T$ with $|I_n|=\norm{q_{n-1}\alpha}/4$ such that, if $y\in I_n$, then 
    \[
    |S(\alpha,h,y,q_n)|\geq 2^{-6}\log q_n.
    \]
\end{lemma}
\begin{proof}
Note $S(\alpha,h,\cdot,q_n)\in L^1(\T)$, so let $1/16>\delta>0$ be such that
\[\left|\int_I S(\alpha,h,y,q_n)\,dy\right|\leq 1,\quad \forall I\subset\T:|I|\leq\delta.\]
Let \[G_n:=\T\setminus\bigcup_{k=0}^{q_n-1}\left(-\frac{\delta}{q_n},\frac{\delta}{q_n}\right),\] then $|G_n|\geq \delta$, so $|\int_{G_n}S(\alpha,h,y,q_n)\,dy|\leq1$. In particular, there exists $y_n\in \T$ such that $|S(\alpha,h,y_n,q_n)|\leq 1$. Moreover, $S(\alpha,h,y_n,q_n)$ is smooth on every interval of $G_n$.

Let $(a_n,b_n)$ be the maximal segment of $G_n$ containing $y_n$. Since $G_n$ is constructed by removing a neighbourhood of a $q_n-$orbit, and all the points on such orbit are at least $\norm{q_{n-1}\alpha}$ appart from each other, we have that
\[ 
b_n-a_n\geq\norm{q_{n-1}\alpha}-\frac{2\delta}{q_n}\geq\norm{q_{n-1}\alpha}(1-4\delta)\geq\frac{3}{4} \norm{q_{n-1}\alpha}.
\]

Let $h_n=h\chi_{\norm{y}>\delta q_n^{-1}}$ and note 
\[
\int_\T h_n'=-\log q_n-\log \delta^{-1}-\log(1-\delta q_n^{-1})\leq-\frac{1}{2}\log q_n,\quad 
\text{Var}(h_n')=O(q_n),
\]
for $n$ large enough. For every $y\in(a_n,b_n)$, $S(\alpha,h,y,q_n)=S(\alpha,h_n,y,q_n)$, so Denjoy-Koksma yields
\[
     |S(\alpha,h',y,q_n)| \geq|q_n\int_{\T} h_n'|-|S(\alpha,h_n',y,q_n)-q_n\int_{\T} h_n'|\geq \frac{1}{2}q_n\log q_n-O(q_n)\geq\frac{1}{4}q_n\log q_n,
\]
for $n$ large enough.

Consider $(y_n+\norm{q_{n-1}\alpha}/8,y_n+3\norm{q_{n-1}\alpha}/8)$ and $(y_n-3\norm{q_{n-1}\alpha}/8,y_n-\norm{q_{n-1}\alpha}/8)$. At least one of them is fully contained in $(a_n,b_n)$; pick such an interval and call it $I_n$. Then for each point $y\in I_n$ there exists $\theta_n\in(a_n,b_n)$ such that
\[
|S(\alpha,h,y,q_n)|=|S(\alpha,h',\theta_n,q_n)||y-y_n|-1\geq\frac{1}{4}q_n\log q_n\frac{1}{8}\norm{q_{n-1}\alpha}-1\geq\frac{1}{64}\log q_n.
\]
\end{proof}

\begin{proof}[Proof of \eqref{eq:large} in Lemma \ref{claim:boundsEsums}]  
Let $N_0\geq N$, for $N$ as in \eqref{eq:dioph1}, such that Lemma \ref{lemma:logEsum} holds for $n\geq N_0$. After passing to a subsequence, we may assume $N_{k-1}\leq q_{n_k}$ for each $k\geq 1$. Fix $k$ and set $(a_k,b_k):=I_{n_k}$ as in Lemma \ref{lemma:logEsum}.

Fix $m_k:=\lfloor \psi(q_{n_k})/4\rfloor$ and note that, by \eqref{eq:dioph2},
\[
m_k\leq a_{n_k+1}/4<\frac{\norm{q_{n_k-1}\alpha}}{4\norm{q_{n_k}\alpha}}.
\]
Hence, letting
\[J_k:=\begin{cases}
(a_k, a_k+\norm{q_{n_k}\alpha}),&q_{n_k}\alpha>p_{n_k},
\\
    (b_k-\norm{q_{n_k}\alpha},b_k),& q_{n_k}\alpha<p_{n_k},
\end{cases}
\]
we have that, for any $0\leq n<m_k$, $R_\alpha^{nq_{n_k}}J_k=J_k+n{q_{n_k}\alpha}\subset I_{n_k}$.

By the rigidity of rotations and \eqref{eq:dioph1}, there exists $C'>0$ such that for any $y\in\T$,
\[
N(y,I_k):=\min\{n> N_{k-1}:R_\alpha^ny\in J_k\}\leq N_{k-1}+C'q_{n_k+1}\leq 2C'q_{n_k}(\log q_{n_k})^2.
\]
Hence, if 
\[|S(\alpha,h,x,N(x,I_k))|\geq{2^{-9}}\log N(x,I_k)\psi(N(x,I_k)),\] 
we can take $N_k:=N(x,I_k)$.
Otherwise, set $N_k:=N(x,I_k)+q_{n_k}m_k-1$ and note that, by Lemma \ref{lemma:logEsum},
\[\begin{aligned}
|S(\alpha,h,x,N_k)|&\geq -|S(\alpha,h,x,N(x,I_k))|+\sum_{m=0}^{m_k-1}|S(\alpha,h,R_\alpha^{N(x,I_k)+mq_{n_k}}x,q_{n_k})|
\\
&\geq- 2^{-9}\log N(x,I_k)\psi(N(x,I_k))+m_k 2^{-6}\log q_{n_k}\geq 2^{-10}\log(N_k)\psi(N_k).
\end{aligned}
\]
\end{proof}

\section*{Acknowledgments} The author would like to thank Dima Dolgopyat for bringing attention to this problem, for many helpful insights, and for his suggestions during the preparation of the paper. 

\vspace{10pt}
\noindent \textbf{Disclosure:} All the ideas and proofs in this manuscript were produced without the use of AI. AI was used as a linguistic tool to polish the text. 
\printbibliography
\end{document}